\documentclass[11pt]{article}

\usepackage{amsmath,amssymb,latexsym, amsfonts, amscd, amsthm}
\usepackage[small,nohug,heads=littlevee]{diagrams}
\diagramstyle[labelstyle=\scriptstyle]

\theoremstyle{plain}

\newtheorem{theorem}{Theorem}[section]

\newtheorem{proposition}[theorem]{Proposition}
\newtheorem{corollary}[theorem]{Corollary}

\newtheorem{lemma}[theorem]{Lemma}

\theoremstyle{definition}
\newtheorem{definition}[theorem]{Definition}

\newcommand{\mat}[4]{\left( \begin{array}{cc} {#1} & {#2} \\ {#3} & {#4}
\end{array} \right)}

\def\bdf{\begin{defn}}
\def\edf{\end{defn}}

\usepackage{color}

\begin{document}
\title{A New Construction of the Tame Local Langlands Correspondence for $GL(\ell,F)$, $\ell$ a prime}

\maketitle

\
\begin{center}

\Large{Moshe Adrian}

\end{center}
\

\

\

\

\

\

\

\

\

\

\

\

\

\

\

\

\

\

\

\

\

\

\

\

\begin{center}
Moshe Adrian

University of Utah

Department of Mathematics

155 S 1400 E ROOM 233

Salt Lake City, UT 84112-0090

\

madrian@math.utah.edu
\end{center}

\pagebreak

\tableofcontents

\pagebreak

\begin{abstract}
In this paper, we give a new construction of the tame local Langlands correspondence for $GL(n,F)$, $n$ a prime, where $F$ is a $p$-adic field.   In the tame case, supercuspidal representations of $GL(n,F)$ are parameterized by characters of elliptic tori, but the local Langlands correspondence is unnatural because it involves a twist by some character of the torus. Taking the cue from real groups, supercuspidal representations should instead be parameterized by characters of covers of tori. Over the reals, Harish-Chandra described the characters of discrete series restricted to compact tori. They are naturally written in terms of functions on a double cover of real tori. We write down a natural analogue of Harish-Chandra's character for $GL(n,F)$, and show that it is the character of a unique supercuspidal representation, away from the local character expansion. This paves the way for a natural construction of the local Langlands correspondence for $GL(n,F)$.
\end{abstract}

\section{Introduction}

In this paper, we reexamine the local Langlands correspondence for $GL(\ell, F)$, where $\ell$ is prime and $F$ is a non-Archimedean local field of characteristic zero, using character theory and ideas from the theory of real reductive groups.  Our main result is a construction of the tame local Langlands correspondence which circumvents some of the difficulties of \cite{bushnellhenniart}.  In particular, there are certain technical choices in the construction of the local Langlands correspondence which are explained from our point of view.  As a result, the construction of the local Langlands correspondence can be made to appear more natural.

Our results illuminate some new ideas about character theory of $p$-adic groups and local Langlands for $p$-adic groups not known before.  In particular, irreducible Weil group representations $W_F \rightarrow GL(\ell, \mathbb{C})$ and supercuspidal representations of $GL(\ell,F)$ are naturally parameterized not by certain characters of elliptic tori known as \emph{admissible pairs}, but by genuine characters of double covers of elliptic tori, as is the case for admissible representations of real groups.  We show that the supercuspidal representations of $GL(\ell,F)$ are naturally parameterized by genuine characters of double covers of elliptic tori using character theory.  To do this we rewrite supercuspidal characters in terms of double covers of elliptic tori as in Harish-Chandra's discrete series character formula and as in the Weyl character formula.  Rewriting the supercuspidal character formulas in this way paves the way for a natural construction of local Langlands for $GL(\ell,F)$.  In particular, it eliminates the need for any finite order character twists in the local Langlands correspondence for $GL(\ell,F)$ that arise in \cite{bushnellhenniart} and \cite{moy}.  As we shall see, our results and formulas also give justification and reason to the character formulas that first appeared in \cite{sallyshalika}, which may look like they came out of nowhere.

Let us recall the classical construction of the tame local Langlands correspondence for $GL(n,F)$.  In the tame case, Howe constructs a map (see \cite{howe}) $$ \{isomorphism \ classes \ of \ admissible \ pairs \} \rightarrow \{supercuspidal \ representations \ of \ GL(n,F) \}$$ $$\hspace{-.5in} (E/F, \chi) \mapsto \pi_{\chi}$$ where $\chi$ is a character of $E^*$ satisfying certain conditions to be described later, and $E/F$ is an extension of degree $n$.  This map is a bijection (see \cite{moy}).  Moreover, we have a bijection $$\{admissible \ pairs \ (E/F, \chi) \} \rightarrow \{irreducible \ W_F \rightarrow GL(n,\mathbb{C}) \}$$ $$(E/F, \chi) \mapsto Ind_{W_E}^{W_F}(\chi) =: \varphi(\chi)$$ (see \cite{moy}).  The problem is that the obvious map, $$\varphi(\chi) \mapsto \pi_{\chi},$$ the so-called ``naive correspondence'', is not the local Langlands correspondence because $\pi_{\chi}$ has the wrong central character.  Instead, the local Langlands correspondence is given by $$\varphi(\chi) \mapsto \pi_{\chi \Delta_{\chi}}$$ for some subtle finite order character $\Delta_{\chi}$ of $E^*$.  The presence of the character twist $\Delta_{\chi}$ makes the correspondence look unnatural.  We will show that if one considers genuine characters of a canonical double cover of elliptic tori rather than characters of elliptic tori, then one obtains a natural local Langlands correspondence.  We do this in the following way.

Taking the cue from the theory of real groups, we use genuine characters $\tilde\chi$ of certain double covers of elliptic tori, denoted $T(F)_{\tau \circ \rho}$, instead of characters of elliptic tori $T(F)$, to parameterize both representations of $W_F$ and supercuspidal representations of $GL(\ell,F)$ using character theory.  We give a method for attaching a genuine character of a double cover of elliptic tori satisfying certain regularity conditions, to a supercuspidal Weil parameter of $GL(\ell,F)$:

\begin{equation}
\{irreducible \ W_F \rightarrow GL(\ell, \mathbb{C}) \} \leftrightarrow \{regular \ genuine \ characters \ of \ T(F)_{\tau \circ \rho} \} \label{eq:parameters1}
\end{equation}
$$\hspace{-.75in}\varphi \mapsto \tilde\chi$$ Moreover, as we shall show, supercuspidal characters of $GL(\ell, F)$ correspond naturally to regular genuine characters of $T(F)_{\tau \circ \rho}$ rather than admissible pairs $(E/F, \chi)$.  Given a regular genuine character $\tilde\chi$ of $T(F)_{\tau \circ \rho}$, we write down a conjectural Harish-Chandra type character formula, denoted $F(\tilde\chi)$.  We show that this naturally gives a bijection

\begin{equation}
\{regular \ genuine \ characters \ of \ T(F)_{\tau \circ \rho} \} \leftrightarrow \{supercuspidal \ representations \ of \ GL(\ell,F) \} \label{eq:parameters2}
\end{equation}
$$\tilde\chi \mapsto \pi(\tilde\chi)$$ where $\pi(\tilde\chi)$ is the unique supercuspidal representation of $GL(\ell,F)$, whose character, restricted to a certain natural subset of $T(F)$ (to be described later), is $F(\tilde\chi)$.

Then the composition of bijections (\ref{eq:parameters1}) and (\ref{eq:parameters2}), $$\varphi \mapsto \tilde\chi \mapsto \pi(\tilde\chi),$$ is the local Langlands correspondence for $GL(\ell,F)$

Let us explain why double covers of tori play a role.  We start by considering the group $PGL(2,F)$.  Recall that the representations of $PGL(2,F)$ are precisely the representations of $GL(2,F)$ with trivial central character.  One of the conditions of the local Langlands correspondence for $GL(n,F)$ says that if $\varphi :W_F \rightarrow GL(n,\mathbb{C})$ is irreducible, then $det(\varphi) = \omega_{\pi(\varphi)}$, where $\omega_{\pi(\varphi)}$ denotes the central character of $\pi(\varphi)$, and where $\pi(\varphi)$ denotes the supercuspidal representation of $GL(n,F)$ that corresponds to $\varphi$ under the local Langlands correspondence (see \cite[Chapter 34]{bushnellhenniart}).  Let $\varphi$ be a supercuspidal Weil parameter for $PGL(2,F)$ (that is, an irreducible representation $W_F \rightarrow GL(2,\mathbb{C})$ that parameterizes a supercuspidal representation of $GL(2,F)$ with trivial central character).  Then $\varphi = Ind_{W_E}^{W_F}(\chi)$, for some admissible pair $(E/F, \chi)$.  The condition $det(\varphi) = \omega_{\pi(\varphi)}$ implies that $\chi|_{F^*} \otimes \aleph_{E/F} = 1$ (see \cite[Chapter 34]{bushnellhenniart}) where $\aleph_{E/F}$ denotes the local class field theory character of $F^*$ relative to the extension $E/F$, so $\chi|_{F^*} = \aleph_{E/F}$. Therefore, the supercuspidal Weil parameters of $PGL(2,F)$ are parameterized by the admissible pairs $(E/F, \chi)$ where $\chi|_{F^*} = \aleph_{E/F}$.  Now, recall again that supercuspidal Weil parameters of $GL(2,F)$ are parameterized by characters of elliptic tori.  One might ask whether the supercuspidal Weil Parameters of $PGL(2,F)$ are parameterized by characters of its elliptic tori $E^* / F^*$.  However, we have just seen that they are parameterized by characters $\chi$ of $E^*$ whose restriction to $F^*$ is $\aleph_{E/F}$.  Such a $\chi$ is not a character of the elliptic torus $E^* / F^*$ in $PGL(2,F)$.  Rather, it is a genuine character of the double cover $E^* / N(E^*)$ arising from the canonical exact sequence (note that $ker(\aleph_{E/F}) = N(E^*)$) $$1 \longrightarrow F^* / N(E^*) \longrightarrow E^* / N(E^*) \rightarrow E^* / F^* \longrightarrow 1$$ where $N$ denotes the norm map from $E$ to $F$.  Since $F^* / N(E^*) \cong \mathbb{Z}/ 2 \mathbb{Z}$ by Local Class Field Theory, we have that $E^* / N(E^*)$ is a double cover of $E^* / F^*$.  Moreover, $\chi$ is not trivial on all of $F^*$, so it doesn't factor to a character of $E^* / F^*$.  This means that $\chi$ factors to a genuine character, denoted $\tilde\chi$, of $E^* / N(E^*)$.  Therefore, we are getting naturally that the supercuspidal Weil parameters of $PGL(2,F)$ are paramterized by genuine characters of a double cover of the elliptic torus $E^* / F^*$ in $PGL(2,F)$.  In fact, this double cover $E^* / N(E^*)$ is none other an the analogue of the $\rho$-cover that appears in the theory over the reals, which is a natural double cover of a real torus.  The same above reasoning applies to the setting of $PGL(\ell,F)$, where $\ell$ is an odd prime.  We will show that double covers arise naturally as well for $GL(2,F)$ and $GL(\ell,F)$, where $\ell$ is an odd prime.

In the theory of real groups, admissible homomorphisms $W_{\mathbb{R}} \rightarrow {}^L G$ naturally produce genuine characters of the \emph{$\rho$-cover of $T(\mathbb{R})$}, denoted $T(\mathbb{R})_{\rho}$, a certain double cover of $T(\mathbb{R})$, which we now define.  First, we need to make the following definition.

\begin{definition}\label{rhocoverintro}
Let $A$, $B$, and $C$ be groups, and suppose we have homomorphisms $\phi_1 : A \rightarrow C$, $\phi_2 : B \rightarrow C$.  Then the \emph{pullback} of these two homomorphisms is the group $$A \times_C B := \{ (a,b) \in A \times B \ | \ \phi_1(a) = \phi_2(b) \}$$ together with projections $$\pi_1 : A \times_C B \rightarrow A \ \ \ \pi_2 : A \times_C B \rightarrow B$$ $$(a,b) \mapsto a \ \ \ \ (a,b) \mapsto b$$  Then the following diagram commutes:

$$
\begin{CD}
A \times_C B @> \pi_1 >> A\\
@VV\pi_2V @VV\phi_1V\\
B @>\phi_2>> C
\end{CD}
$$

\end{definition}

\begin{definition}\label{rhocoverdefinition}
Let $G$ be a connected reductive group over $\mathbb{R}$, and let $T \subset G$ a maximal torus over $\mathbb{R}$. Let $X^*(T)$ be the character group of $T$.  Let $\Delta^+$ be a set of positive roots of $G$ with respect to $T$.  Let $\rho = \frac{1}{2} \displaystyle\sum_{\alpha \in \Delta^+} \alpha$.  Then $2 \rho \in X^*(T)$.  Viewing $2 \rho$ as also a character of $T(\mathbb{R})$ by restriction, we define the \emph{$\rho$-cover of $T(\mathbb{R})$} as the pullback of the two homomorphisms $$2 \rho : T(\mathbb{R}) \rightarrow \mathbb{C}^* \ \ \ \ \ \ \Upsilon_{\mathbb{R}} : \mathbb{C}^* \rightarrow \mathbb{C}^*$$ $$t \mapsto 2 \rho(t) \ \ \ \ \ \ \ \ \ \ \ z \mapsto z^2$$  We denote the $\rho$-cover by $T(\mathbb{R})_{\rho}$, and so the following diagram commutes:

$$
\begin{CD}
T(\mathbb{R})_{\rho} @> \rho >> \mathbb{C}^*\\
@VV \Pi_{\mathbb{R}} V @VV \Upsilon_{\mathbb{R}} V\\
T(\mathbb{R}) @>2 \rho>> \mathbb{C}^*
\end{CD}
$$
\end{definition}

Note that because of the commutativity of the diagram, although $\rho$ is not necessarily a character of $T(\mathbb{R})$,  $\rho$ is naturally a character of $T(\mathbb{R})_{\rho}$.  Moreover, $\Pi_{\mathbb{R}}$ is the canonical projection $\Pi_{\mathbb{R}}(t, \lambda) = t$.

The genuine characters of $T(\mathbb{R})_{\rho}$ that naturally arise from Weil parameters are used to form $L$-packets. In the case of $GL(n,\mathbb{R})$, $L$-packets are singletons, and we have that the irreducible admissible representations of $GL(n,\mathbb{R})$ and admissible homomorphisms $W_{\mathbb{R}} \rightarrow GL(n, \mathbb{C})$ are in natural bijection with genuine characters $\tilde\chi$ of $T(\mathbb{R})_{\rho}$.  The composition of these two parameterizations is in fact the local Langlands correspondence for $GL(n,\mathbb{R})$.

One can write down more explicitly the correspondence for relative discrete series representations, and this is the motivation for our work.  To state the theorem, we make a few preliminary remarks.

\begin{definition}\label{delta0}
Let $G$ be a connected reductive group over $\mathbb{R}$, $T \subset G$ a maximal torus over $\mathbb{R}$.  Let $\Delta^+$ be a set of positive roots of $G$ with respect to $T$.  Define $$\Delta^0(h, \Delta^+) := \prod_{\alpha \in \Delta^+} (1 - \alpha^{-1}(h)), \ h \in T(\mathbb{R})$$
\end{definition}

Recall that in general, $\rho$ is not in $X^*(T)$.  Therefore, $\rho(h)$ does not make sense if $h \in T(\mathbb{R})$.  However, $\rho$ is a well-defined character of $T(\mathbb{R})_{\rho}$.  If $\tilde h \in T(\mathbb{R})_{\rho}$ is any element such that $\Pi_{\mathbb{R}}(\tilde h) = h$, we may consider the function $\Delta^0(h, \Delta^+) \rho(\tilde h)$.  This function lives on $T(\mathbb{R})_{\rho}$, and we have the following theorem.

\begin{theorem}\label{Harish-Chandra}(Harish-Chandra)
Let G be a connected reductive group, defined over $\mathbb{R}$.  Suppose that G contains a Cartan subgroup T  that is defined over $\mathbb{R}$ and that is compact mod center.  Let $\Delta^+$ be a set of positive roots of G with respect to T.  Let $\rho := \frac{1}{2} \displaystyle\sum_{\alpha \in \Delta^+} \alpha$.  Let $\tilde\chi$ be a genuine character of $T(\mathbb{R})_{\rho}$ that is regular.  Let $W := W(G(\mathbb{R}),T(\mathbb{R}))$ be the relative Weyl group of $G(\mathbb{R})$ with respect to $T(\mathbb{R})$.  Let $\epsilon(s) := (-1)^{\ell(s)}$ where $\ell(s)$ is the length of the Weyl group element $s \in W$. Let $T(\mathbb{R})^{reg}$ denote the regular set of $T(\mathbb{R})$.  Then there exists a unique constant $\epsilon(\tilde\chi, \Delta^+) = \pm 1$, depending only on $\tilde\chi$ and $\Delta^+$, and a unique relative discrete series representation of $G(\mathbb{R})$, denoted $\pi(\tilde\chi)$, such that
$$\theta_{\pi(\tilde\chi)}(h) = \epsilon(\tilde\chi, \Delta^+) \frac{\displaystyle\sum_{s \in W} \epsilon(s) \tilde\chi({}^s \tilde h)}{\Delta^0(h, \Delta^+) \rho(\tilde h)}, \ h \in T(\mathbb{R})^{reg}$$ where $\tilde h \in T(\mathbb{R})_{\rho}$ is any element such that $\Pi_{\mathbb{R}}(\tilde h) = h$.  Moreover, every relative discrete series character of $G(\mathbb{R})$ is of this form.
\end{theorem}

This character formula is a variant of the formula found in the literature, using $\rho$-covers.  We can be more specific about the constant $\epsilon(\tilde\chi, \Delta^+)$.  In particular, $\epsilon(\tilde\chi, \Delta^+) = (-1)^{\ell(s)}$ where $s \in W$ makes $d \tilde\chi$ dominant for $\Delta^+$, $d \tilde\chi$ denoting the differential of $\tilde\chi$.

For $GL(2,\mathbb{R})$, the local Langlands correspondence for relative discrete series representations is as follows.  Fix a positive set of roots $\Delta^+$ of $G$ with respect to $T$.  Let $\varphi : W_{\mathbb{R}} \rightarrow GL(2, \mathbb{C})$ be a relative discrete series parameter, and let $T(\mathbb{R})$ be the compact mod center torus of $GL(2,\mathbb{R})$.  Then $\varphi$ naturally gives rise to a genuine character $\tilde\chi$ of $T(\mathbb{R})_{\rho}$.  By Harish-Chandra's discrete series theorem, $\tilde\chi$ gives rise to a unique relative discrete series representation, denoted $\pi(\tilde\chi)$, whose character, restricted to the regular elements of $T(\mathbb{R})$, is $$F(\tilde\chi)(h) := \epsilon(\tilde\chi, \Delta^+) \frac{\displaystyle\sum_{s \in W} \epsilon(s) \tilde\chi({}^s \tilde h)}{\Delta^0(h, \Delta^+) \rho(\tilde h)}, \ h \in T(\mathbb{R})^{reg}$$ where $\tilde h \in T(\mathbb{R})_{\rho}$ is any element such that $\Pi_{\mathbb{R}}(\tilde h) = h$.  The map

\begin{equation}
\varphi \mapsto \pi(\tilde\chi)  \label{eq:naivereal}
\end{equation}

\noindent is the local Langlands correspondence for relative discrete series representations of $GL(2,\mathbb{R})$.
Thus, one can write down the correspondence for relative discrete series in terms of character theory.  This is the approach we take in this paper, and we will show that the above correspondence (\ref{eq:naivereal}) carries over naturally to the $p$-adic setting.

One of the results that we will prove is an analogue of Harish-Chandra's theorem for $GL(\ell,F)$, where $F$ is a $p$-adic field of characteristic zero.  In doing this, we give a new realization of the tame local Langlands correspondence for $GL(\ell,F)$, and the character twists $\Delta_{\chi}$ go away.  Before we present our main results, we need to define the covers of tori that will be essential, which are an analogue of the $\rho$-cover that appears in the theory over the reals.

Let $G$ be a connected reductive group defined over $F$, and $T$ a maximal torus in $G$ defined over $F$.  Let $\Delta^+$ be a choice of positive roots of $G$ with respect to $T$.  Let $\rho = \frac{1}{2} \displaystyle\sum_{\alpha \in \Delta^+} \alpha$.  Let $\lambda$ be a character of $K^*$, where $K/F$ contains the minimal splitting field of $T$. Note that the image of $2 \rho$, restricted to $T(F)$, lies in $K^*$.

\begin{definition}\label{tauofrhocover}
We define the \emph{$\lambda \circ \rho$-cover of $T(F)$}, denoted $T(F)_{\lambda \circ \rho}$, as the pullback of the two homomorphisms $$\lambda \circ 2 \rho : T(F) \rightarrow \mathbb{C}^* \ \ \ \ \ \ \Upsilon : \mathbb{C}^* \rightarrow \mathbb{C}^*$$ $$\ \ \ \ \ t \mapsto \lambda \circ 2 \rho(t) \ \ \ \ \ \ \ z \mapsto z^2$$

$$
\begin{CD}
T(F)_{\lambda \circ \rho} @>\lambda \circ \rho>> \mathbb{C}^*\\
@VV \Pi V @VV \Upsilon V\\
T(F) @>\lambda \circ 2 \rho>> \mathbb{C}^*
\end{CD}
$$

\noindent That is, $T(F)_{\lambda \circ \rho} = \{(z,w) \in T(F) \times \mathbb{C}^* : \lambda(2 \rho(z)) = w^2 \}$
\end{definition}

Note that the above map $T(F)_{\lambda \circ \rho} \rightarrow \mathbb{C}^*$ sends $(z, w)$ to $w$.  We have denoted this map by $\lambda \circ \rho$, even though this map is not literally $\lambda$ composed with $\rho$.  Moreover, $\Pi$ is the canonical projection $\Pi(z, w) = z$.  We will use these maps repeatedly.

Our main results will be the following theorems.

\begin{theorem}\label{gl2theorem1}
Let $\ell$ be prime.  If $\ell$ is $2$, we assume that the residual characteristic of $F$ is not $2$.  If $\ell \neq 2$, we assume that the residual characteristic of $F$ is greater than $2 \ell$.  Let $G(F) = GL(\ell,F)$, and let $T(F) = E^*$ be an elliptic torus in $GL(\ell,F)$, so $E = F(\sqrt[\ell]{\Delta})$ for some $\Delta \in F^*$.  Let $L$ be the unique unramified extension of $F$ of degree $\ell - 1$.  Let $\tau_o$ be any character of $(EL)^*$ whose restriction to $L^*$ is $\aleph_{EL/L}$, where $\aleph_{EL/L}$ is the local class field theory character of $L^*$ relative to $EL/L$.  Let $\tau := \tau_o \ | \ |_{EL}$ where $| \ |_{EL}$ denotes the $EL$-adic absolute value.  Let $\Delta^+$ be a set of positive roots of G with respect to T.  Let $\rho := \frac{1}{2} \displaystyle\sum_{\alpha \in \Delta^+} \alpha$.  Let $T(F)_{\tau \circ \rho}$ be the $\tau \circ \rho$ cover of $T(F)$.  Let $\tilde\chi$ be a genuine character of $T(F)_{\tau \circ \rho}$ that is regular.  Let $W = W(G(F),T(F))$ denote the relative Weyl group of $G(F)$ with respect to $T(F)$.  If $s \in W(G(F),T(F))$, let $\epsilon(s) := (-1, \Delta)^{\ell(s) (\ell+1)}$, where $(,)$ denotes the Hilbert symbol of $F$ and $\ell(s)$ denotes the length of $s$.  Let $T(F)^{reg}$ denote the regular elements of $T(F)$.

Then there exists a unique constant $\epsilon(\tilde\chi, \Delta^+, \tau)$, depending only on $\tilde\chi, \Delta^+$, and $\tau$, and a unique supercuspidal representation of $GL(\ell,F)$ denoted $\pi(\tilde\chi)$, such that
$$\theta_{\pi(\tilde\chi)}(z) = \epsilon(\tilde\chi, \Delta^+, \tau) \frac{\displaystyle\sum_{s \in W} \epsilon(s)\tilde\chi({}^s w)}{\tau(\Delta^0(z,\Delta^+)) (\tau \circ \rho)(w)}, \ \forall z \in T(F)^{reg} \ : 0 \leq n(z) \leq r/2$$ where $w \in T(F)_{\tau \circ \rho}$ is any element such that $\Pi(w) = z$ and $r$ is the depth of $\pi(\tilde\chi)$.  Moreover, every supercuspidal character of $GL(\ell,F)$ is of this form.
\end{theorem}

We will define all of the notation in the above theorem in Chapters 5-8, including $n(z)$, $\epsilon(\tilde\chi, \Delta^+, \tau)$, and regularity.  We remark that $n(z)$ comes from a canonical filtration on the torus $T(F)$, and is defined in \cite{debacker}.  Notice that when we treat the case of depth zero representations (i.e. $r = 0$), the range $\{z \in T(F)^{reg} : 0 \leq n(z) \leq r/2 \}$ becomes $\{z \in T(F)^{reg} : n(z) = 0 \}$.  We wish to make a few comments about the constant $\epsilon(\tilde\chi, \Delta^+, \tau)$.  Firstly, $\epsilon(\tilde\chi, \Delta^+, \tau)^4 \in \mathbb{R}^*$.  Moreover, $|\epsilon(\tilde\chi, \Delta^+, \tau)|$ is a known real number in that it has to do with a canonical measure on the Lie algebra.  The subtlety of $\epsilon(\tilde\chi, \Delta^+, \tau)$ is in the value of $\frac{\epsilon(\tilde\chi, \Delta^+, \tau)}{|\epsilon(\tilde\chi, \Delta^+, \tau)|} \in \{ \pm 1, \pm i \}$.

Now let $\varphi$ be a supercuspidal Weil parameter for $GL(\ell,F)$.  We will show in Chapters 5 and 8 how to construct a regular genuine character, $\tilde\chi$, of $T(F)_{\tau \circ \rho}$, from $\varphi$.  We will then prove the following theorem.

\begin{theorem}\label{gl2theorem2}
In the setting of Theorem \ref{gl2theorem1}, the assignment $$\varphi \mapsto \tilde\chi \mapsto \pi(\tilde\chi)$$ is the Local Langlands correspondence for $GL(\ell,F)$.
\end{theorem}

Let us be a bit more explicit about the representation $\pi(\tilde\chi)$.  In particular, if $\varphi : W_F \rightarrow GL(\ell,\mathbb{C})$ is a supercuspidal Weil parameter for $GL(\ell,F)$, and $\varphi = Ind_{W_E}^{W_F}(\chi)$ for some admissible pair $(E/F, \chi)$, then the $\pi(\tilde\chi)$ that $\varphi$ maps to under the previous theorem is $\pi_{\chi \Delta_{\chi}}$.

In Theorems \ref{gl2theorem1} and \ref{gl2theorem2}, we have made our assumptions on $\ell$ for three reasons.  Firstly, because our methods require the knowledge of the supercuspidal character formulas for $GL(\ell,F)$, and these have only been completely computed so far in the cases where the residual characteristic of $F$ is greater than $\ell$ (see \cite{debacker}).  Secondly, it is unclear whether our methods will generalize to wildly ramified situations.  Thirdly, we will use some results from \cite{spice}, which assumes in the case of $\ell > 2$ that the residual characteristic of $F$ is greater than $2 \ell$.  Once the supercuspidal characters for $GL(n,F)$ become available for arbitrary $n$, we expect that our methods will generalize to the case where the residual characteristic of $F$ is coprime to $n$.

We now briefly present an outline of the paper.  In section 2, we introduce some notation that will be used throughout. In section 3, we recall the necessary theory from real groups that we need.  In particular, we describe some of the basic ingredients of the local Langlands correspondence for real reductive groups, following \cite{adamsvogan}.  In section 4, we recall the necessary background to describe the local Langlands correspondence  for $GL(2,F)$, following \cite{bushnellhenniart}.  In section 5, we introduce the relevant double covers that play a role in our theory.  We then rewrite the supercuspidal characters of $GL(2,F)$ in terms of regular genuine characters of double covers of elliptic tori, and show that the distribution characters are determined by the values on the range $\{ z \in T(F)^{reg} : 0 \leq n(z) \leq r/2 \}$. Finally, we present a natural construction of the tame local Langlands correspondence for positive depth supercuspidal representations of $GL(2,F)$ using the above theory.  In section 6, we treat the case of depth zero representations of $GL(2,F)$, and the theory is similar.  In section 7, we recall the necessary background to describe the local Langlands correspondence  for $GL(\ell,F)$ where $\ell$ is an odd prime, following \cite{moy}.  In section 8, we develop our general theory for $GL(\ell,F)$, which carries over directly from the theory we developed for $GL(2,F)$ in section 5.  In section 9, we treat the case of depth zero representations of $GL(\ell,F)$, and the theory is analogous.  In section 10, we list some useful facts that we need periodically in the paper that would have otherwise been distracting to the flow of the paper.

\emph{Acknowledgements}. I would like to thank Jeffrey Adams, my advisor, for directing my attention to this question and for being a tremendous support throughout my career; Loren Spice for many very helpful conversations and for answering many questions about his work.  I would also like to thank Stephen DeBacker for a very helpful conversation during my visit to Michigan.  I would also like to thank Jeffrey Adler for several very helpful conversations about my work.  I also had several useful conversations with Hunter Brooks and Sean Rostami.

\section{Notation and Definitions}

Let $F$ denote a local field of characteristic zero, $\mathfrak{o}_F$ its ring of integers, and $\mathfrak{p}_F$ the maximal ideal of $\mathfrak{o}_F$. We let $p$ denote a uniformizer of $F$.  Let $k_F$ denote the residue field of $F$ with cardinality $q$. We choose an element $\Phi \in Gal(\overline{F}/F)$ whose inverse induces on $\overline{k_F}$ the map $x \mapsto x^q$. Throughout, we fix once and for all a nontrivial additive character $\psi$ of $F$ of level one.  If $E/F$ is a separable extension, $N$ will denote the norm map from $E$ to $F$, $Tr_{E/F}$ will denote the trace map from $E$ to $F$, and $Aut(E/F)$ will denote the group of automorphisms of $E$ that fix $F$ pointwise.  When we write a decomposition $w = p^n u$ where $w \in F^*$, we mean that $u \in \mathfrak{o}_F^*$.  If $E/F$ is quadratic and $E = F(\delta)$, we will frequently decompose an element $w \in E$ as $w = p^n u + p^m v \delta$ where we are viewing $E$ as a vector space over $F$ with basis $1, \delta$, and $u,v \in \mathfrak{o}_F^*$.  If $E/F$ is quadratic, then we will write $\overline{w}$ instead of $\upsilon(w)$ where $1 \neq \upsilon \in Gal(E/F)$. Let $( , )_F$ denote the Hilbert symbol of $F$; most of the time we will write $( , )$ when there is no confusion about the field.  We also set $U_F^n := 1 + \mathfrak{p}_F^n$ and  $U_F = \mathfrak{o}_F^*$.  If $E/F$ is Galois, we let $\aleph_{E/F}$ denote the local class field theory character of $F^*$ relative to the extension $E/F$.  If $K$ is a local non-archimedean field of characteristic zero, we let $| \ |_K$ denote the $K$-adic absolute value of $K$.  In Chapters 5-6, $\tau_o$ will denote any character of $E^*$ whose restriction to $F^*$ is $\aleph_{E/F}$, where $E/F$ is a tame quadratic extension, and we will set $\tau := \tau_o \ | \ |_{E}$ where $| \ |_E$ is the $E$-adic absolute value.  In Chapters 8-9, $\tau_o$ will denote any character of $(EL)^*$ whose restriction to $L^*$ is $\aleph_{EL/L}$, where $E/F$ is a tame degree $\ell$ extension and $L/F$ is the degree $\ell - 1$ unramified extension, and we will set $\tau := \tau_o \ | \ |_{EL}$ where $| \ |_{EL}$ is the $EL$-adic absolute value.  We will generally write $| \ |$ when it is clear which field we are referring to.

If $$1 \rightarrow \mathbb{Z} / 2 \mathbb{Z} \rightarrow A \rightarrow B \rightarrow 1$$ is an exact sequence of groups, then a character $\beta$ of $A$ is said to be \emph{genuine} if $\beta|_{\mathbb{Z} / 2 \mathbb{Z}}$ is not trivial (that is, $\beta$ does not arise from a character $A \rightarrow B$).  When we say that a $2$-fold cover of a group (as above) \emph{splits}, we mean that the exact sequence splits.  If $G$ denotes any group, then $G^{ab}$ denotes its abelianization.  If $G$ is a connected reductive group defined over $F$ and $T$ is a maximal torus in $G$ defined over $F$, then we will frequently write the relative Weyl group as $W$ instead of $W(G(F),T(F))$.  We will write $T(F)^{reg}$ for the set of regular elements in $T(F)$.  If $B$ is a normal subgroup of $A$ and $a \in A$, then we will write $[a]$ to denote the class of $a$ in $A/B$.

If $\pi$ is a representation of a group, let $\omega_{\pi}$ denote its central character.

\section{Background from real groups}\label{realgroupschapter}

In order to motivate the theory that we wish to develop for $p$-adic groups, we describe the corresponding theory over $\mathbb{R}$ since this is what our theory is based upon.  We will briefly recall the relevant theory of the local Langlands correspondence over $\mathbb{R}$.  More information can be found in \cite{adamsvogan}.

\subsection{Covers of Tori}

It will be important to describe a part of the local Langlands correspondence having to do with discrete series representations.  Recall Definition (\ref{rhocoverdefinition}). The Weyl group acts on $T(\mathbb{R})_{\rho}$ as follows: If $(t, \lambda) \in T(\mathbb{R})_{\rho}$, then define

\begin{equation}
s(t, \lambda) := (st, e^{s^{-1} \rho - \rho}(t) \lambda)  \ \ \forall s \in W(G(\mathbb{R}), T(\mathbb{R})) \ \ \label{eq:weylgroupreal}
\end{equation}

\begin{definition}
A genuine character $\tilde\chi$ of $T(\mathbb{R})_{\rho}$ is called \emph{regular} if ${}^s \tilde\chi \neq \tilde\chi \ \forall s \in W(G(\mathbb{R}),T(\mathbb{R}))$ where ${}^s \tilde\chi(t, \lambda) := \tilde\chi(s^{-1}(t, \lambda))$.
\end{definition}

\subsection{Discrete series Langlands paramaters and character formulas}

Let $G$ be a connected reductive group over $\mathbb{R}$ that contains a compact torus defined over $\mathbb{R}$.  It is known that this is equivalent to $G(\mathbb{R})$ having discrete series representations.

\begin{definition}
Let $t$ be an indeterminate and let $k$ denote the rank of $G$.  For $h \in G$, define $D_G(h)$ by
$$det(t + 1 - Ad(h)) = D_G(h)t^k + ... (terms \ of \ higher \ order)$$
\end{definition}

\noindent Then if $\Delta$ is the set of roots of $T$ in $G$, $$D_G(h) = \prod_{\alpha \in \Delta} (1 - \alpha(h)).$$

\noindent Let $\Delta^+$ be a set of positive roots.  Then if the cardinality of $\Delta^+$ is $n$, we have $$(-1)^n D_G(h) = \Delta^0(h, \Delta^+)^2 (2 \rho)(h)$$ where $\Delta^0(h, \Delta^+)$ is as in Definition \ref{delta0}.
Then, if we define $|\rho(h)| := |2 \rho(h)|^{\frac{1}{2}}$ we get that $$|D_G(h)|^{\frac{1}{2}} = |\Delta^0(h, \Delta^+)| |\rho(h)|$$
If $\tilde h \in T(\mathbb{R})_{\rho}$ maps to $h \in T(\mathbb{R})$ via the canonical projection, then $$|D_G(h)|^{\frac{1}{2}} = |\Delta^0(h, \Delta^+)| |\rho(h)| = |\Delta^0(h, \Delta^+)| |\rho(\tilde h)|$$
Let $G(\mathbb{R}) = GL(2,\mathbb{R})$ and $T(\mathbb{R})$ its compact maximal torus.  Fix a positive set of roots $\Delta^+$ of $G$ with respect to $T$.  Let $\varphi : W_{\mathbb{R}} \rightarrow GL(2, \mathbb{C})$ be a discrete series Weil parameter.  By the theory in \cite{adamsvogan}, $\varphi$ caonically gives rise to a genuine character $\tilde\chi$ of $T(\mathbb{R})_{\rho}$.  By Harish-Chandra's discrete series theorem, $\tilde\chi$ canonically gives rise to a unique discrete series representation, denoted $\pi(\tilde\chi)$, whose distribution character is $$\theta_{\pi(\tilde\chi)}(h) :=  \epsilon(\tilde\chi, \Delta^+) \frac{\displaystyle\sum_{s \in W} \epsilon(s) \tilde\chi({}^s \tilde h)}{\Delta^0(h, \Delta^+) \rho(\tilde h)}, \ h \in T(\mathbb{R})^{reg}$$ where $\tilde h \in T(\mathbb{R})_{\rho}$ is any element such that $\Pi_{\mathbb{R}}(\tilde h) = h$.  Then the map $\varphi \mapsto \pi(\tilde\chi)$ is the local Langlands correspondence for discrete series representations of $GL(2,\mathbb{R})$.  The rest of this paper will be devoted to proving the analogous result for $GL(\ell,F)$, where $F$ is a local non-Archimedean field of characteristic zero, and $\ell$ is prime.

\section{Existing Description of Local Langlands Correspondence for $GL(2,F)$}

In this chapter, we describe the construction of the local Langlands correspondence for $GL(2,F)$ as explained in \cite{bushnellhenniart}.

\subsection{Admissible Pairs}\label{admissiblepairs}

Let $E/F$ be a tamely ramified quadratic extension and $\chi$ a character of $E^*$.  Recall that $N$ denotes the norm map from $E$ to $F$.

\begin{definition}
The pair $(E/F, \chi)$ is called an \emph{admissible pair} if

(i) $\chi$ does not factor through $N$ and

(ii) If $\chi|_{1 + \mathfrak{p}_E}$ factors through $N$, then $E/F$ is unramified.

\end{definition}

We write $\mathbb{P}_2(F)$ for the set of $F$-isomorphism classes of admissible pairs.  For more information about admissible pairs, see \cite[Section 18]{bushnellhenniart}.

\begin{definition} Let $(E/F, \chi)$ be an admissible pair where $\chi$ is level $n$.  We say that $(E/F, \chi)$ is \emph{minimal} if $\chi|_{U_E^n}$ does not factor through $N$.
\end{definition}

\subsection{Depth zero supercuspidal representations of $GL(2,F)$}

Let $(E/F, \chi)$ be an admissible pair where $\chi$ has level 0.  By definition of admissible pair, this implies that $E/F$ is unramified.  Then $k_E / k_F$ is a quadratic extension.  Moreover, since $\chi|_{1 + \mathfrak{p}_E} = 1$, $\chi|_{U_E}$ is the inflation of a character, call it $\chi$ again, of $k_E^*$.  By the theory of finite groups of Lie type, the character $\chi$ then gives rise to an irreducible cuspidal representation $\lambda'$ of $GL(2, k_F)$.  Let $\lambda$ be the inflation of $\lambda'$ to $GL(2, \mathfrak{o}_F)$.  We may extend $\lambda$ to a representation $\Lambda$ of $K := F^* GL(2, \mathfrak{o}_F)$ by setting $\Lambda|_{F^*} = \chi|_{F^*}$, and induce the resulting representation to all of $G$.  Set $$\pi_{\chi} = cInd_K^G \Lambda$$ where $cInd$ denotes compact induction.

Let $\mathbb{P}_2(F)_0$ denote the set of admissible pairs of level zero, and let $\mathbb{A}_2^0(F)_0$ denote the set of equivalence classes of depth zero supercuspidal representations of $GL(2,F)$.

\begin{proposition}\label{depthzerogl2}
The map (E/F, $\chi$) $\mapsto \pi_{\chi}$ induces a bijection $$\mathbb{P}_2(F)_0 \rightarrow \mathbb{A}_2^0(F)_0$$  Furthermore, if (E/F, $\chi) \in \mathbb{P}_2(F)_0$, then:

(i) if $\phi$ is a character of $F^*$ of level zero, then $\pi_{\chi \phi_E} = \phi \pi_{\chi}$

(ii) if $\pi = \pi_{\chi}$, then $\omega_{\pi} = \chi|_{F^*}$

\end{proposition}

\proof

See \cite[Section 19.1]{bushnellhenniart}.

\qed

\subsection{Positive depth supercuspidal representations of $GL(2,F)$}\label{alpha}

Let $(E/F, \chi)$ be a minimal admissible pair such that $\chi$ has level $n \geq 1$.  We set $\psi_E = \psi \circ Tr_{E/F}$.  Let $\alpha(\chi) \in \mathfrak{p}_E^{-n}$ such\footnote{In \cite{bushnellhenniart}, the notation $\alpha$ is used.  We prefer to use the notation $\alpha(\chi)$ since this element depends on the character $\chi$.} that $\chi(1+x) = \psi_E(\alpha(\chi) x) \ \ \forall x \in \mathfrak{p}_E^{\llcorner n/2 \lrcorner + 1}$.  To the pair of data $(E/F, \chi)$ and $\alpha(\chi)$, one can attach a supercuspidal representation $\pi_{\chi}$ of $GL(2,F)$ (see \cite[Chapter 20]{bushnellhenniart}).  In general, let $(E/F, \chi)$ be an arbitrary admissible pair of level $n \geq 1$.  There is a character $\phi$ of $F^*$ and a character $\chi'$ of $E^*$ such that $(E/F, \chi')$ is a minimal admissible pair and $\chi = \chi' \phi_E$.  We define $\pi_{\chi} = \phi \pi_{\chi'}$.  Let $\mathbb{A}_2^0(F)$ denote the set of equivalence classes of all irreducible supercuspidal representations of $GL(2,F)$. Then together with Proposition (\ref{depthzerogl2}), we have the following theorem.

\begin{theorem}{\label{positivedepth}}
The map $(E/F, \chi)$ $\mapsto \pi_{\chi}$ induces a bijection $$\mathbb{P}_2(F) \rightarrow \mathbb{A}_2^0(F) \ \ if \
\ p \neq 2$$

Furthermore, if $(E/F, \chi) \in \mathbb{P}_2(F)$, then:

(i) $\omega_{\pi_{\chi}} = \chi|_{F^*}$

(ii) if $\phi$ is a character of $F^*$, then $\pi_{\chi \phi_E} = \phi \pi_{\chi}$.
\end{theorem}

\proof

See \cite[Section 20.2]{bushnellhenniart}.

\qed

\subsection{Weil parameters}\label{weilparametersgl2}

Let $\mathbb{G}_2^0(F)$ be the set of equivalence classes of irreducible smooth two-dimensional representations of $W_F$. Recall that there is a local Artin reciprocity isomorphism given by $W_E^{ab} \cong E^*$.  Then, if $(E/F, \xi) \in \mathbb{P}_2(F)$, $\xi$ gives rise to a character of $W_E^{ab}$, which we can pullback to a character, also denoted $\xi$, of $W_E$.  We can then form the induced representation $\varphi_{\xi} = Ind_{W_E}^{W_F} \xi$ of $W_F$.

\begin{theorem}\label{G_2(F)}
Suppose the residual characteristic of $F$ is not $2$.  If $(E/F, \xi)$ is an admissible pair, the representation $\varphi_{\xi}$ of $W_F$ is irreducible.  The map $(E/F, \xi) \mapsto \varphi_{\xi}$ induces a bijection $$\mathbb{P}_2(F) \rightarrow \mathbb{G}_2^0(F)$$
\end{theorem}

\proof
See \cite[Chapter 33]{bushnellhenniart}.
\qed

We therefore have canonical bijections

\begin{equation}
\mathbb{P}_2(F) \rightarrow \mathbb{A}_2^0(F), \ \ \ \ \ \ \mathbb{P}_2(F) \rightarrow \mathbb{G}_2^0(F) \ \label{bijections}
\end{equation}
  $$(E/F, \xi) \mapsto \pi_{\xi}, \ \ \ \ \ \ \ (E/F, \xi) \mapsto \varphi_{\xi},$$

\noindent given by Theorem (\ref{G_2(F)}) and Theorem (\ref{positivedepth}).  Combining both of these bijections, we obtain a bijection

\begin{equation}
\mathbb{G}_2^0(F) \rightarrow \mathbb{A}_2^0(F) \ \label{eq:naive}
\end{equation}
  $$\varphi_{\xi} \mapsto \pi_{\xi}$$

However, this bijection is not the local Langlands correspondence.  The reason is as follows. If $(E/F, \xi) \in \mathbb{P}_2(F)$, then by \cite[Proposition 29.2]{bushnellhenniart}, representation $\varphi_{\xi}$ has determinant $\aleph_{E/F} \otimes \xi|_{F^*}$, whereas $\pi_{\xi}$ has central character $\xi|_{F^*}$.  This contradicts one of the requirements of the local Langlands correspondence which says that if $\varphi : W_F \rightarrow GL(2,\mathbb{C})$ is an irreducible Weil parameter, then $\omega_{\pi(\varphi)} =$ det($\varphi$) where $\pi(\varphi)$ is the supercuspidal representation of $GL(2,F)$ that corresponds to $\varphi$ under the local Langlands correspondence.  We must therefore systematically modify the bijection (\ref{eq:naive}), which we proceed to do now.

If $K/F$ is a finite separable extension, let $\lambda_{K/F}(\psi)$ denote the Langlands constant, as in [2, 34.3].

\begin{definition}
Let $(E/F, \xi)$ be an admissible pair in which $E/F$ is unramified.  Define $\Delta_{\xi}$ to be the unique quadratic unramified character of $E^*$.
\end{definition}

Let $\mu_F$ denote the group of roots of unity in $F$ of order prime to to the residual characteristic of $F$. Let $E/F$ be a totally tamely ramified quadratic extension, let $\varpi$ be a uniformizer of $E$, and let $\beta \in E^*$.  Since $U_E = \mu_E U_E^1 = \mu_F U_E^1$, there is a unique root of unity $\zeta(\beta, \varpi) \in \mu_F$ such that $$\beta \varpi^{-v_E(\beta)} = \zeta(\beta, \varpi) \ (mod \ U_E^1).$$

\begin{definition}\label{deltachi}

(i) Let $(E/F, \xi) \in \mathbb{P}_2(F)$ be a minimal admissible pair such that $E/F$ is totally tamely ramified.  Let $n$ be the level of $\xi$ and let $\alpha(\chi) \in \mathfrak{p}_E^{-n}$ satisfy $\xi(1+x) = \psi_E(\alpha(\chi) x), x \in \mathfrak{p}_E^n$.  There is a unique character $\Delta = \Delta_{\xi}$ of $E^*$ such that:

$$\Delta|_{U_E^1} = 1, \ \ \Delta|_{F^*} = \aleph_{E/F},$$ $$\Delta(\varpi) = \aleph_{E/F}(\zeta(\alpha(\chi), \varpi)) \lambda_{E/F}(\psi)^n$$

(ii) Let $(E/F, \xi) \in \mathbb{P}_2(F)$ and suppose that $E/F$ is totally tamely ramified.  Write $\xi = \xi' \chi_E$ for a minimal admissible pair $(E/F, \xi')$ and a character $\chi$ of $F^*$.  Define $$\Delta_{\xi} = \Delta_{\xi'}$$
\end{definition}

\begin{theorem}{$\mathbf{Tame \ Local \ Langlands \ Correspondence}$}{\label{tamellc}}

Suppose the residual characteristic of $F$ is not $2$.  For $\varphi \in \mathbb{G}_2^0(F)$, define $\pi(\varphi) = \pi_{\xi \Delta_{\xi}}$ in the notation of (\ref{bijections}) for any $(E/F, \xi) \in \mathbb{P}_2(F)$ such that $\varphi \cong \varphi_{\xi}$.  The map $$\pi : \mathbb{G}_2^0(F) \rightarrow \mathbb{A}_2^0(F)$$ is the local Langlands correspondence for $GL(2,F)$.
\end{theorem}

\proof
See \cite[Chapter 34]{bushnellhenniart}.
\qed

\begin{proposition}\label{centralcharacter}
If $\varphi \in \mathbb{G}_2^0(F)$ and $\pi = \pi(\varphi)$, then $\omega_{\pi} =$ det($\varphi$).
\end{proposition}

\proof
See \cite[Chapter 33]{bushnellhenniart}
\qed

\section{Positive depth supercuspidal character formulas for $PGL(2,F)$ and $GL(2,F)$}

\subsection{Setup}\label{setup}

In this chapter we prove Theorems \ref{gl2theorem1} and \ref{gl2theorem2} for the positive depth supercuspidal representations of $GL(2,F)$.  Note that since we working here with $GL(2,F)$, we have that $L = F$ and so $\aleph_{EL/L} = \aleph_{E/F}$ and $| \ |_{EL} = | \ |_E$.  Moreover, $\epsilon(s) = (-1, \Delta)^{\ell(s) (\ell+1)} = (-1, \Delta)^{\ell(s)}$ since $\ell = 2$ here.  Finally, $\tau_o$ is a character of $E^*$ such that $\tau_o|_{F^*} = \aleph_{E/F}$.

\begin{definition}
Let $\tilde\chi$ be a genuine character of $T(F)_{\tau \circ \rho}$.  We define the function $F(\tilde\chi) : T(F)^{reg} \rightarrow \mathbb{C}$ by $$F(\tilde\chi)(z) = \epsilon(\tilde\chi, \Delta^+, \tau) \frac{\displaystyle\sum_{s \in W} \epsilon(s)\tilde\chi({}^s w)}{\tau(\Delta^0(z,\Delta^+)) (\tau \circ \rho)(w)}, \ \forall z \in T(F)^{reg}$$
\end{definition}

We will define all of the notation in Theorem (\ref{gl2theorem1}) in the next several sections, including $n(z)$, $\epsilon(\tilde\chi, \Delta^+, \tau)$, and the definition of regular.  We will also regularly use the fact that $W(G(F),T(F)) = Aut(E/F)$.  Notice that there is only one difference between this proposed character formula $F(\tilde\chi)$ and the character formula for discrete series of real reductive groups.  If we were to literally transport the character formula of Theorem (\ref{Harish-Chandra}) to the $p$-adic case, then the denominator $\Delta^0(h, \Delta^+) \rho(\tilde h)$ would take values in $K^*$, where $K/F$ is the minimal splitting field of the elliptic torus $T$. This would be problematic since characters must take values in $\mathbb{C}$.  Therefore, a natural thing to try is to introduce a $\mathbb{C}^*$-valued character $\tau$ of $K^*$ into the denominator in order that the denominator takes values in $\mathbb{C}^*$.

We note that all of our calculations in the next two chapters will assume that we have chosen the standard positive set of roots of $GL(2,\overline{F})$ with respect to the standard split maximal torus.  Our main results, however, will be seen to be independent of any choice of positive roots.

Now let $\varphi$ be a supercuspidal Weil parameter for $GL(2,F)$.  We will show later in this section how to construct a regular genuine character, $\tilde\chi$, of $T(F)_{\tau \circ \rho}$, from $\varphi$.  We will then proceed to prove Theorem (\ref{gl2theorem2}).

\noindent Let $F$ denote a non-Archimedean local field of characteristic zero with residual characteristic coprime to $2$.  Let $E/F$ be a quadratic extension.  Write $E = F(\sqrt{\Delta})$ for some $\Delta \in F^*$ and let $\delta := \sqrt{\Delta}$.  $E^*$ embeds as an elliptic torus in $GL(2,F)$ via the map $$E^* \hookrightarrow GL(2,F)$$ $$a + b \delta \mapsto \mat{a}{b}{b \Delta}{a}$$ and therefore $E^*/F^*$ embeds in $PGL(2,F)$ as an elliptic torus as well.

We now introduce a notion of regularity that we will need.  Let $E/F$ be a tamely ramified quadratic extension and $\chi$ a character of $E^*$.  Recall that $N$ denotes the norm map from $E$ to $F$.

\begin{definition}
$\chi$ is called \emph{regular} if $\chi$ does not factor through $N$.  If $\chi$ is regular, we call the pair $(E/F, \chi)$ a \emph{regular pair}.
\end{definition}

All definitions we have made in the previous chapter for admissible pairs, we also make for regular pairs and regular characters.  For example, as we defined the notion of minimal admissible pair, we make the same definition for minimal regular pair.  In particular, we also define the character twists $\Delta_{\chi}$ for a regular pair $(E/F, \chi)$ exactly the same way they were defined for admissible pairs.  For example, if $(E/F, \chi)$ is a regular pair where $E/F$ is ramified, then $\Delta_{\chi}$ is the character of $E^*$ from Definition \ref{deltachi}.  Given a regular pair $(E/F, \chi)$, one may also construct a supercuspidal representation $\pi_{\chi}$ as in the previous chapter, but this construction is not one to one.  Our constructions and results do not require the stronger notion of admissible pair.  We will sometimes say that $\chi$ is regular when the field $E$ is understood.

We now briefly reexplain why double covers of tori play a role, now in terms of the notion of regular pairs.  Let $\varphi$ be a supercuspidal Weil parameter for $PGL(2,F)$.  Then $\varphi = Ind_{W_E}^{W_F}(\chi)$, for some regular pair $(E/F, \chi)$.  Since we are using the notion of regular pair here rather than admissible pair, there may be a choice involved here.  That is, there may be another regular pair $(E_1/F, \chi_1)$ such that $\varphi = Ind_{W_{E_1}}^{W_F}(\chi_1)$ as well.  However, this will not matter, and we will show that our results and constructions are independent of all choices.  By Proposition (\ref{centralcharacter}), $\chi|_{F^*} = \aleph_{E/F}$. Therefore, the supercuspidal Weil parameters for $PGL(2,F)$ naturally give rise to regular pairs $(E/F, \chi)$ where $\chi|_{F^*} = \aleph_{E/F}$.  Such a $\chi$ is a genuine character of the double cover $E^* / N(E^*)$ of $E^* / F^*$ as in the introduction.  We note that the double cover $E^* / N(E^*)$ splits if and only if $(-1, \Delta) = 1$ (see \cite{adams}).  In fact, this double cover $E^* / N(E^*)$ is none other than an analogue of the $\rho$-cover that appears in the theory over the reals (see Definition (\ref{rhocoverdefinition})).  We explain this now.

Relative to the standard positive root of $PGL(2,F)$, let $\rho$ be half the positive root.  An elliptic torus in $PGL(2,F)$ is of the form $T(F) = E^* / F^*$.  Fix a character $\tau_o$ of $E^*$ whose restriction to $F^*$ is $\aleph_{E/F}$, and set $\tau := \tau_o | \ |_E$.  Recall the denominator $$\tau(\Delta^0(z, \Delta^+)) (\tau \circ \rho)(w)$$ that was defined in Theorem (\ref{gl2theorem1}).  Although $\tau \circ \rho$ is not naturally a function on $E^* / F^*$ since in particular $\rho$ is not naturally a function on $E^* / F^*$, it is by definition a function on the $\tau \circ \rho$-cover of $E^* / F^*$.   Recall Definition (\ref{tauofrhocover}).  Then in our case, $T(F)_{\tau \circ \rho} = \{(z,\lambda) \in E^* / F^* \times \mathbb{C}^* : \tau(2 \rho(z)) = \lambda^2 \}$.

\begin{lemma}\label{doublecovertoriisomorphism}
$E^* / N(E^*) \cong T(F)_{\tau \circ \rho}$.
\end{lemma}

\proof
Define the map $$E^* / N(E^*) \xrightarrow{\kappa} T(F)_{\tau \circ \rho}$$ $$[w] \mapsto ([w], \tau_o(w) |2 \rho([w])|^{1/2})$$ where on the right hand side, $[w]$ lives in $E^* / F^*$. To show injectivity, suppose $\kappa([w]) = ([1], 1)$, where $w \in E^*$.  Since $[w] = [1]$, we get $w \in F^*$.  But since $\tau_o(w) |2 \rho([w])|^{1/2} = 1$, we get that $w \in N(E^*)$ since $\tau_o|_{F^*} = \aleph_{E/F}$ and since $|2 \rho([w])| = 1$ since $[w] = 1$.  To show surjectivity, suppose $([w], \lambda) \in T(F)_{\tau \circ \rho}$, where $w \in E^*$.  Then, by definition of $T(F)_{\tau \circ \rho}$, we get that $\tau(2 \rho([w])) = \lambda^2$.  This means that $\tau_o(w / \overline{w}) |2 \rho([w])| = \lambda^2$.  But $\tau_o$ is trivial on the norms, so we have that $\tau_o(w / \overline{w}) = \tau_o(w^2 / N(w)) = \tau_o(w)^2$.  Therefore, $\lambda = \pm \tau_o(w) |2 \rho([w])|^{1/2}$.  If $\lambda = \tau_o(w) |2 \rho([w])|^{1/2}$, then we get that $\kappa([w]) = ([w], \lambda)$.  If $\lambda = - \tau_o(w) |2 \rho([w])|^{1/2}$, then let $x \in F^* \setminus N(E^*)$.  Then $\kappa([x w]) = ([w], \lambda)$.  Therefore, $\kappa$ is surjective.
\qed

We note that the importance of the term $|2 \rho([w])|^{1/2}$ comes from the fact that $$|D([w])|^{\frac{1}{2}} = |\Delta^0([w], \Delta^+)| |2 \rho([w])|^{1/2},$$ an observation made in Chapter 3.  The reason why this is important is that the term $|D(w)|^{1/2}$ appears in the supercuspidal characters (see Section (\ref{theconstantepsilon})).  We will need this fact for the character formulas for $GL(n,F)$, where $n$ is prime.

Now let's write down the character formula for a supercuspidal representation of $PGL(2,F)$.  Let $\varphi : W_F \rightarrow GL(2,\mathbb{C})$ be a supercuspidal parameter for $PGL(2,F)$ so that $\varphi = Ind_{W_E}^{W_F}(\chi)$ for some regular pair $(E/F, \chi)$. As discussed earlier, this gives us a genuine character $\tilde\chi$ of $E^* / N(E^*)$.

\begin{definition}
A genuine character $\tilde\eta$ of $E^* / N(E^*)$ is called \emph{regular} if $(E/F, \eta)$ is regular, where $\eta$ is the pullback of $\tilde\eta$ to $E^*$.  A genuine character $\tilde\lambda$ of $T(F)_{\tau \circ \rho}$ is called \emph{regular} if $\tilde\lambda \circ \kappa$ is regular.
\end{definition}

Now recall from Theorem (\ref{gl2theorem1}) the proposed character formula $F(\tilde\chi)$.  We naturally constructed a genuine character $\tilde\chi$ of $E^* / N(E^*)$.  However, the functions in $F(\tilde\chi)$ have domain $T(F)_{\tau \circ \rho}$.  Recall that $T(F)_{\tau \circ \rho} \cong E^* / N(E^*)$ by Lemma (\ref{doublecovertoriisomorphism}), so we can pull the function $(\tau \circ \rho)(w)$ and the Weyl group action in $F(\tilde\chi)$ back to $E^* / N(E^*)$ via this isomorphism, and leave our constructed $\tilde\chi$ as living on $E^* / N(E^*)$.  That is, we consider $$F(\tilde\chi)(z) = \epsilon(\tilde\chi, \Delta^+, \tau) \frac{\displaystyle\sum_{s \in W} \epsilon(s)\tilde\chi({}^s [w])}{\tau(\Delta^0(z,\Delta^+)) (\tau \circ \rho)(\kappa([w]))}, \ \ z \in T(F)^{reg}$$ where $[w] \in E^* / N(E^*)$ such that $\Pi(\kappa([w])) = z$.  Unwinding the definitions, we see that \\ $(\tau \circ \rho)(\kappa([w])) = \tau_o(w) |2 \rho([w])|^{1/2} \ \forall [w] \in E^* / N(E^*)$, where we also write $[w]$ as the element in $E^* / F^*$.

We also need to define the Weyl group action.  The Weyl group action on the $\tau \circ \rho$-cover is obtained as follows. If $([w],\lambda)$ is an element of $T(F)_{\tau \circ \rho}$, then analogously to the real case (recall equation (\ref{eq:weylgroupreal}) in Chapter \ref{realgroupschapter}), define $s([w],\lambda) = (s[w],\lambda \tau((s^{-1} \rho - \rho)([w])))$ for $s \in W = W(G(F),T(F)) = Aut(E/F)$.  Simplifying this expression, we get $s([w],\lambda) = ([\overline{w}], \lambda \tau(\overline{w}/w))$ when $s \in W$ is nontrivial.  Then, since our character formula lives on $E^* / N(E^*)$, we must pull back this action from $T(F)_{\tau \circ \rho}$ to $E^* / N(E^*)$ via $\kappa$.  Doing this, we see that we get $s [w] = \kappa^{-1}(s \kappa([w])) = \kappa^{-1}(s ([w], \tau_o(w) |2 \rho([w])|^{1/2})) = \kappa^{-1}([\overline{w}], \tau_o(w) |2 \rho([w])|^{1/2} \tau(\overline{w}/w)) = \kappa^{-1}([\overline{w}], \tau_o(\overline{w}) |2 \rho([\overline{w}])|^{1/2}) = [\overline{w}] \ \forall [w] \in E^* / N(E^*)$ when $s \in W = Aut(E/F)$ is nontrivial, since $|2 \rho([w])| = |w / \overline{w}| = 1 \ \forall w \in E^*$.

We note that the definition of regularity for a genuine character of $T(F)_{\tau \circ \rho}$ is analogous to the definition of regularity for a genuine character $\tilde\lambda$ of $T(\mathbb{R})_{\rho}$ for real groups.

Then, pulling $(\tau \circ \rho)(w)$ and the Weyl group action back to $E^* / N(E^*)$ via $\kappa$, we get $$F(\tilde\chi)(z) = \epsilon(\tilde\chi, \Delta^+, \tau) \frac{\tilde\chi([w]) + (-1, \Delta) \tilde\chi([\overline{w}])}{\tau_o(1-1/z) |\Delta^0(z, \Delta^+)| \tau_o(w) |2 \rho(z)|^{1/2} }$$ where $z \in E^*/F^*$ and $[w] \in E^* / N(E^*)$ is some element that maps to $z$ under the map \\ $E^* / N(E^*) \rightarrow E^* / F^*$.  Note that our formula simplifies: $$F(\tilde\chi)(z) = \epsilon(\tilde\chi, \Delta^+, \tau) \frac{\chi(w) + (-1, \Delta) \chi(\overline{w})}{\tau_o(w-\overline{w}) |D(z)|^{1/2}}, \ z \in T(F)^{reg}$$ where $w \in E^*$ is any element that maps to $z \in T(F) = E^* / F^*$ under canonical map $E^* \rightarrow E^* / F^*$.  The reason is that $|D(z)|^{1/2} = |\Delta^0(z, \Delta^+)| |2 \rho(z)|^{1/2}$ from Chapter 3.

We also note that if we had made the other choice of $\Delta^+$, the denominator in our character formula would include the term $\tau_o(\overline{w} - w)$ instead of $\tau_o(w - \overline{w})$.  However, because our definition of $\epsilon(\tilde\chi, \Delta^+, \tau)$ includes the term $\epsilon(\Delta^+)$ (see Section (\ref{theconstantepsilon})), our overall character formula $F(\tilde\chi)$ remains the same regardless of the choice of positive root.  The same line of reasoning is true for the case of $GL(2,F)$, which we present next.

\

Let us now compute our proposed character formula for $GL(2,F)$.  Let $\rho$ be half the standard positive root of $GL(2,F)$.  An elliptic torus in $GL(2,F)$ is of the form $T(F) = E^*$.  We now introduce a cover which is isomorphic to $T(F)_{\tau \circ \rho}$.

\begin{definition}
Let $E^* / N(E^*) \rightarrow E^* / F^*$ be the canonical projection map.  We define $E^* \times_{E^*/F^*} E^* / N(E^*)$ as the group arising in the following pullback diagram:
$$
\begin{CD}
E^* \times_{E^*/F^*} E^* / N(E^*) @>>> E^* / N(E^*)\\
@VVV @VVV\\
E^* @>w \mapsto [w]>> E^*/F^*
\end{CD}
$$
That is, $E^* \times_{E^*/F^*} E^* / N(E^*) = \{(w,z) \in E^* \times E^* / N(E^*) : [w] = [z] \in E^* / F^* \}$
\end{definition}

\begin{lemma}\label{doublecovertoriisomorphismgl2}
$E^* \times_{E^*/F^*} E^* / N(E^*) \cong T(F)_{\tau \circ \rho}$
\end{lemma}

\proof
An explicit isomorphism is given by $$E^* \times_{E^*/F^*} E^* / N(E^*) \xrightarrow{\kappa} T(F)_{\tau \circ \rho}$$ $$(w,[z]) \mapsto (w, \aleph_{E/F}(z/w) \tau_o(w) |2 \rho(w)|^{1/2})$$ To see that this is injective, note that if $w = 1$, then $z \in F^*$ by definition of pullback.  But then $\aleph_{E/F}(z) = 1$ implies that $z \in N(E^*)$.  To see surjectivity, suppose that $(w, \lambda) \in T(F)_{\tau \circ \rho}$.  Then by definition of the pullback, we get $\tau_o(w / \overline{w}) |2 \rho(w)| = \lambda^2$.  But $\tau_o(w / \overline{w}) = \tau_o(w^2 / N(w)) = \tau_o(w)^2$.  Thus, $\lambda = \pm \tau_o(w) |2 \rho(w)|^{1/2}$.  If $\lambda = \tau_o(w) |2 \rho(w)|^{1/2}$, then $\kappa(w, [w]) = (w, \lambda)$.  If $\lambda = - \tau_o(w) |2 \rho(w)|^{1/2}$, then $\kappa(w, [x w]) = (w, \lambda)$, where $x \in F^* \setminus N(E^*)$.
\qed

Now let's write down the character formula for a supercuspidal representation of $GL(2,F)$.  Now let $\varphi : W_F \rightarrow GL(2,\mathbb{C})$ be a supercuspidal parameter so that $\varphi = Ind_{W_E}^{W_F}(\chi)$ for some regular pair $(E/F, \chi)$.  Then this canonically gives a genuine character $\tilde\chi$ of $E^* \times_{E^*/F^*} E^* / N(E^*)$ as follows.  Define $\tilde\chi(w,[z]) := \chi(w) \aleph_{E/F}(z/w)$.

\begin{definition}
A genuine character $\tilde\eta$ of $E^* \times_{E^*/F^*} E^* / N(E^*)$ is called \emph{regular} if $(E/F, \eta)$ is regular, where $\eta(w) := \tilde\eta(w,[z]) \aleph_{E/F}(z/w)$.  A genuine character $\tilde\lambda$ of $T(F)_{\tau \circ \rho}$ is called \emph{regular} if $\tilde\lambda \circ \kappa$ is regular.
\end{definition}

We have therefore given a map $\widehat{E^*} \rightarrow (E^* \times_{E^*/F^*} E^* / N(E^*))^{\wedge}$ given by $\eta \mapsto \tilde\eta$, where $\tilde\eta(w,[z]) := \eta(w) \aleph_{E/F}(z/w)$.  Note that we have a canonical map in the other direction, $(E^* \times_{E^*/F^*} E^* / N(E^*))^{\wedge} \rightarrow \widehat{E^*}$, given by $\tilde\eta \mapsto \eta$, where $\eta(w) := \tilde\eta(w,[z]) \aleph_{E/F}(z/w)$.  We will regularly go back and forth between characters of $E^*$ and genuine characters of $E^* \times_{E^*/F^*} E^* / N(E^*)$.  In particular, when we write $\tilde\chi$, a genuine character of $E^* \times_{E^*/F^*} E^* / N(E^*)$, we will sometimes keep in mind that there is a canonical character $\chi$ of $E^*$ that $\tilde\chi$ comes from via the above maps.

Now recall the proposed character formula $F(\tilde\chi)$ from Theorem (\ref{gl2theorem1}).  We have naturally constructed a genuine character $\tilde\chi$ of $E^* \times_{E^*/F^*} E^* / N(E^*)$.  However, the functions in $F(\tilde\chi)$ have domain $T(F)_{\tau \circ \rho}$.  Recall that $T(F)_{\tau \circ \rho} \cong E^* \times_{E^*/F^*} E^* / N(E^*)$, so we can pull the function $(\tau \circ \rho)(w)$ and the Weyl group action in $F(\tilde\chi)$ back to $E^* \times_{E^*/F^*} E^* / N(E^*)$ via this isomorphism $\kappa$ from Lemma (\ref{doublecovertoriisomorphismgl2}), and leave our constructed $\tilde\chi$ as living on $E^* \times_{E^*/F^*} E^* / N(E^*)$.  That is, we consider $$F(\tilde\chi)(w) = \epsilon(\tilde\chi, \Delta^+, \tau) \frac{\displaystyle\sum_{s \in W} \epsilon(s)\tilde\chi({}^s (w,[z]))}{\tau(\Delta^0(w,\Delta^+)) (\tau \circ \rho)(\kappa(w,[z]))}, \ \ w \in T(F)^{reg}$$ where $(w,[z]) \in E^* \times_{E^* / F^*} E^* / N(E^*)$ such that $\Pi(\kappa((w,[z]))) = w$.  Unwinding the definitions, we see that $(\tau \circ \rho)(\kappa((w,[z]))) = \aleph_{E/F}(z/w) \tau_o(w) |2 \rho(w)|^{1/2} \ \ \forall \ (w,[z]) \in E^* \times_{E^* / F^*} E^* / N(E^*)$.

We also need to define the Weyl group action.  As in the case of $PGL(2,F)$, the action is $s(w,\lambda) = (sw,\lambda \tau((s^{-1} \rho - \rho)(w)))$ for $s \in W$.  Pulling this back to $E^* \times_{E^* / F^*} E^* / N(E^*)$ via $\kappa$, this simplifies to $s(w,[z]) = (\overline{w}, [\overline{z}]) \ \forall (w,[z]) \in E^* \times_{E^* / F^*} E^* / N(E^*)$ when $s \in W$ is nontrivial.

Pulling back $(\tau \circ \rho)(w)$ and the Weyl group action to $E^* \times_{E^*/F^*} E^* / N(E^*)$ via $\kappa$, we get
$$F(\tilde\chi)(w) = \epsilon(\tilde\chi, \Delta^+, \tau) \frac{\displaystyle\sum_{s \in W} \epsilon(s) \chi({}^s w) \aleph_{E/F}({}^s (z/w))}{\tau_o(1 - \overline{w}/w) |\Delta^0(w, \Delta^+)| \tau_o(w) \aleph_{E/F}(z/w) |2 \rho(w)|^{1/2}}  =$$  $$ \epsilon(\tilde\chi, \Delta^+, \tau) \frac{\chi(w) + (-1, \Delta) \chi(\overline{w})}{\tau_o(w-\overline{w}) |D(w)|^{1/2}}, \ \ w \in T(F)^{reg}$$ since $|D(w)|^{1/2} = |\Delta^0(w, \Delta^+)| |2 \rho(w)|^{1/2}$ from Chapter 3.
We will see that our proposed character formulas for $GL(2,F)$ and $PGL(2,F)$ are independent of the choice of $\tau$.

Summing up, we have given a method of assigning a conjectural character formula for a supercuspidal representation of $GL(2,F)$ or $PGL(2,F)$ to a supercuspidal Weil parameter $\varphi$ of $GL(2,F)$ or $PGL(2,F)$, respectively, given by

$$\varphi \mapsto \tilde\chi \in \widehat{T(F)}_{\tau \circ \rho} \mapsto F(\tilde\chi)$$

\subsection{The constant $\epsilon(\tilde\chi, \Delta^+, \tau)$}\label{theconstantepsilon}

We now turn to the question of defining the constant $\epsilon(\tilde\chi, \Delta^+, \tau)$.  We recall the main theorem describing the distribution characters of positive depth supercuspidal representations of $GL(n,F)$, where $n$ is prime.  We note that there is an analogous definition of regular pair $(E/F, \chi)$ when $E/F$ has degree $n$, and this is discussed further in Section (\ref{moythesis}).

\begin{theorem}{\label{DeBacker1}}{\cite[Theorem 5.3.2]{debacker}} \footnote{In \cite{debacker}, $X_\pi$ is the notation used instead of $\alpha(\chi)$ (recall the notation $\alpha(\chi)$ from Section 4.3).  The notation $X_\pi$ is a bit misleading, because the element $X_\pi$ depends on $\chi$, not just $\pi$.  Since the notation in \cite{debacker} and \cite{bushnellhenniart} differ, we need to choose a set of notation.  We prefer to use the notation $\alpha(\chi)$.}  Let $(E/F, \chi)$ be a regular pair where $E/F$ has degree $n$ and $\chi$ has positive level, and write $G' = E^*$.  Let $\pi = \pi_{\chi}$ be the associated positive depth supercuspidal representation of $GL(n,F)$ given by Theorem (\ref{positivedepth}).  Then

\begin{equation*}
\frac{\theta_{\pi}(\gamma)}{deg(\pi)} = \left\{
\begin{array}{rl}
C \lambda(\sigma) \displaystyle\sum_{\overline{w} \in W} \chi({}^w t) & \text{if } n(\gamma) = 0 \ and \ \gamma = {}^g t \ with \ g \in G \ and \ t \in G'\\
C \displaystyle\sum_{\overline{w} \in W} \chi({}^w t) \gamma(\alpha(\chi), {}^w Y) & \text{if } 0 < n(\gamma) \leq r/2 \ and \ \gamma = {}^g t \ where \ t = z(1 + Y) \ with \\

 & z \in Z, g \in G, \ and \ Y \in \mathfrak{g}_{n(\gamma)}' \setminus (\mathfrak{z}_{n(\gamma)} + \mathfrak{g}_{n(\gamma)^+}')\\
\chi(z) \mu_{\alpha(\chi)}(Y) & \text{if } n(\gamma) > r/2 \ and \ \gamma = z(1 + {}^g Y) \ with \ g \in G, z \in Z, \\
 & and \ Y \in V_{n(\gamma)} \\
0 & \text{otherwise }
\end{array} \right.
\end{equation*}

\end{theorem}

Here, $C := c_{\psi}(\mathfrak{g'}) c_{\psi}^{-1}(\mathfrak{g}) |D(\gamma)|^{-1/2} |\eta(\alpha(\chi))|^{-1/2}$ is defined in \cite[Section 5.3]{debacker}.\footnote{In \cite{debacker}, $\Lambda$ is used to denote a fixed additive character of $F$.  In \cite{bushnellhenniart}, $\psi$ is used to denote a fixed additive character of $F$.  We prefer to use the notation $\psi$.}  The notation is explained in \cite[Chapters 4-5]{debacker} and we will mostly follow the same notation.  In particular, $D$ and $\eta$ denote the Weyl discriminants of $G$ and $\mathfrak{g}$, respectively.

To define $\epsilon(\tilde\chi, \Delta^+, \tau)
$, we need to calculate the constant $\gamma(\alpha(\chi), Y)$ in the above theorem.  Recall $G = GL(2,F), \mathfrak{g} = \mathfrak{gl}(2,F), G' = E^*, \mathfrak{g}' = E$.  We have a direct sum decomposition $\mathfrak{g} = \mathfrak{g}' + \mathfrak{g}^{\bot}$ where the perpendicular is taken with respect to the trace form $<Z_1, Z_2> := Tr(Z_1 Z_2) \ \forall Z_1, Z_2 \in \mathfrak{g}$.  Let $V,W \in \mathfrak{g}^{\bot}$.  Then we define $$Q_{(\alpha(\chi),Y)}(V,W) := (trace([\alpha(\chi),W][V,Y])/2).$$ $Q_{(\alpha(\chi),Y)}(V,W)$ is a non-degenerate, symmetric, bilinear form on $\mathfrak{g}^{\bot}$.  Then, $\gamma_{(\alpha(\chi), Y)}$ is by definition the Weil Index of $\psi \circ Q_{(\alpha(\chi),Y)}$ (see \cite{rao} and Definition (\ref{rao0})).  Let us calculate $\gamma_{(\alpha(\chi), Y)}$.

First, note that $\mathfrak{g}'$ may be embedded in $M_2(F)$ by

$$\mathfrak{g}' \hookrightarrow \mathfrak{g}$$ $$\ \ \ \ \ \ \ \ \ \ \ \ \ \ \ \ \ \ \ \ \ \ \ \ a + d \delta \mapsto \mat{a}{d}{d \Delta}{a}, \ \ a,d \in F$$ where $E = F(\delta)$ with $\delta = \sqrt{\Delta}$.  The following two Lemmas are elementary.

\begin{lemma}
$$\mathfrak{g}^{\bot} = \left\{ \mat{a}{b}{-b \Delta}{a} : a,b \in F \right\}$$
\end{lemma}

\begin{lemma}\label{randomness}
If $\alpha(\chi) = \mat{a}{x}{x \Delta}{a}$ and $Y = \mat{t}{y}{y \Delta}{t}$ are arbitrary elements in $\mathfrak{g'}$, then the matrix of the quadratic form $Q_{(\alpha(\chi),Y)}$ is $$\mat{4xy \Delta}{0}{0}{-4xy \Delta^2}$$
\end{lemma}

\begin{lemma}
Let $\alpha(\chi) = \mat{a}{x}{x \Delta}{a}$ and $Y = \mat{t}{y}{y \Delta}{t} \in \mathfrak{g'}$.  Then $$\gamma(\alpha(\chi), Y) = (x, \Delta)_F (y, \Delta)_F \ \gamma_F(\Delta, \psi)$$
\end{lemma}

\proof
By Definition (\ref{rao2}), we have $$\gamma(\psi \circ Q_{(\alpha(\chi),Y)}) = h_F(Q_{(\alpha(\chi),Y)})\gamma_F(\psi)^n \gamma_F(det(Q_{(\alpha(\chi),Y)}),\psi)$$  In our case, $n = 2$, and we also have that by Lemma (\ref{rao3}), $$h_F(Q_{(\alpha(\chi),Y)}) = (4xy \Delta, -4xy \Delta^2)_F = (\Delta, -xy)_F$$  Then, Lemmas (\ref{rao1}) and (\ref{randomness}) imply that $$\gamma_F(\psi)^2 = \gamma_F(-1, \psi)^{-1} \ \mathrm{and}$$  $$\gamma_F(det(Q_{(\alpha(\chi),Y)}), \psi) = \gamma_F(-16x^2 y^2 \Delta^3, \psi) = \gamma_F(-\Delta, \psi)$$  Thus, by Lemma (\ref{rao1}), $$\gamma(\psi \circ Q_{(\alpha(\chi),Y)}) = (\Delta, -xy)_F \gamma_F(-1, \psi)^{-1} \gamma_F(-\Delta, \psi) = $$ $$(\Delta, -xy)_F \ \gamma_F(\Delta, \psi) (-1, \Delta)_F   = (x, \Delta)(y, \Delta) \gamma_F(\Delta, \psi)$$
\qed

\begin{definition}\label{definitionofgammafactor}
Let $(E/F, \chi)$ be a regular pair such that $\chi$ has positive level.  Associated to $(E/F, \chi)$ is an element $\alpha(\chi)$ (from Section (\ref{alpha})) and a supercuspidal representation $\pi := \pi_{\chi}$ via Theorem (\ref{positivedepth}).  Now, $\alpha(\chi) \in E^*$, so $\alpha(\chi) = a + x_{\chi} \delta$ for some $a, x_{\chi} \in F$.  Let $deg(\pi)$ denote the formal degree of $\pi$. Let $\Delta^+$ be a choice of a positive root of $GL(2,\overline{F})$ with respect to the diagonal maximal torus $T(\overline{F})$.  Define $\epsilon(\Delta^+)$ to be 1 if $\Delta^+$ is the standard positive root and define $\epsilon(\Delta^+)$ to be $\tau_o(-1)$ if $\Delta^+$ is the opposite root. Then, we define $\epsilon(\tilde\chi, \Delta^+, \tau) := deg(\pi)(x_{\chi}, \Delta) \gamma_F(\Delta, \psi) \tau_o(2 \delta) c_{\psi}(\mathfrak{g'}) c_{\psi}^{-1}(\mathfrak{g})  |\eta(\alpha(\chi))|^{-\frac{1}{2}} \epsilon(\Delta^+)$, where $c_{\psi}(\mathfrak{g'}), c_{\psi}(\mathfrak{g}),$ and  $\eta(\alpha(\chi))$ are defined in \cite[Chapter 5]{debacker}.
\end{definition}

In the calculations we will make throughout the rest of this chapter and the next, we will make a choice of $\Delta^+$ to be the standard set of positive roots.  Therefore, the term $\epsilon(\Delta^+)$ is just $1$, and so this term will not appear in most of our calculations and formulas.  We will show later that all of our results will be independent of the choice of $\Delta^+$.

\subsection{On certain decompositions associated to elements of $E^*$}\label{decompositions}

We need to understand the sets $n(\gamma) = 0$ and $0 < n(\gamma) \leq r/2$ (see Theorem (\ref{DeBacker1})).  We recall some relevant notions and definitions from \cite[Section 5.3, Section 3.2]{debacker}.  In general, we define a filtration on $E^*$ by setting $G_t' := 1 + \mathfrak{p}_E^{\ulcorner te \urcorner}$ for $t > 0$, where $e$ is the ramification index of $E$ over $F$ and $G' := E^*$.  We also define $G_{t^+}' := \bigcup_{s > t} G_s'$ for $t > 0$.  We let $Z(G)$ denote the center of $GL(2,F)$.

\begin{definition}
Let $w \in Z(G) G_{0^+}'$.  Then $n(w)$ is defined by $w \in Z(G) G_{n(w)}' \setminus Z(G) G_{n(w)^+}'$.  The \emph{decomposition of w} by definition is the rewriting of $w$ in the form $w = ab$ where $a \in Z(G), b \in G_{n(w)}'$ and such that $w$ is not in $Z(G) G_{n(w)^+}'$.
\end{definition}

\begin{definition}
If $w$ is not in $Z(G) G_{0^+}' = F^*(1 + \mathfrak{p}_E)$, then define $n(w) = 0$.
\end{definition}

We first deal with the situation where $E/F$ is ramified.  Assume $E = F(\sqrt{p})$.  The analogous results hold for $E = F(\sqrt{dp})$ with the same proofs, where $d \in \mathfrak{o}_F^*$ is not a square.  Note that we are not interested in the case where $w \in F^*$, since distribution characters are only defined on the regular set.

\begin{lemma}
Let $w = p^n u + p^m v \delta \in E^*$, where $u,v \in \mathfrak{o}_F^*$, and $n,m \in \mathbb{Z}$ such that $n \leq m$.  We can rewrite w as $$w = p^n u \left(1 + p^{m-n + \frac{1}{2}} \frac{v}{u} \right)$$  Thus, $w \in F^* U_E^{2m-2n+1}$.
Moreover, w is not in $F^* U_E^{2m-2n+2}$.  Therefore, the decomposition of w is
$$w = p^n u \left(1 + p^{m-n+\frac{1}{2}} \frac{v}{u} \right) = p^n u \left(1 + \sqrt{p}^{2m-2n+1} \frac{v}{u} \right)$$
\end{lemma}

\proof

Suppose by way of contradiction that $w = x(1 + s(\sqrt{p})^{2m-2n+2}), s \in U_E, x \in F^*$.  So \\ $p^n u(1 + v' p^{m-n+1/2}) = x(1+sp^{m-n+1})$, where $v' = \frac{v}{u}$.  Then $$x^{-1}p^n u(1 + v' p^{m-n+1/2}) = 1 + sp^{m-n+1}$$  Well, $x \in F^*$ is arbitrary, therefore $x^{-1} p^n u \in F^*$ is arbitrary, so the proof of the Lemma reduces to showing that there is no $y \in F^*$ such that $y(1 + v' p^{m-n+1/2}) = 1 + sp^{m-n+1}$.  By way of contradiction, suppose such a $y$ existed.  We consider power series expansions of various elements.   Let $y = p^k(y_0 + y_1p + ...),$ where $k \in \mathbb{Z}, y_0 \neq 0$, let $v' = v_0' + v_1' p + ...,$ where $v_0' \neq 0$, and let $s = s_0 + s_1 \sqrt{p} + ...  $, where $s_0 \neq 0$. Then we have $$p^k(y_0 + y_1 p + ...)(1 + (v_0' + v_1' p + ...)p^{m-n+1/2}) = 1 + (s_0 + s_1 \sqrt{p} + ...)p^{m-n+1}$$ Comparing leading coefficients, this implies that $k=0$ and $y_0 = 1$.  Therefore, $$(1 + y_1 p + ...)(1 + v_0' p^{m-n+1/2} + v_1' p^{m-n+3/2} + ...) = 1 + (s_0 + s_1 \sqrt{p} + ...)p^{m-n+1}$$  But, expanding the left hand side, we see that $$(1 + y_1 p + ...)(1 + v_0' p^{m-n+1/2} + v_1' p^{m-n+3/2} + ...) = 1 + v_0' p^{m-n+1/2} + y_1 p + ...$$  On the other hand, $1 + (s_0 + s_1 \sqrt{p} + ...)p^{m-n+1}$ does not have a $p^{m-n+1/2}$ term, and this implies that $v_0' = 0$.  But $v'$ is a unit, so we have a contradiction.
\qed

\begin{lemma}
(1)Let $w = p^n u + p^m v \delta \in E^*$ where $u = 0$ and $v \neq 0$. Then $n(w) = 0$.

(2)Let $w = p^n u + p^m v \delta \in E^*$, where $u,v \neq 0$, $n,m \in \mathbb{Z}$ such that $n > m$. Then $n(w) = 0$.
\end{lemma}

\proof
The proofs are similar as in the previous Lemma.
\qed

\

We now describe the decomposition of elements of $E^*$ when $E/F$ is unramified.

\begin{lemma}
Let $w = p^n u + p^m v \delta \in E^*$, where $u,v \in F$ both non zero, $n,m \in \mathbb{Z}$ such that $n < m$.  We can rewrite w as $$w = p^n u \left(1 + p^{m-n} \frac{\delta v}{u} \right)$$  Thus, $w \in F^* U_E^{m-n}$, so $n(w) > 0$.
Moreover, w is not in $F^* U_E^{m-n+1}$.  Therefore, the decomposition of w is $$w = p^n u \left(1 + p^{m-n} \frac{\delta v}{u} \right)$$
\end{lemma}

\proof
The proof is similar as in the ramified case.
\qed

\begin{lemma}
(1) Let $w = p^n u + p^m v \delta \in E^*$, where $u = 0$ and $v \neq 0$.  Then $n(w) = 0$.

(2) Let $w = p^n u + p^m v \delta \in E^*$, where $u,v \neq 0$ and $n = m$.  Then $n(w) = 0$.

(3) Let $w = p^n u + p^m v \delta \in E^*$, where $u,v \neq 0$ and $n > m$.  Then $n(w) = 0$.
\end{lemma}

\proof
The proofs are similar.
\qed

\subsection{On the proof that our character formulas agree with positive depth supercuspidal characters}

Here we prove that on the range $\{ z \in T(F)^{reg} : 0 \leq n(z) \leq r/2 \}$, our conjectured character formula agrees with the character of a specific positive depth supercuspidal representation of $GL(2,F)$, and this supercuspidal is the one given by the local Langlands correspondence.  Note again that when we study the ranges $\{ w \in E^* : 0 \leq n(w) \leq r/2 \}$, we do not consider elements $w \in F^*$ since distribution characters are only defined on the regular set.

In the remainder of the chapter and the next, we will deal exclusively with $GL(2,F)$, and so we set $T(F) = E^*$.  Recall from Section (\ref{setup}) that our proposed character formula simplifies to $$F(\tilde\chi)(w) = \epsilon(\tilde\chi, \Delta^+, \tau) \
\frac{\chi(w) + (-1, \Delta) \chi(\overline{w})}{ \tau_o(w-\overline{w}) |D(w)|^{1/2}}, \ \ w \in T(F)^{reg}$$ It will be useful for computational purposes to rewrite this formula as $$F(\tilde\chi)(w) = \epsilon(\tilde\chi, \Delta^+, w) \
\frac{\chi(w) + (-1, \Delta) \chi(\overline{w})}{ \tau_o(\frac{w-\overline{w}}{2 \delta})}, \ \ w \in T(F)^{reg}$$ where $\epsilon(\tilde\chi, \Delta^+, w) := deg(\pi) (x_{\chi}, \Delta)_F \gamma(\Delta, \psi) C$ and $C := c_{\psi}(\mathfrak{g'}) c_{\psi}^{-1}(\mathfrak{g}) |D(\gamma)|^{-1/2} |\eta(\alpha(\chi))|^{-1/2}$ is as in Section (\ref{theconstantepsilon}).  We will use this rewritten version for the rest of Chapter 5.  Note that $F(\tilde\chi)$ is independent of the choice of $\tau$ because $\frac{w - \overline{w}}{2 \delta} \in F^*$, and we have required only that $\tau_o|_{F^*} = \aleph_{E/F}$.  We will start by assuming that all of our regular pairs $(E/F, \chi)$ are \emph{minimal} (see Section (\ref{admissiblepairs})), and then we will show that there is no harm in assuming this, and that all of our results are true for arbitrary regular pairs.

\

First we consider the case that $E/F$ is ramified, so we may take $E = F(\sqrt{p})$ without loss of generality.  The same proofs work in the case $E = F(\sqrt{dp})$, where $d \in \mathfrak{o}_F^*$ is not a square.  We must first conduct a careful analysis of the supercuspidal characters in the $0 < n(w) \leq r/2$ range.  Recall that on the regular set, $n(w) > 0$ if and only if $w = p^n u + p^m v \delta$ where $n \leq m$.

\begin{lemma}\label{gammafactors}
$\gamma(\alpha(\chi), Y) = (x_{\chi}, \Delta)_F \gamma(\Delta, \psi) \mu(\frac{w - \overline{w}}{2 \delta}) \mu(w)$
 $\ \forall w \in E^* : n(w) > 0$ for any character $\mu$ of $E^*$ and whose restriction to $F^*$ is $\aleph_{E/F}$ and whose order is a power of 2.  Similarly,

$\gamma(\alpha(\chi), {}^s Y) = (x_{\chi}, \Delta)_F \gamma(\Delta, \psi)\mu(\frac{\overline{w} - w}{2 \delta}) \mu(\overline{w}) \ \forall w \in E^* : n(w) > 0$ where $1 \neq s \in W = W(G(F),T(F))$.
\end{lemma}

\proof
Recall that if $\alpha(\chi), Y$ embed in $\mathfrak{g}'$ as $\alpha(\chi) = \mat{a}{x_{\chi}}{x_{\chi} \Delta}{a}$ and $Y = \mat{t}{y}{y \Delta}{t}$, then $$\gamma(\alpha(\chi), Y) = (x_{\chi}, \Delta)_F (y, \Delta)_F \ \gamma_F(\Delta, \psi)$$  If $w = m + n \delta \in E^*$ where $m,n \in F^*$, and $n(w) > 0$, then the decomposition of $w$ is $w = m(1 + n/m \delta)$.  Thus, since $(n/m, \Delta) = (n, \Delta) (m, \Delta)$, and since $m = \frac{w+\overline{w}}{2}$ and $n = \frac{w - \overline{w}}{2 \delta}$, we have that $$\gamma(\alpha(\chi), Y) = (x_{\chi}, \Delta)_F \gamma(\Delta, \psi) \mu \left(\frac{w - \overline{w}}{2 \delta} \right) \mu \left(\frac{w + \overline{w}}{2} \right)$$

\noindent Now, $\mu(\frac{w + \overline{w}}{2}) = \mu(w) \mu(\frac{1}{2}(1 + \overline{w}/w))$, but since $n(w) > 0$, we have $w \in F^* U_E^1$, so $\overline{w}/w \subset  U_E^1$. So $\frac{1}{2}(1 + \overline{w}/w) \in U_E^1$, and therefore $\mu(\frac{1}{2}(1 + \overline{w}/w)) = 1$ by Lemma (\ref{infinitesquareroot}).
\qed

Therefore, we can simplify the supercuspidal characters in the $0 < n(w) \leq r/2$ range to $$\theta_{\pi}(w) = \epsilon(\tilde\chi, \Delta^+, w)  \ \frac{\phi(w) \mu(w) + (-1, \Delta) \phi(\overline{w}) \mu(\overline{w})}{\mu(\frac{w-\overline{w}}{2 \delta})} \ \forall w \in E^* : 0 < n(w) \leq r/2$$ where $\mu$ is any character of $E^*$ whose restriction to $F^*$ is $\aleph_{E/F}$ and whose order is a power of 2.  We have therefore proven the following proposition.

\begin{proposition}
$F(\tilde\chi)$ agrees with the character of the supercuspidal representation $\pi_{\chi \mu^{-1}}$ in the $0 < n(w) \leq r/2$ range, where $\mu$ is any character of $E^*$ whose restriction to $F^*$ is $\aleph_{E/F}$ and whose order is a power of 2.
\end{proposition}

We will need to investigate the $n(w) = 0$ range to see which such characters $\mu$ can arise, if any.  We will show that our conjectured formula agrees with a supercuspidal character in the $n(w) = 0$ range, for a unique $\mu$.  We will also show that $\mu = \Delta_{\chi}^{-1}$, and so $\phi = \chi \Delta_{\chi}$.

\begin{lemma}
$$F(\tilde\chi)(w) = \epsilon(\tilde\chi, \Delta^+, w) \
(\chi(w)\Omega(\frac{w}{\delta}) + \chi(\overline{w})\Omega(\frac{\overline{w}}{\delta})) \ \ \forall w \in E^* : n(w) = 0$$ for any character $\Omega$ of $E^*$ whose restriction to $F^*$ is $\aleph_{E/F}$ and whose order is a power of 2.
\end{lemma}

\proof
Recall that $n(w) = 0$ if and only if either

i) $w = p^n u + p^m v \delta$, where $u,v$ are both nonzero, $n > m$, or

ii) $w = p^m v \delta$, where $v$ is nonzero.

In case (i), $\frac{w}{\delta} = p^m v + p^{n} \frac{u}{\delta} = p^m v(1 + p^{n-m} \frac{u}{v \delta})$.  Note that $1 + p^{n-m} \frac{u}{v \delta} \in U_E^1$.  $\Omega$ is
trivial on $U_E^1$ by Lemma \ref{infinitesquareroot}, and thus we have that $\Omega(\frac{w}{\delta}) = \Omega(p^m v(1 + p^{n-m} \frac{u}{v \delta})) = \Omega(p^m v)
= \Omega(\frac{w - \overline{w}}{2 \delta}) = \tau_o(\frac{w - \overline{w}}{2 \delta})$.  In case (ii), it's clear that $\Omega(\frac{w}{\delta}) = \Omega(\frac{w-\overline{w}}{2 \delta}) = \tau_o(\frac{w-\overline{w}}{2 \delta})$. Therefore, we get that $$\tau_o(\frac{w-\overline{w}}{2 \delta}) = \Omega(\frac{w}{\delta}) \ \forall w \in E^* : n(w) = 0.$$

\noindent Thus, since $\tau_o(-1) = (-1, \Delta)$, $F(\tilde\chi)$ simplifies in the $n(w) = 0$ range to $$F(\tilde\chi)(w) = \epsilon(\tilde\chi, \Delta^+, w) \
\frac{\chi(w) + (-1, \Delta) \chi(\overline{w})}{ \tau_o(\frac{w-\overline{w}}{2
\delta})} = $$ $$\epsilon(\tilde\chi, \Delta^+, w) \
(\chi(w)\Omega(\frac{w}{\delta}) + \chi(\overline{w})\Omega(\frac{\overline{w}}{\delta})) \ \ \forall w \in E^* : n(w) = 0$$
\qed

\begin{lemma}
$F(\tilde\chi)(w) = \theta_{\pi}(w) \ \forall w \in E^* : n(w) = 0$, for $\phi = \chi \mu^{-1}$, for some character $\mu$ of $E^*$ whose restriction to $F^*$ is $\aleph_{E/F}$ and whose order is a power of 2.  In particular, $F(\tilde\chi)(w) = \theta_{\pi}(w) \ \forall w \in E^* : n(w) = 0$ for $\phi = \chi \Delta_{\chi}$.
\end{lemma}

\proof
Since $E/F$ is ramified, $\lambda(\sigma) = 1$.  Unwinding the definitions, one can see that to prove the lemma, it suffices to show that $\Omega(\delta) = (x_{\chi}, \Delta)_F
\gamma(\Delta, \psi) \Omega(w) \mu(w)$.  If we let $\Omega = \mu^{-1}$, then we are reduced to showing that there is some character $\mu$ of $E^*$ whose restriction to $F^*$ was $\aleph_{E/F}$, and whose order is a power of 2, such that $\mu^{-1}(\delta) = (x_{\chi}, \Delta)_F
\gamma(\Delta, \psi)$.  We claim that $\mu := \Delta_{\chi}^{-1}$ is such a character.

So we want to show that $\Delta_{\chi}(\delta) = (x_{\chi},
\Delta)_F \gamma(\Delta, \psi)$.  We need to investigate the term $\alpha(\phi)$.  Note that $\alpha(\phi) = \alpha(\chi)$ since $\mu|_{1 + \mathfrak{p}_E} \equiv 1$. We prefer to work with $\alpha(\chi)$.  Firstly, since $E/F$ is ramified, $\chi$ has odd
level $n = 2m+1$ \cite[Chapter 19]{bushnellhenniart}.  We also have $\alpha(\chi) \in \mathfrak{p}_E^{-n}$ (cf Section 4.3).  So let $\alpha(\chi) = p^k u + p^{\ell} v
\delta$.  The fact that $n$ is odd clearly implies that $\ell < k$.  Therefore, rewriting $\alpha(\chi)$ as $\alpha(\chi) =
p^{\ell} \sqrt{p}(v + p^{k-\ell-1/2} u)$, we get
$\aleph_{E/F}(\zeta(\alpha(\chi), \varpi)) = (v_0, \Delta)$ where $v_0$ is the leading term of the power series expansion of $v$. But $(v_0, \Delta) = (v, \Delta)$.  Moreover, $\alpha(\chi) \in \mathfrak{p}_E^{-n}$ implies that $-n = 2 \ell + 1$.  Then, by
definition of $x_{\chi}$, we get that $x_{\chi} = p^{\ell} v$.  Now, by definition of $\Delta_{\chi}$, we have $\Delta_{\chi}(\delta) = (v, \Delta) \lambda_{E/F}(\psi)^n$  and  we wish to show that this is equal to
$(x_{\chi}, \Delta) \gamma(\Delta, \psi) = (p^{\ell} v, \Delta) \gamma(\Delta,
\psi)$.  Cancelling out terms, we want to show that $(p^{\ell}, \Delta)
\gamma(\Delta, \psi) = \lambda_{E/F}(\psi)^{-2 \ell-1}$ since $-n = 2 \ell +
1$.  Well, we know that $$\lambda_{E/F}(\psi)^{-2 \ell} = (-1)^{\ell}$$  (cf \cite[page 217]{bushnellhenniart})  Thus we are reduced to showing that $\lambda_{E/F}(\psi) =
\gamma(\Delta, \psi)^{-1}$ since $(p, \Delta) = -1$.  But we prove this in Section (\ref{variousneededconstants}).
\qed

Therefore, we have proven the following, when $E/F$ is ramified.

\begin{theorem}
$F(\tilde\chi)$ agrees with the character of the supercuspidal representation $\pi_{\chi \Delta_{\chi}}$ on the range $\{w \in E^* : 0 \leq n(w) \leq r/2 \}$.
\end{theorem}

What we have actually proven is that if $(E/F, \chi)$ is a \emph{minimal} regular pair with $E/F$ ramified and $\chi$ having positive level, then $F(\tilde\chi)$ agrees with the character of the supercuspidal representation $\pi_{\chi \Delta_{\chi}}$ on the range $\{w \in E^* : 0 \leq n(w) \leq r/2 \}$.  To prove this for an arbitrary regular pair follows from this.  For if $(E/F, \chi)$ is an arbitrary regular pair, then there exists a minimal regular pair $(E/F, \chi')$ such that $\chi = \chi' \phi_E$ where $\phi_E = \phi \circ N_{E/F}$ for some $\phi \in \widehat{F^*}$.  Moreover, $\pi_{\chi} = \phi \pi_{\chi'}$ by definition.  We proved above that $F(\widetilde{\chi'}) = \theta_{\pi_{\chi' \Delta_{\chi'}}}$ on the range $\{w \in E^* : 0 \leq n(w) \leq r/2 \}$.  Therefore, $\theta_{\pi_{\chi \Delta_{\chi}}}(w) = \theta_{\pi_{\chi' \phi_E \Delta_{\chi' \phi_E}}}(w) = \theta_{\pi_{\chi' \phi_E \Delta_{\chi'}}}(w) = \phi_E(w) \theta_{\pi_{\chi' \Delta_{\chi'}}}(w) = \phi_E(w) F(\widetilde{\chi'})(w) = F(\widetilde{\chi'\phi_E})(w) = F(\tilde\chi)(w)$ on the range $\{w \in E^* : 0 \leq n(w) \leq r/2 \}$.

\

Now we consider the case $E/F$ is unramified, so $\delta = \sqrt{\Delta}$, where $\Delta \in \mathfrak{o}_F^*$ is not a square. Note that $(-1, \Delta) = 1$ since $E/F$ is unramified.

We again first conduct a careful analysis of the supercuspidal characters evaluated on the range $\{w \in E^* : 0 < n(w) \leq r/2 \}$.  Recall that on the regular set, $n(w) > 0$ if and only if $w = p^n u + p^m v \delta, \ u,v \in \mathfrak{o}_F^*$ where $n < m$.

\begin{proposition}
$F(\tilde\chi)$ agrees with the character of the supercuspidal representation $\pi_{\chi \mu^{-1}}$ in the $0 < n(w) \leq r/2$ range, where $\mu$ is any character of $E^*$ whose restriction to $F^*$ is $\aleph_{E/F}$ and whose order is a power of 2.
\end{proposition}

\proof
Again we claim that $$\gamma(\alpha(\chi), Y) = (x_{\chi}, \Delta)_F \gamma(\Delta, \psi)
\mu(\frac{w - \overline{w}}{2 \delta}) \mu(w)$$  $$\gamma(\alpha(\chi), {}^s Y) = (x_{\chi}, \Delta)_F \gamma(\Delta, \psi)
\mu(\frac{\overline{w} - w}{2 \delta}) \mu(\overline{w})$$ where $1 \neq s \in W = Aut(E/F)$ and $\alpha(\chi) = a + x_{\chi} \delta$, for $a, x_{\chi} \in F$.  The reasoning is similar as in Lemma (\ref{gammafactors}).  One must prove that $\mu(\frac{w + \overline{w}}{2}) = \mu(w) \ \forall w \in E^* : n(w) > 0$.  This is elementary.

Thus, the supercuspidal character on the $0 < n(w) \leq r/2$ range simplifies to $$\theta_{\pi}(w) = \epsilon(\tilde\chi, \Delta^+, w) \ \frac{\phi(w) \mu(w) + (-1, \Delta) \phi(\overline{w}) \mu(\overline{w})}{\mu(\frac{w-\overline{w}}{2 \delta})} \ \ \forall w \in E^* : 0 < n(w) \leq r/2.$$
\qed

The above analysis shows $F(\tilde\chi)(w) = \theta_{\pi}(w) \ \forall w \in E^* : 0 < n(w) < r/2$ where $\phi = \chi \mu^{-1}$ for any character $\mu$ of $E^*$ whose restriction to $F^*$ is $\aleph_{E/F}$ and whose order is a power of 2.

We will need to investigate the $n(w) = 0$ range to see which such characters $\mu$ can arise, if any.  We will show that our conjectured formula agrees with a supercuspidal character in the $n(w) = 0$ range, for a unique $\mu$.  We will also show that $\mu = \Delta_{\chi}^{-1}$, and so $\phi = \chi \Delta_{\chi}$.

\begin{lemma}
$$F(\tilde\chi)(w) = \epsilon(\tilde\chi, \Delta^+, w) \
(\chi(w)\Omega(\frac{w}{\delta}) + \chi(\overline{w})\Omega(\frac{\overline{w}}{\delta})) \ \forall w \in E^* : n(w) = 0$$ for the unique unramified character $\Omega$ of $E^*$ whose restriction to $F^*$ is $\aleph_{E/F}$.
\end{lemma}

\proof
Recall that $n(w) = 0$ if and only if either

i) $w = p^n u + p^m v \delta$, where $u,v$ are both nonzero, $n > m$, or

ii) $w = p^m v \delta$, where $v$ is nonzero.

iii) $w = p^n u + p^m v \delta$, where $n = m$ and $u,v$ are both non-zero.

We first show that $\tau_o(\frac{w-\overline{w}}{2 \delta}) = \Omega(\frac{w}{\delta}) \ \forall w \in E^* : n(w) = 0$. In case (i), $\frac{w}{\delta} = p^m v + p^n \frac{u}{\delta} = p^m v(1 + p^{n-m} \frac{u}{v \delta})$.  Note that $1 + p^{n-m} \frac{u}{v \delta} \in U_E^1$.  But $\Omega$ is trivial on $U_E^1$ by Lemma (\ref{infinitesquareroot}).  Therefore, we have $\Omega(\frac{w}{\delta}) = \Omega(p^m v(1 + p^{n-m} \frac{u}{v \delta})) = \Omega(p^m v) = \Omega(\frac{w - \overline{w}}{2 \delta}) = \tau_o(\frac{w - \overline{w}}{2 \delta})$.  In case (ii), it's clear that $\Omega(\frac{w}{\delta}) = \Omega(\frac{w-\overline{w}}{2 \delta}) = \tau_o(\frac{w-\overline{w}}{2 \delta})$.  In case (iii), $w = p^n (u + v \delta)$.  But $u + v \delta \in \mathfrak{o}_E^*$, and therefore $\Omega(w) = \Omega(p^n)\Omega(u+v \delta) = \Omega(p^n) = \Omega(p^n v)$ since $\Omega$ is unramified.  In particular, since $\Omega(\delta) = 1$, we get $\Omega(\frac{w}{\delta}) = \tau_o(\frac{w-\overline{w}}{2 \delta})$.  Therefore, $$\tau_o(\frac{w-\overline{w}}{2 \delta}) = \Omega(\frac{w}{\delta}) \ \forall w \in E^* : n(w) = 0.$$

Since $\tau_o(-1) = (-1, \Delta) = 1$, $F(\tilde\chi)$ simplifies in the $n(w) = 0$ range. $$F(\tilde\chi)(z) = \epsilon(\tilde\chi, \Delta^+, w) \
\frac{\chi(w) + \chi(\overline{w})}{ \tau_o(\frac{w-\overline{w}}{2
\delta})} = \epsilon(\tilde\chi, \Delta^+, w) \ (\chi(w)\tau_o(\frac{w-\overline{w}}{2
\delta}) + \chi(\overline{w})\tau_o(\frac{\overline{w}-w}{2
\delta})) =$$ $$ \epsilon(\tilde\chi, \Delta^+, w) \
(\chi(w)\Omega(\frac{w}{\delta}) + \chi(\overline{w})\Omega(\frac{\overline{w}}{\delta})) \ \forall w \in E^* : n(w) = 0$$
\qed

We want to show that $F(\tilde\chi)(w) = \theta_{\pi}(w) \ \forall w \in E^* : n(w) = 0$ for $\phi = \chi \mu^{-1}$ for some character $\mu$ of $E^*$ whose restriction to $F^*$ is $\aleph_{E/F}$ and whose order is a power of 2.  Note that $\Omega(\delta) = 1$ since $\Omega$ is unramified. Unwinding the definitions, one can see that to prove $F(\tilde\chi)(w) = \theta_{\pi}(w) \ \forall w \in E^* : n(w) = 0$ for $\phi = \chi \mu^{-1}$, it suffices to show that $$(x_{\chi}, \Delta) \gamma(\Delta, \psi) \chi(w) \Omega(w) = \chi(w) \mu^{-1}(w) \lambda(\sigma) \ \forall w \in E^* : n(w) = 0$$

\begin{lemma}
$\lambda(\sigma) = (x_{\chi}, \Delta)_F \gamma(\Delta, \psi)$.
\end{lemma}

\proof
Unwinding the definition of $\lambda(\sigma)$ (cf \cite[Section 5.3]{debacker}), one can see that since $E/F$ is unramified, then $\lambda(\sigma) = (-1)^{r+1}$ where $r$ is the depth of the supercuspidal representation $\pi_{\phi}$ via Theorem (\ref{positivedepth}).  Note that the depth of $\pi_{\phi}$ equals the depth of $\pi_{\chi}$ since $\mu|_{1 + \mathfrak{p}_E} \equiv 1$, i.e. $\mu$ has level zero.

Now, consider the term $(x_{\chi}, \Delta)$.  We need to investigate the term $\alpha(\phi)$.  Note that $\alpha(\phi) = \alpha(\chi)$ since again, $\mu|_{1 + \mathfrak{p}_E} \equiv 1$. We prefer to work with $\alpha(\chi)$.  Recall that $\alpha(\chi) \in \mathfrak{p}_E^{-n}$.  Moreover, since $\alpha(\chi) \in \mathfrak{g}_{-r}' \setminus \mathfrak{g}_{-r^+}'$ (cf \cite[page 34]{debacker}), we have that $n = r$.  Now let $\alpha(\chi) = p^k u + p^l v \delta$.  We need a lemma.

\begin{lemma}
$l = -n$.
\end{lemma}

\proof
Since $(E/F, \chi)$ is a minimal regular pair (cf \cite[Section 19.2 line 1]{bushnellhenniart}), we have that $\alpha(\chi)$ is a \emph{minimal element} over $F$ (see \cite[Proposition 18.2]{bushnellhenniart}).  By \cite[Section 13.4]{bushnellhenniart}, $\alpha(\chi)$ is minimal over $F$ if and only if $(\alpha(\chi) + \mathfrak{p}_E^{-n+1}) \bigcap F = \emptyset$, where $n = -v_E(\alpha(\chi))$, where $v_E$ denotes valuation.  Recall that $\alpha(\chi) \in \mathfrak{p}_E^{-n} \setminus \mathfrak{p}_E^{-n+1}$.  It is not difficult to show that if $l \neq -n$, then $(\alpha(\chi) + \mathfrak{p}_E^{-n+1}) \bigcap F \neq \emptyset$, a contradiction, and so the lemma is proven.
\qed

Returning to the proof of the proposition, recall that $x_{\chi} = p^l v$.  Thus, $(x_{\chi}, \Delta) = (p^{-n} v, \Delta)$.  Also, since $r = n$, we have $\lambda(\sigma) = (-1)^{r+1} = (-1)^{n+1}$.  Therefore, we are reduced to showing that

$$(p^{-n} v, \Delta) \gamma(\Delta, \psi) = (-1)^{n+1}$$ so equivalently, $\gamma(\Delta, \psi) = -1$.  But since $\psi$ has level 1, it is a fact that $\gamma(\Delta, \psi) = -1$ (we will show this in Section (\ref{variousneededconstants})).
\qed

\begin{proposition}
$F(\tilde\chi)(w) = \theta_{\pi}(w) \ \forall w \in E^* : n(w) = 0$, for $\phi = \chi \mu^{-1}$, for some character $\mu$ of $E^*$ whose restriction to $F^*$ is $\aleph_{E/F}$ and whose order is a power of 2.  In particular, $F(\tilde\chi)(w) = \theta_{\pi}(w) \ \forall w \in E^* : n(w) = 0$ for $\phi = \chi \Delta_{\chi}$.
\end{proposition}

\proof
By the previous lemma, the conjectured equation $$(x_{\chi}, \Delta) \gamma(\Delta, \psi) \chi(w) \Omega(w) = \chi(w) \mu^{-1}(w) \lambda(\sigma) \ \forall w \in E^*:  n(w) = 0$$ simplifies to $\Omega(w) = \mu^{-1}(w) \ \forall w \in E^* : n(w) = 0$.  But recall that $\mu$ is a character of $E^*$ whose restriction to $F^*$ is $\aleph_{E/F}$ and whose order is a power of 2.  There is only one character $\mu$ of $E^*$ whose restriction to $F^*$ is $\aleph_{E/F}$, whose order is a power of 2, and that equals $\Omega^{-1}$ on the $n(w) = 0$ range.  Indeed, this forces $\mu^{-1} = \Delta_{\chi} = \Delta_{\chi}^{-1}$.
\qed

Therefore, we have proven the following, when $E/F$ is unramified.

\begin{theorem}\label{charactersmatchingup}
$F(\tilde\chi)$ agrees with the character of the supercuspidal representation $\pi_{\chi \Delta_{\chi}}$ on the range $\{w \in E^* : 0 \leq n(w) \leq r/2 \}$.
\end{theorem}

Again, we have assumed at various points that our regular pairs are minimal.  However, our results hold without this assumption for the same reason as in the case of $E/F$ ramified.

\subsection{On whether there are two positive depth character formulas coming from the same Cartan}\label{samecartan}

In the next two sections we show that a positive depth supercuspidal representation of $GL(2,F)$ is uniquely determined by the restriction of its distribution character to the $n(w) = 0$ range.  In this section, we show that if the distribution characters of two positive depth supercuspidal representations, both coming from the same Cartan, agree on the $n(w) = 0$ range, then the supercuspidal representations are isomorphic.  That is, we prove the following theorem.

\begin{theorem}\label{samecartans}
Suppose $(E/F, \chi_1)$ and $(E/F,\chi_2)$ are admissible pairs such that $F(\tilde\chi_1)(w) = F(\tilde\chi_2)(w) \ \forall w \in E^* : n(w) = 0$.  Then, $\chi_1 = \chi_2^{\upsilon}$ for some $\upsilon \in Aut(E/F)$.
\end{theorem}

We wish to make the following important note.  Recall that in the previous sections, we constructed a character formula $F(\tilde\chi)$ from a regular pair $(E/F, \chi)$.  The above theorem and Theorem (\ref{differentcartans}) will together prove that a positive depth supercuspidal representation of $GL(2,F)$ is uniquely determined by its restriction of its distribution character to the $n(w) = 0$ range.  We are claiming in the above theorem and in Theorem (\ref{differentcartans}) that it is sufficient to consider admissible pairs rather than regular pairs in order to prove that a positive depth supercuspidal representation of $GL(2,F)$ is uniquely determined by its restriction of its distribution character to the $n(w) = 0$ range.  This is because the positive depth supercuspidal representations of $GL(2,F)$ are parameterized by admissible pairs, and so it is sufficient to consider just admissible pairs.

\begin{lemma}\label{samecartanramifiedgl2}
Let $E/F$ be ramified.  Suppose $(E/F, \chi_1)$ and $(E/F, \chi_2)$ are admissible pairs such that $F(\tilde\chi_1)(w) = F(\tilde\chi_2)(w) \ \forall w \in E^* : n(w) = 0$.  Then, $\chi_1 = \chi_2^{\upsilon}$ for some $\upsilon \in Aut(E/F)$.
\end{lemma}

\proof
We may assume without loss of generality that $E = F(\sqrt{p})$.  We have assumed that $$\epsilon(\tilde\chi_1, \Delta^+, w) \frac{\chi_1(w) + (-1, \Delta) \chi_1(\overline{w})}{\tau_o(\frac{w-\overline{w}}{2 \delta})} = \epsilon(\tilde\chi_2, \Delta^+, w) \frac{\chi_2(w) + (-1, \Delta) \chi_2(\overline{w})}{\tau_o(\frac{w-\overline{w}}{2 \delta})} \ \ \forall w \in E^* : n(w) = 0$$ Let us write $deg(\pi_i)$ for the $deg(\pi)$ that are associated to the pairs $(E/F, \chi_i)$. Now, let $c_i = deg(\pi_i)(x_{\chi_i}, \Delta) |\eta(\alpha(\chi_i))|^{-1/2}$.  Then, cancelling out like terms, we have that $$c_1 (\chi_1(w) + (-1, \Delta) \chi_1(\overline{w})) = c_2 (\chi_2(w) + (-1, \Delta) \chi_2(\overline{w})) \ \forall w \in E^* : n(w) = 0$$  Then the same proof of \cite[Lemma 5.1]{spice} shows that $\chi_1|_{F^*(1 + \mathfrak{p}_E)} = \chi_2^{\upsilon}|_{F^*(1 + \mathfrak{p}_E)}$ for some $\upsilon \in Aut(E/F)$.  For the following arguments, it suffices without loss of generality to assume $\upsilon = 1$.

Let $c := \frac{c_1}{c_2}$.  Let $[\chi](w) := \chi(w) + (-1, \Delta) \chi(\overline{w})$.  Then we have $c [\chi_1](w) =  [\chi_2](w) \ \forall w \in E^* \setminus F^*(1 + \mathfrak{p}_E)$ and we also have that $\chi_1|_{F^*(1 + \mathfrak{p}_E)} = \chi_2|_{F^*(1 + \mathfrak{p}_E)}$.  Now let $p_E$ be a uniformizer of $E$ and recall that $p$ is a uniformizer of $F$.  We may take $p_E$ so that $p_E^{2} = p$. Since $\chi_1(p) = \chi_2(p)$, we have that $\chi_1(p_E)^{2} = \chi_2(p_E)^{2}$, and so $\chi_2(p_E) = \xi_{2} \chi_1(p_E)$ where $\xi_2$ could be plus or minus 1.  Therefore, $\chi_2(w) = \chi_1(w) \xi_{2}^{val(w)} \ \forall w \in E^*$, where $val(w)$ denotes the $E$-adic valuation of $w \in E^*$.  Therefore, after substituting and noting that $val(w) = val(\upsilon(w))$, we obtain $$c [\chi_1](w) = \xi_{2}^{val(w)} [\chi_1](w) \ \forall w \in E^* \setminus F^*(1 + \mathfrak{p}_E)$$  We will prove in the next section that there exists a $w' \in E^* \setminus F^*(1 + \mathfrak{p}_E)$ such that $[\chi_1](w') \neq 0$. Therefore, we can cancel $[\chi_1](w')$ from both sides to obtain $$c = \xi_{2}^{val(w')}$$  Therefore, $c$ is plus or minus $1$.

Suppose that $c = 1$.  Then, we get $[\chi_1](w) =  [\chi_2](w) \ \forall w \in E^*$. By linear independence of characters, $\chi_1 = \chi_2$.

Suppose $c = -1$.  Then $\chi_2(w) = \chi_1(w) (-1)^{val(w)}$. This implies that $\chi_2 = \chi_1 \otimes \phi_E$, where $\phi_E := \phi \circ N_{E/F}$ where $\phi = \aleph_{L/F}$ where $L/F$ is the unique unramified degree 2 extension of $F$.  In this case, one can check that $\alpha(\chi_1) = \alpha(\chi_2)$.  Therefore, $c = \frac{deg(\pi_1)}{deg(\pi_2)}$.  But formal degrees are positive real numbers, and so we get a contradiction to the supposition that $c = -1$.

Therefore, $\chi_1 = \chi_2$ or $\chi_1 = \chi_2^{\upsilon}$, and so the admissible pairs are isomorphic.
\qed

\begin{lemma}
Let $E/F$ be unramified.  Suppose $(E/F, \chi_1)$ and $(E/F, \chi_2)$ are admissible pairs such that $F(\tilde\chi_1)(w) = F(\tilde\chi_2)(w) \ \forall w \in E^* : n(w) = 0$.  Then, $\chi_1 = \chi_2^{\upsilon}$ for some $\upsilon \in Aut(E/F)$.
\end{lemma}

\proof
See \cite[page 16]{spice}.
\qed

\subsection{On whether there are two positive depth character formulas coming from different Cartans}\label{othercartan}

In this section we show that the distribution characters of two positive depth supercuspidal representations, coming from different Cartans, can't agree on the $n(w) = 0$ range.  This, together with the results from the previous section, shows that if $(E/F, \chi)$ is an admissible pair, then there is a unique positive depth supercuspidal representation whose character agrees with $F(\tilde\chi)$ on the range $\{ w \in E^* : n(w) = 0$ \}.

\begin{theorem}\label{differentcartans}
Suppose $(E/F, \chi)$ and $(E_1/F, \chi_1)$ are admissible pairs with $E \ncong E_1$.  Then $\exists w \in E^* : n(w) = 0$ such that $F(\tilde\chi)(w) \neq \theta_{\pi_{\chi_1 \Delta_{\chi_1}}}(w)$.
\end{theorem}

Let $E = F(\sqrt{\Delta_E})$ and $E_1 = F(\sqrt{\Delta_{E_1}})$.  There are many cases to check, and we split them up in a sequence of propositions.

\begin{proposition}\label{case1}
Suppose $(E/F, \chi)$ and $(E_1/F, \chi_1)$ are admissible pairs with $E$ ramified and $E_1$ unramified.  Then $\exists w \in E^* : n(w) = 0$ such that $F(\tilde\chi)(w) \neq \theta_{\pi_{\chi_1 \Delta_{\chi_1}}}(w)$.
\end{proposition}

\proof
By comparing valuations of determinants of elements, one can show that since the inducing representation of the representation coming from $E_1$ is $E_1^* G_{x,r/2}$ (see \cite[page 34-35]{debacker}), then if $w \in E^* : n(w) = 0$, then one can't conjugate $w$ into $E_1^* G_{x,r/2}$ (note that if two elements are conjugate, then they have the same determinant).  Therefore, $\theta_{\pi_{\chi_1 \Delta_{\chi_1}}}(w) = 0$.  We now need to consider two subcases.

Subcase (a): Suppose $(-1, \Delta_E) = -1$.  We suppose by way of contradiction that $F(\tilde\chi)(w) = 0 \ \forall w \in E^* : n(w) = 0$.  This implies that $\chi(w) + (-1, \Delta_E) \chi(\overline{w}) = 0 \ \forall w \in E^* : n(w) = 0$.  Since $(-1, \Delta_E) = -1$, we have $\chi(w) = \chi(\overline{w}) \ \forall w \in E^* : n(w) = 0$.  But since $E^*$ is generated as a group by the set $\{ w \in E^* : n(w) = 0 \}$, we get $\chi(w) = \chi(\overline{w}) \ \forall w \in E^*$, which contradicts the fact that $(E/F, \chi)$ is an admissible pair.

Subcase (b): Suppose $(-1, \Delta_E) = 1$, and suppose again by way of contradiction that $F(\tilde\chi)(w) = 0 \ \forall w \in E^* : n(w) = 0$.  Then $\chi(w) + \chi(\overline{w}) = 0 \ \forall w \in E^* : n(w) = 0.$  Thus, $\chi(w/\overline{w}) = -1 \ \forall w \in E^* : n(w) = 0$, so $\chi(w/\overline{w})^2 = 1$ for all $w \in E^* : n(w) = 0$.  But since the set $\{w \in E^* : n(w) = 0 \}$ generates all of $E^*$ as a group, we get that $\chi|_{(E^*)^2} = \chi^{\upsilon}|_{(E^*)^2}$.  Now since $U_E^1 \subset (E^*)^2$ by Lemma (\ref{serresquare}), we get $\chi|_{U_E^1} = \chi^{\upsilon}|_{U_E^1}$ , which contradicts the fact that $(E/F, \chi)$ is an admissible pair.
\qed

\begin{proposition}\label{case2}
Suppose $(E/F, \chi)$ and $(E_1/F, \chi_1)$ are admissible pairs with $E \ncong E_1$, with $E$ unramified and $E_1$ ramified.  Then $\exists w \in E^* : n(w) = 0$ such that $F(\tilde\chi)(w) \neq \theta_{\pi_{\chi_1 \Delta_{\chi_1}}}(w)$.
\end{proposition}

We will need to split this proposition into two cases: $(-1,p) = 1$ and $(-1,p) = -1$.  We have $E = F(\sqrt{\Delta})$, where $\Delta \in \mathfrak{o}_F^*$ is not a square, and without loss of generality $E_1 = F(\sqrt{p})$.

\begin{lemma}\label{subcasea''}
Suppose $(-1,p) = 1$.  Suppose $(E/F, \chi)$ and $(E_1/F, \chi_1)$ are admissible pairs with $E \ncong E_1$, with $E$ unramified and $E_1$ ramified.  Then $\exists w \in E^* : n(w) = 0$ such that $F(\tilde\chi)(w) \neq \theta_{\pi_{\chi_1 \Delta_{\chi_1}}}(w)$.
\end{lemma}

\proof
We use the same strategy as in Proposition (\ref{case1}).  Recall that $w \in E^* : n(w) = 0$ if and only if
(i) $w = p^m v \delta$,
(ii) $w = p^n u + p^m v \delta, \ n = m$, or
(iii) $w = p^n u + p^m v \delta, \ n > m$

By comparing determinants again, one can show that if $w$ is of the form (i) or (iii), then $w$ can't be conjugated into $E_1^* G_{x,r/2}$ (see \cite[page 34-35]{debacker}).  Therefore, $\theta_{\pi_{\chi_1 \Delta_{\chi_1}}}(w) = 0$ for all $w$ of the form (i) and (iii).  Now assume that $F(\tilde\chi)(w) = 0$ for all $w$ of the form (i) and (iii).  Therefore, since $(-1, \Delta) = 1$, $\chi(w) + \chi(\overline{w}) = 0$ for all $w$ of the form (i) and (iii).  We need the following lemma.

\begin{lemma}\label{randomlemma1}
$\chi|_{F^* U_E^1} = \chi^{\upsilon}|_{F^* U_E^1}$, where $\upsilon$ generates $Aut(E/F)$.
\end{lemma}

\proof
Let $z \in U_E^1$.  It is easy to see that one can write $z = w_1 \overline{w_2}$, where $w_1$ is of the form (i) and $w_2$ is of the form (iii).  Now, since $F(\tilde\chi)(w) = 0$ on all $w$ of the form (i) and (iii), we have $$\chi(w_1) + \chi(\overline{w_1}) = 0$$ $$\chi(w_2) + \chi(\overline{w_2}) = 0$$  Multiplying the first equation by $\chi(w_2)$ and the second equation by $\chi(w_1)$, we conclude that $$\chi(w_2 \overline{w_1}) = \chi(w_1 \overline{w_2})$$  Therefore, $\chi(z) = \chi(\overline{z}) \ \forall z \in U_E^1$.  Since $\chi(x) = \chi(\overline{x}) \ \forall x \in F^*$, we get $\chi|_{F^* U_E^1} = \chi^{\upsilon}|_{F^* U_E^1}$.
\qed

We need the following lemma, which we shall not prove, as it is not difficult.

\begin{lemma}\label{intermediatelemma}
Every element b of the form (ii) can be written as a product ac, where a is either an element of the form (i) or (iii), and c is an element of $F^* U_E^1$.
\end{lemma}

Because of lemma (\ref{intermediatelemma}), we have that if $b$ is of the form (ii), then if we write $b = ac$ as in Lemma (\ref{intermediatelemma}), then $\chi(b) + \chi(\overline{b}) = \chi(ac) + \chi(\overline{ac}) = \chi(a) \chi(c) + \chi(\overline{a}) \chi(\overline{c}) = \chi(a) \chi(c) + \chi(\overline{a}) \chi(c)$ by Lemma (\ref{randomlemma1}).  Thus, $\chi(b) + \chi(\overline{b}) = \chi(c)(\chi(a) + \chi(\overline{a}))$.  But we assumed that $\chi(a) + \chi(\overline{a}) = 0$ for all elements $a$ of the form (i) and of the form (iii). Thus, $F(\tilde\chi)(b) = 0$ for all $w$ in of the form (ii).  Thus, $F(\tilde\chi)(w) = 0 \ \forall w \in E^* : n(w) = 0$.

Recall that we also have that $\chi|_{F^* U_E^1} = \chi^{\upsilon}|_{F^* U_E^1}$ where $1 \neq \upsilon \in Aut(E/F)$.  By a similar argument as in the proof of Lemma (\ref{randomlemma1}), we have that since $F(\tilde\chi)(w) = 0 \ \forall w \in E^* : n(w) = 0$, then $\chi(w_1 \overline{w_2}) = \chi(\overline{w_1} w_2) \ \forall w_1, w_2 \in E^* : n(w_1) = n(w_2) = 0$.  Now, if $w_1$ is in of the form (i) and $w_2$ is in of the form (ii), then $w_1 \overline{w_2}$ is of the form (ii).  This shows that $\chi(z) = \chi(\overline{z})$ for some element $z$ of the form (ii).  But notice that the group generated by $U_E^1, F^*$, and any single element in of the form (ii) is all of $E^*$.  Thus, we get that $\chi = \chi^{\upsilon}$ on all of $E^*$, and therefore $(E/F, \chi)$ is not an admissible pair. Finally, we are done with proving Lemma (\ref{subcasea''}).
\qed

\begin{lemma}\label{subcaseb''}
Suppose $(-1,p) = -1$.  Suppose $(E/F, \chi)$ and $(E_1/F, \chi_1)$ are admissible pairs with $E$ unramified and $E_1$ ramified.  Then $\exists w \in E^* : n(w) = 0$ such that $F(\tilde\chi)(w) \neq \theta_{\pi_{\chi_1 \Delta_{\chi_1}}}(w)$.
\end{lemma}

\proof

Recall again that $w \in E^* : n(w) = 0$ if and only if

(i) $w = p^m v \delta$

(ii) $w = p^n u + p^m v \delta, \ n = m$, or

(iii) $w = p^n u + p^m v \delta, \ n > m$.

We want to show like in Lemma (\ref{subcasea''}) that any character formula $\theta_{\pi_{\chi_1 \Delta_{\chi_1}}}$ vanishes on elements of the form (i) and (iii).  After showing this, rest of the proof of Lemma (\ref{subcaseb''}) goes exactly the same way as in Lemma (\ref{subcasea''}).  Well, the values on the ramified torus $E_1$ of the character formula of a supercuspidal representation coming from an admissible pair $(E_1/F, \chi_1)$ can be computed using \cite[Proposition 2, Proposition (\ref{case1}), $\ell > 0$, p. 101]{shimizu}.  This formula shows that $\theta_{\pi_{\chi_1 \Delta_{\chi_1}}}(w) = 0$ for all $w \in E^*$ such that $w$ is of the form (i) and (iii).
\qed

Therefore, we have finished the proof of Proposition (\ref{case2}).

\begin{proposition}\label{case3}
Suppose $(E/F, \chi)$ and $(E_1/F, \chi_1)$ are admissible pairs with $E \ncong E_1$, with $E$ ramified and $E_1$ ramified.  Then $\exists w \in E^* : n(w) = 0$ such that $F(\tilde\chi)(w) \neq \theta_{\pi_{\chi_1 \Delta_{\chi_1}}}(w)$.
\end{proposition}

\proof
Suppose without loss of generality that $E = F(\sqrt{p})$ and $E_1 = F(\sqrt{\Delta p})$, $\Delta \in \mathfrak{o}_F^*$ not a square.  The proof for the case where $E = F(\sqrt{\Delta p})$ and $E = F(\sqrt{p})$ is similar.  We first claim that elements $w \in E^*$ such that $n(w) = 0$ can't be conjugated into $E_1^* G_{x,r/2}$ (see \cite[page 34-35]{debacker}).

Recall that $w \in E^*$ such that $n(w) = 0$ if and only if either $w = p^n u + p^m v \sqrt{p}$ with $n > m$, $u,v \neq 0$, or $w = p^m v \sqrt{p}$ with $v \neq 0$.  Suppose $w = p^m v \sqrt{p}$. Then $det(w) = N(w) = - p^{2m+1} v^2$.  Recall that $det(E_1^* G_{x, r/2}) \subset N_{E_1 / F}(E_1^*)$.  Well, $(det(w), \Delta p) = (-p^{2m+1}v^2, \Delta p) = (-p, \Delta p) = (-1, \Delta p) (p, \Delta p) = (\Delta p, \Delta p) (p, \Delta p) = (\Delta p^2, \Delta p) = (\Delta, \Delta p) = (\Delta, \Delta) (\Delta, p) = (\Delta, p) = -1$.  Therefore, det($w$) can't be a norm from $E_1$.
Suppose $w = p^n u + p^m v \sqrt{p}$ with $n > m, v \neq 0$.  Then the same method of proof works as above to show that $w$ can't be conjugated into $E_1^* G_{x,r/2}$, but the calculation is more tedious.  We omit the details.  The rest of the proof follows as in subcases (a) and (b) from Proposition (\ref{case1}) above.
\qed

We have now finished the proof of Theorem (\ref{differentcartans}).  Summing up, we have altogether shown that if $(E/F, \chi)$ is a regular pair such that $\chi$ has positive level, then there is a unique positive depth supercuspidal representation, $\pi_{\chi \Delta_{\chi}}$, whose character, on the range $\{ z \in T(F)^{reg} : 0 \leq n(z) \leq r/2 \}$, agrees with $F(\tilde\chi)$.  There is one minor point here to resolve.  Is there possibly a depth zero supercuspidal representation whose character, on the range $\{ z \in T(F)^{reg} : 0 \leq n(z) \leq r/2 \}$, also equals $F(\tilde\chi)$?  We will prove in the next chapter that if $(E_1/F, \chi_1)$ is a regular pair corresponding to a depth zero supercuspidal representation $\pi$, then its character formula, on the range $\{ z \in T(F)^{reg} : 0 \leq n(z) \leq r/2 \}$, is $$F(\tilde\chi_1)(w) = - \epsilon(\Delta^+) \frac{deg(\pi)}{deg(\sigma)} \left(\frac{\chi_1(w) + (-1, \Delta) \chi_1(\overline{w})}{ \tau_o(\frac{w - \overline{w}}{2 \delta})}\right), \ \ w \in E_1^* \setminus F^*(1 + \mathfrak{p}_{E_1})$$ Then, the same arguments as in Theorems (\ref{samecartans}) and (\ref{differentcartans}) show that the character of $\pi$ cannot equal $F(\tilde\chi)$, on the range $\{ z \in T(F)^{reg} : 0 \leq n(z) \leq r/2 \}$, unless $\pi \cong \pi(\tilde\chi)$.

Therefore, combining Theorems (\ref{samecartans}), (\ref{differentcartans}), and (\ref{charactersmatchingup}), we obtain the following result.

\begin{theorem}
The assignment

\begin{eqnarray}
\left\{
\begin{array}{rl}
 irreducible \ \varphi : W_F \rightarrow GL(2,\mathbb{C})
\end{array} \right\} & \mapsto & \tilde\chi \in \widehat{T(F)}_{\tau \circ \rho} \mapsto \pi(\tilde\chi) \nonumber
\end{eqnarray}

\noindent from Section (\ref{setup}) is the Local Langlands correspondence for positive depth supercuspidal representations of $GL(2,F)$, where $\pi(\tilde\chi)$ is the unique supercuspidal representation whose character, on the range $\{ z \in T(F)^{reg} : 0 \leq n(z) \leq r/2 \}$, is $F(\tilde\chi)$.
\end{theorem}

\subsection{Calculation of various needed constants}\label{variousneededconstants}

In this section we calculate various constants needed throughout Chapter 5.  Let $\Delta \in \mathfrak{o}_F^*$ be a non-square, and let $\psi$ be an additive character of $F$.  Let $\ell(\psi)$ denote the level of $\psi$.  One can prove the following Lemma using the machinery from \cite[Appendix]{rao}.

\begin{lemma}
$\gamma_F(\Delta, \psi) = (-1)^{\ell(\psi)}$.
\end{lemma}

\begin{lemma}
$\lambda_{E/F}(\psi) = \gamma_F(\Delta, \psi) (-1, \Delta)_F$
\end{lemma}

\proof
Recall that $\lambda_{E/F}(\psi)$ is the Langlands constant (cf \cite[Pages 216, 217, 241]{bushnellhenniart}), where $\psi$ is an additive character of $F$.  By \cite[page 240-241]{bushnellhenniart}, $\lambda_{E/F}$ is defined to be the Weil index (cf Definition (\ref{rao0})) of $$q : E^* \rightarrow \mathbb{C}^*$$ where $q(z) := \psi(N(z))$.  The quadratic form in question is $N : E^* \rightarrow F^*$.  The associated symmetric bilinear form associated is $$(z,w)_q := \frac{Tr(z \overline{w})}{2}$$  We wish to therefore calculate the Weil index $\gamma(\psi \circ N)$.

Well, a basis for $E/F$ is $1, \delta$, where $E = F(\delta), \delta = \sqrt{\Delta}$.  Then, $(1,1)_q= 1, (1, \delta)_q = (\delta, 1)_q = 0,$ and $(\delta, \delta)_q = -\Delta$.  Thus, the matrix of the quadratic form is $$\mat{1}{0}{0}{-\Delta}$$  Then, by Lemmas (\ref{rao1}), (\ref{rao2}), and (\ref{rao3}), and since $E/F$ is quadratic, $$\gamma(\psi \circ N) = (1, -\Delta)_F \gamma_F(\psi)^2 \gamma_F(-\Delta, \psi) = \gamma_F(-1, \psi) (-1, -1)_F \gamma_F(-\Delta, \psi) =$$ $$ \gamma_F(\Delta, \psi) (-1, -\Delta)_F = \gamma_F(\Delta, \psi) (-1, \Delta)_F$$
\qed

\section{Depth zero supercuspidal character formulas for $PGL(2,F)$ and $GL(2,F)$}

\subsection{On the proof that our conjectural character formulas agree with depth zero supercuspidal characters}

In the following two sections, we prove Theorems (\ref{gl2theorem1}) and (\ref{gl2theorem2}) for the case of depth zero supercuspidal representations of $GL(2,F)$.

Let us recall from the previous chapter that the proposed character formula simplifies to $$F(\tilde\chi)(w) = \epsilon(\tilde\chi, \Delta^+, \tau)
 \left(\frac{\chi(w) + (-1, \Delta) \chi(\overline{w})}{ \tau_o(w - \overline{w}) |D(w)|^{1/2}}\right), \ w \in T(F)^{reg}$$
For depth zero representations we define $\epsilon(\tilde\chi, \Delta^+, \tau) := -\frac{deg(\pi) \tau_o(2 \delta)}{deg(\sigma)} \epsilon(\Delta^+)$, where $\epsilon(\Delta^+)$ is as in Section (\ref{theconstantepsilon}) and $deg(\pi), deg(\sigma)$ are as in Theorem (\ref{depthzerocharacters}).  Let us recall the following theorem from \cite{debacker}.

\begin{theorem}\label{depthzerocharacters}\footnote{Notice that this theorem is slightly different than the one from \cite{debacker}.  It is because there are a few typos in \cite{debacker}.}
\cite[Theorem 5.4.1]{debacker}

Let $(E/F, \chi)$ be a regular pair where $E/F$ has degree $n$, $\chi$ has level zero, and $n$ is prime.  Let $\pi = \pi_{\chi}$ be the associated depth zero supercuspidal representation of $GL(n,F)$ given by Proposition (\ref{depthzerogl2}).  Suppose $\gamma \in F^* K_0^{reg}$.  Then

\begin{equation*}
\frac{\theta_{\pi}(\gamma)}{deg(\pi)} = \left\{
\begin{array}{rl}
\chi_{\pi}(z) \frac{\chi_{\sigma}(\gamma)}{deg(\sigma)} & \text{if } \gamma = zw \ is \ unramified \ elliptic \ and \ \gamma \ is \ not \ in \ F^* K_1, z \in Z, w \in K_0\\
\chi_{\pi}(z) L.C.E. & \text{if } \gamma = z(1 + {}^g X) \ with \ X \in \mathfrak{b}_1, g \in G, \ and \ z \in Z\\
0 & \text{otherwise}
\end{array} \right.
\end{equation*}

\end{theorem}

Let $(E/F, \chi)$ be a regular pair such that $\chi$ has level zero and $E/F$ is degree $2$.  Thus, $E/F$ is unramified and $\chi|_{U_E}$ gives rise to a character $\theta$ of the multiplicative group
of the residue field $\mathbb{F}_{q^2}$ of $E$.  Note that when $E/F$ is unramified, $(E/F, \chi)$ is regular if and only if $(E/F, \chi)$ is admissible.  Let $\mathbb{G} := GL(2, \overline{\mathbb{F}_q})$.  Let
$\mathbb{T}$ be the maximal torus of $\mathbb{G}$ defined over $\mathbb{F}_q$ such that $\mathbb{T}^{\Phi} = \mathbb{F}_{q^2}^*$ is the elliptic
torus in $GL(2, \mathbb{F}_q)$.  Then, by Deligne-Lusztig theory, the pair $(\mathbb{T}, \theta)$
yields a generalized character $R_{\mathbb{T}, \theta}$ of $\mathbb{G}(\mathbb{F}_q) = GL(2, \mathbb{F}_q)$.

\begin{proposition}\label{carterformula}
If $s \in \mathbb{G}^{\Phi}$ is semisimple, then $$R_{\mathbb{T}, \theta}(s) = \frac{\epsilon_{\mathbb{T}} \epsilon_{C^0(s)}}{|\mathbb{T}^{\Phi}| |C^0(s)^{\Phi}|_p} \sum_{g \in \mathbb{G}^{\Phi} : g^{-1}sg \in \mathbb{T}^{\Phi}} \theta(g^{-1} s g)$$
\end{proposition}

\proof
\cite[Proposition 7.5.3]{carter}
\qed

\begin{proposition}\label{gl2characterformuladepthzerolala}
$$R_{\mathbb{T}, \theta}(s) = - \displaystyle\sum_{i=0}^1 \theta(\upsilon^i(s))$$ for all regular semisimple s in $\mathbb{T}^{\Phi}$, where $\upsilon$ is the generator of $Gal(\mathbb{F}_{q^2}/\mathbb{F}_q)$
\end{proposition}

\proof
Let $s \in \mathbb{T}^{\Phi}$ be regular semisimple.  Note that if $g \in \mathbb{G}^{\Phi}$ satisfies $g^{-1} s g \in \mathbb{T}^{\Phi}$, then $g \in N_{\mathbb{G}^{\Phi}}(\mathbb{T}^{\Phi})$.  Therefore, $$\sum_{g \in \mathbb{G}^{\Phi} \ : \ g^{-1}sg \in \mathbb{T}^{\Phi}} \theta(g^{-1} s g) = \sum_{g \in N_{\mathbb{G}^{\Phi}}(\mathbb{T}^{\Phi})} \theta(g^{-1} s g) = |\mathbb{T}^{\Phi}| \sum_{w \in N_{\mathbb{G}^{\Phi}}(\mathbb{T}^{\Phi}) / \mathbb{T}^{\Phi}} \theta({}^w s)$$

\noindent Therefore,

$$R_{\mathbb{T}, \theta}(s) = \frac{\epsilon_{\mathbb{T}} \epsilon_{C^0(s)}}{ |C^0(s)^{\Phi}|_p}  \displaystyle\sum_{i=0}^1 \theta(\upsilon^i(s))$$
since the relative Weyl group is $W(\mathbb{G}(\mathbb{F}_q),\mathbb{T}(\mathbb{F}_q)) = Aut(\mathbb{F}_{q^2}/\mathbb{F}_q)$.  It remains to calculate the constants in front.

Now, since $s \in \mathbb{T}^{\Phi}$ is regular semisimple, $|C^0(s)^{\Phi}|_p = 1$.  Moreover, $\epsilon_{\mathbb{T}} = -1$ and $\epsilon_{C^0(s)} = 1$.  Therefore, $$R_{\mathbb{T}, \theta}(s) = \epsilon_{\mathbb{T}} \epsilon_{C^0(s)} \displaystyle\sum_{i=0}^1 \theta(\upsilon^i(s)) = - (\theta(s) + \theta(\overline{s})$$
\qed

Our character formula is defined on the unramified elliptic torus $E^*$.  We wish to show that our character formula agrees with a depth zero supercuspidal character on the sets where they are both defined, i.e.  the set $(F^* K_0 \setminus F^* K_1) \bigcap E^*$.

\begin{lemma}\label{lemmaweird}
$(F^* K_0 \setminus F^* K_1) \bigcap E^* = F^* A = E^* \setminus F^*(1 + \mathfrak{p}_E) = \{z \in T(F)^{reg} : n(z) = 0 \}$, where $A:= \{p^n u + v \delta : n \geq 0 \}$
\end{lemma}

\proof
The proof is elementary.
\qed

\begin{theorem}\label{depthzerocharactersmatchingup}
$F(\tilde\chi)$ agrees with the supercuspidal character of $\pi_{\chi \Delta_{\chi}}$ on $F^* A$.
\end{theorem}

\proof
Recall that $\Delta_{\chi}$ is the unique quadratic unramified character of $E^*$.  Therefore, we need to show that

$$-\frac{deg(\pi)}{deg(\sigma)} \left(\frac{\chi(w) + (-1, \Delta) \chi(\overline{w})}{ \tau_o(\frac{w - \overline{w}}{2 \delta}) |D(w)|^{1/2})}\right) = -\frac{deg(\pi)}{deg(\sigma)}(\chi(w) \Delta_{\chi}(w) + \chi(\overline{w})\Delta_{\chi}(\overline{w})) \ \ \forall w \in F^* A$$

\noindent Let $w \in A$, so $w = p^n u + v \delta, n \geq 0$.  Then $|D(w)| = 1$.  Moreover, $\tau_o( \frac{w - \overline{w}}{2 \delta} ) = \tau_o(v) = 1$.  But $\Delta_{\chi}(w) = 1 \ \forall w \in A$ since $\Delta_{\chi}$ is unramified.  Therefore, both sides agree on $A$.  Finally, since $\tau_o(x) = \Delta_{\chi}(x) \ \forall x \in F^*$, we have that both sides agree on $F^* A$.
\qed

As in the case of positive depth supercuspidal representations, our overall character formula $F(\tilde\chi)$ remains the same regardless of the choice of positive root.

\subsection{On whether there are two character formulas coming from the same Cartan}

In this section, we show that if the distribution characters of two depth zero supercuspidal representations, both coming from the unramified Cartan, agree on the $n(w) = 0$ range, then the supercuspidal representations are isomorphic.

\begin{theorem}\label{depthzerosamecartans}
Suppose $(E/F, \chi_1)$, $(E/F, \chi_2)$ are admissible pairs such that $F(\tilde\chi_1)(w) = F(\tilde\chi_2)(w)$ on the set $E^* \setminus F^*(1 + \mathfrak{p}_E)$.  Then $\chi_1 = \chi_2^{\upsilon}$ for some $\upsilon \in Aut(E/F)$.
\end{theorem}

\proof
See \cite[page 16]{spice}.
\qed

Summing up, we have altogether shown that if $(E/F, \chi)$ is a regular pair such that $\chi$ has level zero, then there is a unique depth zero supercuspidal representation, $\pi_{\chi \Delta_{\chi}}$, whose character, on the range $\{ z \in T(F)^{reg} : 0 \leq n(z) \leq r/2 \} = \{z \in T(F)^{reg} : n(z) = 0 \}$, agrees with $F(\tilde\chi)$.  As in the argument at the end of Section \ref{othercartan}, one also sees that there is no positive depth supercuspidal representation whose character, on the range $\{z \in T(F)^{reg} : n(z) = 0 \}$, also equals $F(\tilde\chi)$. Therefore, combining Theorems (\ref{depthzerosamecartans}) and (\ref{depthzerocharactersmatchingup}), we obtain the following result.

\begin{theorem}
The assignment

\begin{eqnarray}
\left\{
\begin{array}{rl}
 irreducible \ \varphi : W_F \rightarrow GL(2,\mathbb{C})
\end{array} \right\} & \mapsto & \tilde\chi \in \widehat{T(F)}_{\tau \circ \rho} \mapsto \pi(\tilde\chi) \nonumber
\end{eqnarray}

\noindent from Section (\ref{setup}) is the Local Langlands correspondence for depth zero supercuspidal representations of $GL(2,F)$, where $\pi(\tilde\chi)$ is the unique supercuspidal representation whose character, on the range $\{ z \in T(F)^{reg} : 0 \leq n(z) \leq r/2 \} = \{z \in T(F)^{reg} : n(z) = 0 \}$, is $F(\tilde\chi)$.
\end{theorem}

\section{Existing Description of Local Langlands Correspondence for $GL(\ell,F)$, $\ell$ an odd prime}

In this chapter, we describe the construction of the local Langlands correspondence for $GL(\ell,F)$ as explained in \cite{moy}.

\subsection{Admissible Pairs}\label{moythesis}

Let $E/F$ be a tamely ramified degree $\ell$ extension and $\chi$ a character of $E^*$, where $\ell$ is an odd prime.  Recall that $N$ denotes the norm map from $E$ to $F$.

\begin{definition}
The pair $(E/F, \chi$) is called an \emph{admissible pair} if

(i) $\chi$ does not factor through $N$ and

(ii) If $\chi|_{1 + \mathfrak{p}_E}$ factors through $N$, then $E/F$ is unramified.

\end{definition}

If $(E/F, \chi)$ is an admissible pair, then there is a Howe factorization $\chi = \chi' \phi_E$, where $\phi_E = \phi \circ N_{E/F}$, for some $\phi \in \widehat{F^*}$ and where $\chi' \in \widehat{E^*}$ is of minimal conductor.  We write $\mathbb{P}_{\ell}(F)$ for the set of $F$-isomorphism classes of admissible pairs.  For more notions on admissible pairs, see \cite{moy}.

\subsection{Depth zero and positive depth supercuspidal representations of $GL(\ell,F)$}

Let $\mathbb{A}_{\ell}^0(F)$ denote the set of all irreducible supercuspidal representations of $GL(\ell,F)$.

\begin{theorem}\label{bijection2}
Suppose the residual characteristic of $F$ is not equal to $\ell$.  There is a map $(E/F, \chi)$ $\mapsto \pi_{\chi}$ that induces a bijection $$\mathbb{P}_{\ell}(F) \rightarrow \mathbb{A}_{\ell}^0(F)$$

If $(E/F, \chi) \in \mathbb{P}_{\ell}(F)$, then:

(i) $\omega_{\pi_{\chi}} = \chi|_{F^*}$

(ii) if $\phi$ is a character of $F^*$, then $\pi_{\chi \phi_E} = \phi \pi_{\chi}$.
\end{theorem}

\proof
See \cite{moy}.
\qed

\subsection{Weil parameters}

Let $\mathbb{G}_{\ell}^0(F)$ be the set of equivalence classes of irreducible smooth $\ell$-dimensional representations of $W_F$.  If $(E/F, \xi) \in \mathbb{P}_{\ell}(F)$, we can form $\varphi_{\xi} := Ind_{W_E}^{W_F}(\xi)$ as in Section (\ref{weilparametersgl2}). Let $\delta_{E/F} := det(Ind_{W_E}^{W_F}(1))$.

\begin{definition}
Let $(E/F, \xi) \in \mathbb{P}_{\ell}(F)$ such that $\delta_{E/F} = 1$.  Define $\Delta_{\xi}$ to be the trivial character of $E^*$.
\end{definition}

\begin{definition}
Let $(E/F, \xi) \in \mathbb{P}_{\ell}(F)$ such that $\delta_{E/F} \neq 1$.  Define $\Delta_{\xi}$ to be the unique quadratic unramified character of $E^*$.
\end{definition}

\begin{theorem}{$\mathbf{Tame \ Local \ Langlands \ Correspondence}$}{\label{tamellc}}

Suppose the residual characteristic of $F$ is not equal to $\ell$.  For $\varphi \in \mathbb{G}_{\ell}^0(F)$, define $\pi(\varphi) = \pi_{\xi \Delta_{\xi}}$ in the notation of Theorem (\ref{bijection2}) for any $(E/F, \xi) \in \mathbb{P}_{\ell}(F)$ such that $\varphi \cong \varphi_{\xi}$.  The map $$\pi : \mathbb{G}_{\ell}^0(F) \rightarrow \mathbb{A}_{\ell}^0(F)$$ is the local Langlands correspondence for $GL(\ell,F)$, where $\ell$ is an odd prime.
\end{theorem}

\proof
See \cite{moy}.
\qed

\begin{proposition}\label{centralcharactergll}
If $\varphi \in \mathbb{G}_{\ell}^0(F)$ and $\pi = \pi(\varphi)$, then $\omega_{\pi} =$ det($\varphi$).
\end{proposition}

\proof
See \cite{moy}.
\qed

\section{Positive depth supercuspidal character formula for $GL(\ell,F)$ and $PGL(\ell,F)$, $\ell$ an odd prime}

\subsection{Preliminaries}\label{preliminaries}

In this Chapter we prove Theorems (\ref{gl2theorem1}) and (\ref{gl2theorem2}) for the positive depth supercuspidal representations of $GL(\ell,F)$, where $\ell$ is an odd prime.  We will define $\epsilon(\tilde\chi, \Delta^+, \tau)$ and the notion of regular in the next section.  We will also regularly use the fact that $W(G(F),T(F)) = Aut(E/F)$.

We note that all of our calculations in the next two chapters will assume that we have chosen the standard positive set of roots of $GL(\ell,\overline{F})$ with respect to the standard split maximal torus.  Our main results, however, will be seen to be independent of any choice of positive roots.

Now let $\varphi$ be a supercuspidal Weil parameter for $GL(\ell,F)$.  We will show later in this section how to construct a regular genuine character, $\tilde\chi$, of $T(F)_{\tau \circ \rho}$, from $\varphi$.  We will then prove Theorem (\ref{gl2theorem2}) for $GL(\ell,F)$ where $\ell$ is an odd prime.

We now introduce a notion of regularity that we will need.  Let $E/F$ be a tamely ramified degree $\ell$ extension and $\chi$ a character of $E^*$.  Recall that $N$ denotes the norm map from $E$ to $F$.

\begin{definition}
$\chi$ is called \emph{regular} if $\chi$ does not factor through $N$.  If $\chi$ is regular, we call the pair $(E/F, \chi)$ a \emph{regular pair}.
\end{definition}

All definitions we have made in the previous chapter for admissible pairs, we also make for regular pairs and regular characters as in the case of $GL(2,F)$.  In particular, we also define the character twists $\Delta_{\chi}$ for a regular pair $(E/F, \chi)$ exactly the same way they were defined for admissible pairs.  For example, if $(E/F, \chi)$ is a regular pair where $\delta_{E/F} \neq 1$, then $\Delta_{\chi}$ is the unique unramified quadratic character of $E^*$.  Given a regular pair $(E/F, \chi)$, one may also construct a supercuspidal representation $\pi_{\chi}$ as in the previous chapter, but this construction is not one to one.

Our constructions and results do not require the stronger notion of admissible pair.  We will sometimes say that $\chi$ is regular when the field $E$ is understood.

We first explain why double covers of tori play a role.  We start by considering the group $PGL(\ell,F)$.  First recall that the representations of $PGL(\ell,F)$ are precisely the representations of $GL(\ell,F)$ with trivial central character.  Let $\varphi$ be a supercuspidal Weil parameter for $PGL(\ell,F)$.  Then $\varphi = Ind_{W_E}^{W_F}(\chi)$ for some regular pair $(E/F, \chi)$.  Since we are using the notion of regular pair here rather than admissible pair, there may be a choice involved here.  That is, there may be another regular pair $(E_1/F, \chi_1)$ such that $\varphi = Ind_{W_{E_1}}^{W_F}(\chi_1)$ as well.  However, this will not matter, and we will show that our results and constructions are independent of all choices.  It is a fact (see \cite[Proposition 29.2]{bushnellhenniart}) that $det(Ind_{W_E}^{W_F}(\chi)) = \chi|_{F^*} \otimes \delta_{E/F}$, where $\delta_{E/F} = det(Ind_{W_E}^{W_F}(1))$.  Thus, by Proposition (\ref{centralcharactergll}) $\chi|_{F^*} = \delta_{E/F}$. Therefore, the supercuspidal Weil parameters for $PGL(\ell,F)$ naturally give rise to regular pairs $(E/F, \chi)$ where $\chi|_{F^*} = \delta_{E/F}$.  Such a $\chi$ is not necessarily a character of the elliptic torus $E^*/F^*$ in $PGL(\ell,F)$.  Rather, it is a genuine character of a cover $E^* / ker(\delta_{E/F})$ of $E^*/F^*$ arising from the canonical exact sequence $$1 \longrightarrow F^* / ker(\delta_{E/F}) \longrightarrow E^* / ker(\delta_{E/F}) \rightarrow E^*/F^* \longrightarrow 1$$  In the case that $E/F$ is degree $\ell$, $\delta_{E/F}^2 = 1$ (see \cite[Corollary 2.5.15]{moy}). We first consider the case where $\delta_{E/F} \neq 1$.  Then since $F^* / ker(\delta_{E/F}) \cong \mathbb{Z}/ 2 \mathbb{Z}$, we have that $E^* / ker(\delta_{E/F})$ is a nontrivial double cover of $E^*/F^*$.  Then the character $\chi$ of $E^*$ naturally factors to a genuine character $\tilde\chi$ of $E^* / ker(\delta_{E/F})$, given by $\tilde\chi([w]) := \chi(w) \ \forall [w] \in E^* / ker(\delta_{E/F})$.  Therefore, the supercuspidal Weil parameters for $PGL(\ell,F)$ naturally give rise to genuine characters of the double cover $E^* / ker(\delta_{E/F})$.  In the case of $\delta_{E/F} = 1$, there is no cover, since $E^* / ker(\delta_{E/F}) = E^*/F^*$.  Therefore, in this case, supercuspidal Weil parameters for $PGL(\ell,F)$ naturally give rise to characters $\chi$ of the elliptic torus $E^* / F^*$ inside $PGL(\ell,F)$.  There is still a small subtlety in this case, as the proposed character formula in Theorem (\ref{gl2theorem1}) is written in terms of a double cover of the elliptic tori.  Since the denominator of the proposed character formula lives on a double cover, we must have that the numerator of the proposed character formula lives on a double cover as well, so we will need a method of going from the character $\chi$ of $E^*/F^*$ to a genuine character $\tilde\chi$ of a double cover.

We first consider the case where $\delta_{E/F} \neq 1$.  Relative to the standard positive system of roots of $PGL(\ell,F)$, let $\rho$ be half the sum of the positive roots.  An elliptic torus in $PGL(\ell,F)$ is of the form $T(F) = E^* / F^*$.  Fix a character $\tau_o$ of $(EL)^*$ whose restriction to $L^*$ is $\aleph_{EL/L}$, where $L$ is the unramified extension of $F$ of degree $\ell-1$ and set $\tau := \tau_o \ | \ |_{EL}$.  Recall the denominator $$\tau(\Delta^0(z, \Delta^+)) (\tau \circ \rho)(w)$$ that was defined in Theorem (\ref{gl2theorem1}).  Although $\tau \circ \rho$ is not naturally a function on $E^* / F^*$ since in particular $\rho$ is not naturally a function on $E^* / F^*$, it is naturally a function on the $\tau \circ \rho$-cover of $E^* / F^*$.  Recall Definition (\ref{tauofrhocover}).  Then, in our case, $T(F)_{\tau \circ \rho} = \{([z],w) \in E^* / F^* \times \mathbb{C}^* : \tau(2 \rho([z])) = w^2 \}$

We make an important note here.  Note that since $\ell$ is odd, not only do we have $2 \rho \in X^*(T)$, but we also have $\rho \in X^*(T)$.  Therefore, if we consider $\rho$ as a function on $T$ as such, we may apply $\tau$ to the element $\rho(w)$ where $w \in E^*$.  We will denote this resulting function $\rho_{\tau}(w)$, as this is a different function than the function $(\tau \circ \rho)(w)$ which naturally lives on $T(F)_{\tau \circ \rho}$.  We make this a formal definition.

\begin{definition}
We define $$\rho_{\tau}(w) := \tau(\rho(w)), \ w \in E^*$$
Here, we are viewing $\rho$ as an element of $X^*(T)$, which we may do since $\ell$ is odd.
\end{definition}

\begin{lemma}\label{doublecovertorigll}
$E^* / ker(\delta_{E/F}) \cong T(F)_{\tau \circ \rho}$, and this isomorphism is unique as an isomorphism of covering groups (i.e. as covers of $T(F)$).
\end{lemma}

\proof
The isomorphism is given by $$E^* / ker(\delta_{E/F}) \xrightarrow{\kappa} T(F)_{\tau \circ \rho}$$ $$\ \ \ \ \ \ \ \ \ \ [w] \mapsto ([w], \Delta_{\chi}(w) \rho_{\tau}(w))$$
\noindent where on the right hand side, $[w]$ lives in $E^* / F^*$.  It is easy to show that this covering isomorphism is unique, since $E^* = (E^*)^2 F^*$ (here we are using that $\ell$ is an odd prime).
\qed

Now let's write down the character formula for a supercuspidal representation of $PGL(\ell,F)$ in the case where $\delta_{E/F} \neq 1$.  Let $\varphi : W_F \rightarrow GL(\ell,\mathbb{C})$ be a supercuspidal parameter for $PGL(\ell,F)$ so that $\varphi = Ind_{W_E}^{W_F}(\chi)$ for some regular pair $(E/F, \chi)$.  As discussed earlier, this gives us a genuine character $\tilde\chi$ of $E^* / ker(\delta_{E/F})$.

\begin{definition}
A genuine character $\tilde\eta$ of $E^* / ker(\delta_{E/F})$ is called \emph{regular} if $(E/F, \eta)$ is regular, where $\eta$ is the pullback of $\tilde\eta$ to $E^*$.  A genuine character $\tilde\lambda$ of $T(F)_{\tau \circ \rho}$ is called \emph{regular} if $\tilde\lambda \circ \kappa$ is regular.
\end{definition}

Now recall from Theorem (\ref{gl2theorem1}) the proposed formula $F(\tilde\chi)$.  As in the cases of $PGL(2,F)$ and $GL(2,F)$, we pull the function $(\tau \circ \rho)(w)$ and the Weyl group action in $F(\tilde\chi)$ back to $E^* / ker(\delta_{E/F})$ via $\kappa$, and leave our constructed $\tilde\chi$ as living on $E^* / ker(\delta_{E/F})$.  That is, we consider $$F(\tilde\chi)(z) = \epsilon(\tilde\chi, \Delta^+, \tau) \frac{\displaystyle\sum_{s \in W} \epsilon(s)\tilde\chi({}^s [w])}{\tau(\Delta^0(z,\Delta^+)) (\tau \circ \rho)(\kappa([w]))}, \ \ z \in T(F)^{reg}$$ where $[w] \in E^* / ker(\delta_{E/F})$ such that $\Pi(\kappa([w])) = z$.  Unwinding the definitions, we see that $(\tau \circ \rho)(\kappa([w])) = \Delta_{\chi}(w) \rho_{\tau}(w) \ \forall [w] \in E^* / ker(\delta_{E/F})$.

We need to define the Weyl group action.  As in the case of $GL(2,F)$, define $s([w],\lambda) = (s[w],\lambda \tau((s^{-1} \rho - \rho)([w])))$ for $s \in W$.  Pulling back this action from $T(F)_{\tau \circ \rho}$ to $E^* / ker(\delta_{E/F})$ via $\kappa$, this simplifies to $s [w] = \kappa^{-1}(s \kappa([w])) = [s w] \ \forall [w] \in E^* / ker(\delta_{E/F})$ for $s \in W = Aut(E/F)$.

We note that the definition of regularity for a genuine character of $T(F)_{\tau \circ \rho}$ is analogous to the definition of regularity for a genuine character $\tilde\lambda$ of $T(\mathbb{R})_{\rho}$ for real groups when $E/F$ is Galois.

Pulling back $(\tau \circ \rho)(w)$ and the Weyl group action to $E^* / ker(\delta_{E/F})$ via $\kappa$, we get $$F(\tilde\chi)(z) = \epsilon(\tilde\chi, \Delta^+, \tau) \frac{\displaystyle\sum_{s \in W} \epsilon(s) \tilde\chi({}^s [w]) }{\tau(\Delta^0(z, \Delta^+)) \Delta_{\chi}(w) \rho_{\tau}(z)} = \epsilon(\tilde\chi, \Delta^+, \tau) \frac{\displaystyle\sum_{s \in W}  \tilde\chi({}^s [w]) \Delta_{\chi}({}^s w) }{\tau(\Delta^0(z, \Delta^+)) \rho_{\tau}(z)}$$ where $z \in E^*/F^*$ and $[w] \in E^* / ker(\delta_{E/F})$ is any element such that $[w]$ maps to $z$ under the canonical map $E^* / ker(\delta_{E/F}) \rightarrow E^* / F^*$.

Let us now compute the proposed character formula for $GL(\ell,F)$.  Relative to the standard positive system of roots of $GL(\ell,F)$, let $\rho$ be half the sum of the positive roots.  An elliptic torus in $GL(\ell,F)$ is of the form $T(F) = E^*$.  We now introduce a cover which is isomorphic to $T(F)_{\tau \circ \rho}$.

\begin{definition}
Let $E^* / ker(\delta_{E/F}) \rightarrow E^* / F^*$ be the canonical projection map.  We define $E^* \times_{E^*/F^*} E^* / ker(\delta_{E/F})$ as the group arising in the following pullback diagram:

$$
\begin{CD}
E^* \times_{E^*/F^*} E^* / ker(\delta_{E/F}) @>>> E^* / ker(\delta_{E/F})\\
@VVV @VVV\\
E^* @>w \mapsto [w]>> E^*/F^*
\end{CD}
$$

\noindent That is, $E^* \times_{E^*/F^*} E^* / ker(\delta_{E/F}) = \{(w,[z]) \in E^* \times E^* / ker(\delta_{E/F}) : [w] = [z] \in E^* / F^* \}$
\end{definition}

\begin{lemma}\label{kappaglell1}
$E^* \times_{E^*/F^*} E^* / ker(\delta_{E/F}) \cong T(F)_{\tau \circ \rho}$
\end{lemma}

\proof
We naturally have an isomorphism $$E^* \times_{E^*/F^*} E^* / ker(\delta_{E/F}) \xrightarrow{\kappa} T(F)_{\tau \circ \rho}$$ $$(w,[z]) \mapsto (w, \Delta_{\chi}(z) \rho_{\tau}(w))$$ The proof that this is an isomorphism follows as in the previous cases.
\qed

Now let's write down the character formula for a supercuspidal representation of $GL(\ell,F)$.  Let $\varphi : W_F \rightarrow GL(\ell, \mathbb{C})$ be a supercuspidal parameter so that $\varphi = Ind_{W_E}^{W_F}(\chi)$ for some regular pair $(E/F, \chi)$.  Then this canonically gives a genuine character $\tilde\chi$ of $E^* \times_{E^*/F^*} E^* / ker(\delta_{E/F})$ as follows.  Define $\tilde\chi(w,[z]) := \chi(w) \delta_{E/F}(z/w)$.

\begin{definition}
A genuine character $\tilde\eta$ of $E^* \times_{E^*/F^*} E^* / ker(\delta_{E/F})$ is called \emph{regular} if $(E/F, \eta)$ is regular, where $\eta(w) := \tilde\eta(w,[z]) \delta_{E/F}(z/w)$.  A genuine character $\tilde\lambda$ of $T(F)_{\tau \circ \rho}$ is called \emph{regular} if $\tilde\lambda \circ \kappa$ is regular.
\end{definition}

We have therefore given a map $\widehat{E^*} \rightarrow (E^* \times_{E^*/F^*} E^* / ker(\delta_{E/F}))^{\wedge}$ given by $\eta \mapsto \tilde\eta$, where $\tilde\eta(w,[z]) := \eta(w) \delta_{E/F}(z/w)$.  Note that we have a canonical map in the other direction, $(E^* \times_{E^*/F^*} E^* / ker(\delta_{E/F}))^{\wedge} \rightarrow \widehat{E^*}$, given by $\tilde\eta \mapsto \eta$, where $\eta(w) := \tilde\eta(w,[z]) \delta_{E/F}(z/w)$.  We will regularly go back and forth between characters of $E^*$ and genuine characters of $E^* \times_{E^*/F^*} E^* / ker(\delta_{E/F})$.  In particular, when we write $\tilde\chi$, a genuine character of $E^* \times_{E^*/F^*} E^* / ker(\delta_{E/F})$, we will sometimes keep in mind that there is a canonical character $\chi$ of $E^*$ that $\tilde\chi$ comes from via the above maps.

Now recall from Theorem (\ref{gl2theorem1}) the proposed formula $F(\tilde\chi)$. As in the previous cases, we pull the function $(\tau \circ \rho)(w)$ and the Weyl group action in $F(\tilde\chi)$ back to $E^* \times_{E^*/F^*} E^* / ker(\delta_{E/F})$ via $\kappa$, and leave our constructed $\tilde\chi$ as living on $E^* \times_{E^*/F^*} E^* / ker(\delta_{E/F})$.  That is, we consider $$F(\tilde\chi)(w) = \epsilon(\tilde\chi, \Delta^+, \tau) \frac{\displaystyle\sum_{s \in W} \epsilon(s)\tilde\chi({}^s (w,[z]))}{\tau(\Delta^0(w,\Delta^+)) (\tau \circ \rho)(\kappa(w,[z]))} \ \ w \in T(F)^{reg}$$ where $(w,[z]) \in E^* \times_{E^* / F^*} E^* / ker(\delta_{E/F})$ such that $\Pi(\kappa((w,[z]))) = w$.  Unwinding the definitions, we see that $(\tau \circ \rho)(\kappa((w,[z]))) = \Delta_{\chi}(z) \rho_{\tau}(w) \ \forall (w,[z]) \in E^* \times_{E^* / F^*} E^* / ker(\delta_{E/F})$.

We also need to define the Weyl group action.  As in the case of $GL(2,F)$, define $s(w,\lambda) = (sw,\lambda \tau((s^{-1} \rho - \rho)(w)))$ for $s \in W$.  Pulling back this action from $T(F)_{\tau \circ \rho}$ to $E^* \times_{E^*/F^*} E^* / ker(\delta_{E/F})$ via $\kappa$, this simplifies to $s(w,[z]) = (sw, [sz]) \forall (w,[z]) \in E^* \times_{E^*/F^*} E^* / ker(\delta_{E/F})$ for $s \in W$.

Pulling back $\tau \circ \rho$ and the Weyl group action to $E^* \times_{E^* / F^*} E^* / ker(\delta_{E/F})$ via $\kappa$, our character formula simplifies to
$$F(\tilde\chi)(w) = \epsilon(\tilde\chi, \Delta^+, \tau) \frac{\displaystyle\sum_{s \in W} \epsilon(s) \chi({}^s w) \delta_{E/F}(z/w)}{\tau(\Delta^0(w, \Delta^+)) \rho_{\tau}(w) \Delta_{\chi}(z)}, \ \ w \in T(F)^{reg}$$ where $(w,[z]) \in E^* \times_{E^*/F^*} E^* / ker(\delta_{E/F})$ is any element that maps to $w$ under the canonical projection $E^* \times_{E^*/F^*} E^* / ker(\delta_{E/F}) \rightarrow E^*$. Note that $\Delta_{\chi}|_{F^*} = \delta_{E/F}$, so $\delta_{E/F}(z/w) = \Delta_{\chi}(z/w)$.  Therefore, $$F(\tilde\chi)(w) = \epsilon(\tilde\chi, \Delta^+, \tau) \frac{\displaystyle\sum_{s \in W} \epsilon(s) \chi({}^s w) \Delta_{\chi}(z/w)}{\tau(\Delta^0(w, \Delta^+)) \rho_{\tau}(w) \Delta_{\chi}(z)}, \ \ w \in T(F)^{reg} $$ where $(w,[z]) \in E^* \times_{E^*/F^*} E^* / ker(\delta_{E/F})$ is any element that maps to $w$ under the canonical projection $E^* \times_{E^*/F^*} E^* / ker(\delta_{E/F}) \rightarrow E^*$. Then, since $\Delta_{\chi}^2 = 1$ and $\epsilon(s) = 1 \ \forall s \in W$, this formula reduces to $$F(\tilde\chi)(w) = \epsilon(\tilde\chi, \Delta^+, \tau) \frac{\displaystyle\sum_{s \in W} \chi({}^s w) \Delta_{\chi}({}^s w)}{\tau(\Delta^0(w, \Delta^+)) \rho_{\tau}(w)}, \ \ w \in T(F)^{reg}$$

\

We now consider the case of $PGL(\ell,F)$ where $\delta_{E/F} = 1$.  Relative to the standard positive system of roots of $PGL(\ell,F)$, let $\rho$ be half the sum of the positive roots.  An elliptic torus in $PGL(\ell,F)$ is of the form $T(F) = E^* / F^*$.  Recall that in this case, a supercuspidal parameter for $PGL(\ell,F)$ does not naturally yield a genuine character of a double cover of $E^* / F^*$.  Rather, we are naturally handed a character of $E^* / ker(\delta_{E/F}) = E^* / F^*$, since $\delta_{E/F} = 1$.  Therefore, the only natural double cover to consider in this case is the canonical split cover $E^* / F^* \times \mathbb{Z} / 2 \mathbb{Z}$. We will show that this setting still naturally fits into our theory.

\begin{lemma}
$E^* / F^* \times \mathbb{Z} / 2 \mathbb{Z} \cong T(F)_{\tau \circ \rho}$, and this isomorphism is unique as an isomorphism of covering groups (i.e. as covers of $T(F)$).
\end{lemma}

\proof
The explicit isomorphism is given by
$$E^* / F^* \times \mathbb{Z} / 2 \mathbb{Z} \xrightarrow{\kappa} T(F)_{\tau \circ \rho}$$ $$\ \ \ \ \ \ \ (z,\epsilon) \mapsto (z, \epsilon \rho_{\tau}(z))$$
It is easy to see that this isomorphism is unique, using the fact that $E^* = (E^*)^2 F^*$ (here we are using that $\ell$ is an odd prime.
\qed

Now let's write down the character formula for a supercuspidal representation of $PGL(\ell,F)$ in the case where $\delta_{E/F} = 1$.  Let $\varphi : W_F \rightarrow GL(\ell,\mathbb{C})$ be a supercuspidal parameter for $PGL(\ell,F)$ so that $\varphi = Ind_{W_E}^{W_F}(\chi)$ for some regular pair $(E/F, \chi)$.  Now recall that we are trying to make sense of the proposed character formula $F(\tilde\chi)$.  We have that $\chi$ factors to a character of $E^*/F^*$, which we will also denote $\chi$.  However, the functions in $F(\tilde\chi)$ have domain $T(F)_{\tau \circ \rho}$.  As in previous cases, we pull the function $(\tau \circ \rho)(w)$ and the Weyl group action in $F(\tilde\chi)$ back to $E^* / F^* \times \mathbb{Z} / 2 \mathbb{Z}$ via $\kappa$.  That is, we consider $$F(\tilde\chi)(z) = \epsilon(\tilde\chi, \Delta^+, \tau) \frac{\displaystyle\sum_{s \in W} \epsilon(s)\tilde\chi({}^s (z, \epsilon))}{\tau(\Delta^0(z,\Delta^+)) (\tau \circ \rho)(\kappa((z,\epsilon)))}, \ \ z \in T(F)^{reg}$$ where $(z,\epsilon) \in E^* / F^* \times \mathbb{Z} / 2 \mathbb{Z}$ such that $\Pi(\kappa((z, \epsilon))) = z$.  Unwinding the definitions, we see that $(\tau \circ \rho)(\kappa((z, \epsilon))) = \epsilon \rho_{\tau}(z) \ \forall (z, \epsilon)  \in E^*/F^* \times \mathbb{Z} / 2 \mathbb{Z}$.  We can then canonically assign a genuine character $\tilde\chi$ of $E^* / F^* \times \mathbb{Z} / 2 \mathbb{Z}$ from the regular character $\chi$ by setting $$\tilde\chi := \chi \otimes sgn$$

\begin{definition}
A genuine character $\tilde\eta$ of $E^* \times \mathbb{Z} / 2 \mathbb{Z}$ is called \emph{regular} if $(E/F, \eta)$ is regular, where $\eta := \tilde\eta \otimes sgn$.  A genuine character $\tilde\lambda$ of $T(F)_{\tau \circ \rho}$ is called \emph{regular} if $\tilde\lambda \circ \kappa$ is regular.
\end{definition}

We also need to define the Weyl group action.  As in the case of $GL(2,F)$, define $s([w],\lambda) = ([sw],\lambda \tau((s^{-1} \rho - \rho)(w)))$ for $s \in W$.  Pulling back this action from $T(F)_{\tau \circ \rho}$ to $E^* / F^* \times \mathbb{Z} / 2 \mathbb{Z}$ via $\kappa$, this simplifies to $s([z],\epsilon) = ([sz], \epsilon) \ \forall ([z], \epsilon) \in E^* / F^* \times \mathbb{Z} / 2 \mathbb{Z}$ for $s \in W$.

Pulling back $(\tau \circ \rho)(w)$ and the Weyl group action to $E^* / F^* \times \mathbb{Z} / 2 \mathbb{Z}$ via $\kappa$, incorporating $\tilde\chi$, and noting that $\epsilon(s) = 1 \ \forall s \in W$, we get $$F(\tilde\chi)(z) = \epsilon(\tilde\chi, \Delta^+, \tau) \frac{\displaystyle\sum_{s \in W} \epsilon(s) \epsilon \chi({}^s z) }{\tau(\Delta^0(z, \Delta^+)) \epsilon \rho_{\tau}(z)} = \epsilon(\tilde\chi, \Delta^+, \tau) \frac{\displaystyle\sum_{s \in W}  \chi({}^s z) }{\tau(\Delta^0(z, \Delta^+)) \rho_{\tau}(z)}, \ \ z \in E^* / F^*$$

\

Let us now compute the proposed character formula for $GL(\ell,F)$ in the case that $\delta_{E/F} = 1$.  Relative to the standard system of roots in $GL(\ell,F)$, let $\rho$ be half the sum of the positive roots.  An elliptic torus in $GL(\ell,F)$ is of the form $T(F) = E^*$.  We now introduce a cover which is isomorphic to $T(F)_{\tau \circ \rho}$, completely analogously to the case of $\delta_{E/F} \neq 1$.

\begin{definition}
Let $E^* / F^* \times \mathbb{Z} / 2 \mathbb{Z} \rightarrow E^* / F^*$ be the canonical projection map.  We define $E^* \times_{E^*/F^*} (E^* / F^* \times \mathbb{Z} / 2 \mathbb{Z})$ as the group arising in the following pullback diagram:

$$
\begin{CD}
E^* \times_{E^*/F^*} (E^* / F^* \times \mathbb{Z} / 2 \mathbb{Z})@>>> E^* / F^* \times \mathbb{Z} / 2 \mathbb{Z}\\
@VVV @VVV\\
E^* @>w \mapsto [w]>> E^*/F^*
\end{CD}
$$

\noindent That is, $E^* \times_{E^*/F^*} (E^* / F^* \times \mathbb{Z} / 2 \mathbb{Z}) = \{(w,([z],\epsilon)) : [w] = [z] \in E^* / F^* \}$
\end{definition}

\begin{lemma}
$E^* \times_{E^*/F^*} (E^*/F^* \times \mathbb{Z} / 2 \mathbb{Z}) \cong T(F)_{\tau \circ \rho}$
\end{lemma}

\proof
An explicit isomorphism is given by $$E^* \times_{E^*/F^*} (E^*/F^* \times \mathbb{Z} / 2 \mathbb{Z}) \xrightarrow{\kappa} T(F)_{\tau \circ \rho}$$ $$(w,([z],\epsilon)) \mapsto (w, \epsilon \rho_{\tau}(w))$$
\qed

Now let's write down the character formula for a supercuspidal representation of $GL(\ell,F)$.  Let $\varphi : W_F \rightarrow GL(\ell, \mathbb{C})$ be a supercuspidal parameter so that $\varphi = Ind_{W_E}^{W_F}(\chi)$ for some regular pair $(E/F, \chi)$.  Then this gives a genuine character $\tilde\chi$ of $E^* \times_{E^*/F^*} (E^*/F^* \times \mathbb{Z} / 2 \mathbb{Z})$ as follows.  Define $\tilde\chi(w,([z],\epsilon)) := \chi(w) \epsilon$.

\begin{definition}
A genuine character $\tilde\eta$ of $E^* \times_{E^*/F^*} (E^*/F^* \times \mathbb{Z} / 2 \mathbb{Z})$ is called \emph{regular} if $(E/F, \eta)$ is regular, where $\eta(w) := \tilde\eta(w,([z],\epsilon)) \epsilon$.  A genuine character $\tilde\lambda$ of $T(F)_{\tau \circ \rho}$ is called \emph{regular} if $\tilde\lambda \circ \kappa$ is regular.
\end{definition}

We have therefore given a map $\widehat{E^*} \rightarrow (E^* \times_{E^*/F^*} (E^*/F^* \times \mathbb{Z} / 2 \mathbb{Z}))^{\wedge}$ given by $\eta \mapsto \tilde\eta$, where $\tilde\eta(w,([z],\epsilon)) := \eta(w) \epsilon$.  Note that we have a canonical map in the other direction, $(E^* \times_{E^*/F^*} (E^*/F^* \times \mathbb{Z} / 2 \mathbb{Z}))^{\wedge} \rightarrow \widehat{E^*}$, given by $\tilde\eta \mapsto \eta$, where $\eta(w) := \tilde\eta(w,([z],\epsilon)) \epsilon$.  We will regularly go back and forth between characters of $E^*$ and genuine characters of $E^* \times_{E^*/F^*} (E^*/F^* \times \mathbb{Z} / 2 \mathbb{Z})$.  In particular, when we write $\tilde\chi$, a genuine character of $E^* \times_{E^*/F^*} (E^*/F^* \times \mathbb{Z} / 2 \mathbb{Z})$, we will sometimes keep in mind that there is a canonical character $\chi$ of $E^*$ that $\tilde\chi$ comes from via the above maps.

As in previous cases, we pull the function $(\tau \circ \rho)(w)$ and the Weyl group action in $F(\tilde\chi)$ back to $E^* \times_{E^*/F^*} (E^* / F^* \times \mathbb{Z} / 2 \mathbb{Z})$ via this isomorphism, and leave our constructed $\tilde\chi$ as living on $E^* \times_{E^*/F^*} (E^* / F^* \times \mathbb{Z} / 2 \mathbb{Z})$.  That is, we consider $$F(\tilde\chi)(w) = \epsilon(\tilde\chi, \Delta^+, \tau) \frac{\displaystyle\sum_{s \in W} \epsilon(s) \tilde\chi({}^s (w,([z], \epsilon)))}{\tau(\Delta^0(w, \Delta^+)) (\tau \circ \rho)(\kappa(w,([z], \epsilon)))} \ \ w \in T(F)^{reg} $$ where $(w,([z], \epsilon)) \in E^* \times_{E^*/F^*} (E^*/F^* \times \mathbb{Z} / 2 \mathbb{Z})$ such that $\Pi(\kappa((w,([z], \epsilon)))) = w$.  Unwinding the definitions, we see that $(\tau \circ \rho)(\kappa((w,([z], \epsilon)))) = \epsilon \rho_{\tau}(w) \ \forall (w,([z], \epsilon)) \in E^* \times_{E^*/F^*} (E^*/F^* \times \mathbb{Z} / 2 \mathbb{Z})$.

We also need to define the Weyl group action.  As in the case of $GL(2,F)$, define $s(w,\lambda) = (sw,\lambda \tau((s^{-1} \rho - \rho)(w)))$ for $s \in W$.  Pulling back this action from $T(F)_{\tau \circ \rho}$ to $E^* \times_{E^*/F^*} (E^*/F^* \times \mathbb{Z} / 2 \mathbb{Z})$ via $\kappa$, this simplifies to $s(w,([z], \epsilon)) = (sw, ([sz], \epsilon)) \ \forall (w, ([z], \epsilon)) \in E^* \times_{E^*/F^*} (E^*/F^* \times \mathbb{Z} / 2 \mathbb{Z})$ for $s \in W = Aut(E/F)$.

Pulling back $\tau \circ \rho$ and the Weyl group action to $E^* \times_{E^*/F^*} E^* / F^* \times \mathbb{Z} / 2 \mathbb{Z}$ via $\kappa$, we get $$F(\tilde\chi)(w) = \epsilon(\tilde\chi, \Delta^+, \tau) \frac{\displaystyle\sum_{s \in W} \epsilon(s) \chi({}^s w) \epsilon}{\tau(\Delta^0(w, \Delta^+)) \epsilon \rho_{\tau}(w)} = \epsilon(\tilde\chi, \Delta^+, \tau) \frac{\displaystyle\sum_{s \in W}  \chi({}^s w) }{\tau(\Delta^0(w, \Delta^+)) \rho_{\tau}(w)}, \ \ w \in T(F)^{reg} $$ where $(w,(z, \epsilon)) \in E^* \times_{E^*/F^*} (E^*/F^* \times \mathbb{Z} / 2 \mathbb{Z})$ is any element that maps to $w$ under the canonical projection $E^* \times_{E^*/F^*} (E^*/F^* \times \mathbb{Z} / 2 \mathbb{Z}) \rightarrow E^*$.

Summing up, we have given a method of assigning a conjectural character formula for a supercuspidal representation of $GL(\ell,F)$ or $PGL(\ell,F)$ to a supercuspidal Weil parameter $\varphi$ of $GL(\ell,F)$ or $PGL(\ell,F)$, respectively, given by

$$\varphi \mapsto \tilde\chi \in \widehat{T(F)}_{\tau \circ \rho} \mapsto F(\tilde\chi)$$

We will show that the proposed character formula $F(\tilde\chi)$ constructed in this section is independent of the choise of $\tau$ and the choice of positive roots $\Delta^+$.

We wish to make an important note. The case $\delta_{E/F} = 1$ is the only case where the ``naive correspondence'' of \cite{bushnellhenniart}, \cite{moy}, is the actual local Langlands correspondence.  This is precisely because there is nothing interesting to see in the double cover.  It is clear that in this special case, we didn't need to use double covers in order to obtain the local Langlands correspondence, and this is the only case where we could avoid using double covers of tori.  What we are showing in this paper, however, is that if we move to the setting of double covers of elliptic tori, then we can obtain the local Langlands correspondence in all cases by a ``naive correspondence''.

\subsection{The constant $\epsilon(\tilde\chi, \Delta^+, \tau)$}\label{epsilonforell}

In this section we define the constant $\epsilon(\tilde\chi, \Delta^+, \tau)$.  First recall Theorem (\ref{DeBacker1}). Recall the constant $C := c_{\psi}(\mathfrak{g'}) c_{\psi}^{-1}(\mathfrak{g}) |D(\gamma)|^{-1/2} |\eta(\alpha(\chi))|^{-1/2}$ that is defined in Theorem (\ref{DeBacker1}).

\begin{definition}
Define $\epsilon(\Delta^+)$ to be 1 if $\Delta^+$ is the standard choice of positive roots of $GL(\ell,\overline{F})$ with respect to the diagonal maximal torus $T(\overline{F})$.  Any other set of positive roots is of the form $s \Delta^+$ for some $s \in W(G(\overline{F}), T(\overline{F}))$.  We then define $\epsilon(s \Delta^+) = \tau_o((-1)^{\ell(s)})$ where $\ell$ is the length of $s$.  Let $\Delta^+$ be any set of positive roots.  We set $$\epsilon(\tilde\chi, \Delta^+, \tau) := deg(\pi) \lambda(\sigma) c_{\psi}(\mathfrak{g'}) c_{\psi}^{-1}(\mathfrak{g}) |\eta(\alpha(\chi))|^{-1/2} \tau_o((-1)^{\sum_{k=1}^{\frac{\ell-1}{2}} k}) \epsilon(\Delta^+)$$.
\end{definition}

It will be useful to define another constant.  We set $$\epsilon(\tilde\chi, \Delta^+, \tau)' := deg(\pi) \lambda(\sigma) c_{\psi}(\mathfrak{g'}) c_{\psi}^{-1}(\mathfrak{g}) |\eta(\alpha(\chi))|^{-1/2} \epsilon(\Delta^+).$$

In the calculations we will make throughout the rest of this chapter and the next, we will make a choice of $\Delta^+$ to be the standard set of positive roots.  Therefore, the term $\epsilon(\Delta^+)$ is just $1$, and therefore this term will not appear in most of our calculations and formulas.  We will show later that all of our results will be independent of the choice of $\Delta^+$.

\subsection{The case $\delta_{E/F} \neq 1$}

In the next two sections, we show that the proposed character formula $F(\tilde\chi)$ agrees with the character of the positive depth supercuspidal representation $\pi_{\chi \Delta_{\chi}}$ occuring in the local Langlands correspondence, on the range $\{z \in T(F)^{reg} : 0 \leq n(z) \leq r/2 \}$.  We will assume again without loss of generality, as in Chapter 5, that our regular pairs $(E/F, \chi)$ are such that $\chi$ has minimal conductor (In \cite{spice}, the terminology conductor is used instead of the terms ``minimal regular pair''.  We follow the terminology in \cite{spice} since we will use results from there).  The same argument as in the end of Chapter 5 shows that this doesn't matter, and that all of our results are true for arbitrary regular pairs.  We need a result from \cite{spice}.  Unwinding all the definitions, it is shown in \cite[Sections 6 and 7]{spice} that if $(E/F, \chi)$ is a regular pair with positive level, then the character, $\theta$, of the associated supercuspidal representation $\pi_{\chi}$ via Theorem (\ref{bijection2}), satisfies $$\theta(w) = deg(\pi) \lambda(\sigma) c_{\psi}(\mathfrak{g'}) c_{\psi}^{-1}(\mathfrak{g}) |\eta(\alpha(\chi))|^{-1/2} \frac{\displaystyle\sum_{s \in Aut(E/F)} \chi({}^s w)}{|D(w)|^{1/2}} \ \ \ \ \ \forall w \in E^* : 0 \leq n(w) \leq r/2$$

Let $\Delta$ be the standard set of roots for $GL(\ell,F)$ with respect to the diagonal torus $T$, where $\ell$ is an odd prime, and let $\Delta^+$ be the standard set of positive roots.  In this section, we will view $\rho$ as a character of $T$, so that we may eventually compute $\rho_{\tau}$ (recall that since $\ell$ is odd, $\rho \in X^*(T)$).

\begin{lemma}
If $w$ is the diagonal matrix

\[\left( \begin{array}{ccccc}
w_1 & 0 & 0 & 0 & 0 \\
0 & w_2 & 0 & 0 & 0\\
0 & 0 & w_3 & 0 & 0 \\
\vdots & \vdots & \vdots & \vdots & \vdots \\
0 & 0 & 0 & 0 & w_{\ell}
\end{array} \right) \]

then $$\Delta^0(w, \Delta^+) \rho(w) = \frac{\prod_{i < j} (w_i - w_j)}{(w_1 w_2 ... w_{\ell})^{\frac{\ell-1}{2}}}$$
\end{lemma}

\proof
This is elementary.
\qed

\begin{corollary}
If $w \in E^*$, then $$\Delta^0(w, \Delta^+) \rho(w) = (-1)^{\sum_{k=1}^{\frac{\ell-1}{2}} k} \ \frac{N_{EL/L}((w-\upsilon(w))(w - \upsilon^2(w))(w - \upsilon^3(w))...(w - \upsilon^{\frac{\ell-1}{2}}(w)))}{N_{E/F}(w)^{\frac{\ell-1}{2}}}$$ where $L$ is the unramified extension of $F$ of degree $\ell - 1$, and $\upsilon$ is an embedding of $E$ into $\overline{F}$.
\end{corollary}

\proof
This is elementary.
\qed

Recall that our proposed character formula reduces in the case $\delta_{E/F} \neq 1$ to $$F(\tilde\chi)(w) = \epsilon(\tilde\chi, \Delta^+, \tau) \frac{\displaystyle\sum_{s \in W} \chi({}^s w) \Delta_{\chi}({}^s w)}{\tau(\Delta^0(w, \Delta^+)) \rho_{\tau}(w)}, \ \ w \in T(F)^{reg}$$

\begin{theorem}\label{charactersmatchingupell1}
$F(\tilde\chi)$ agrees with the character of the supercuspidal representation $\pi_{\chi \Delta_{\chi}}$ on the range $\{ w \in E^* : 0 \leq n(w) \leq r/2 \}$.
\end{theorem}

\proof
Since $\tau_o$ is trivial on $N_{EL/L}((EL)^*)$, and since $N_{E/F}(w) = N_{EL/L}(w)$, we have that

\noindent $\tau_o(\Delta^0(w, \Delta^+)) \rho_{\tau_o}(w)  = \tau_o((-1)^{\sum_{k=1}^{\frac{\ell-1}{2}} k})$.  Thus,  $$F(\tilde\chi)(w) = \epsilon(\tilde\chi, \Delta^+, \tau)' \frac{\displaystyle\sum_{s \in W} \chi({}^s w) \Delta_{\chi}({}^s w)}{|\Delta^0(w, \Delta^+) \rho(w)|} = $$ $$\epsilon(\tilde\chi, \Delta^+, \tau)' \frac{\displaystyle\sum_{s \in W} \chi({}^s w) \Delta_{\chi}({}^s w)}{|D(w)^{1/2}|}$$ since recall that $|\Delta^0(w, \Delta^+) \rho(w)| = |D(w)|^{1/2}$ from Chapter 3.
\qed

\noindent Note that $F(\tilde\chi)$ is independent of the choice of $\tau$.  That is, all that matters is $\tau_o|_{L^*}$, which we have required from the outset is $\aleph_{EL/L}$.

Note that in the above, we have chosen $\Delta^+$ to be the standard set of positive roots, which implies that $\epsilon(\Delta^+) = 1$.  We wish to make the following observation.  Suppose we made another choice of positive roots.  Any other choice is of the form $s \Delta^+$ where $\Delta^+$ is the standard choice of positive roots and $s \in W(G(\overline{F}),T(\overline{F}))$.  Let $\rho$ be half the sum of positive roots in $\Delta^+$ and let $\rho_s$ be half the sum of positive roots in $s \Delta^+$.  Then $\Delta^0(w, s \Delta^+) \rho_s(w) = (-1)^{\ell(s)} \Delta^0(w, \Delta^+) \rho(w)$ where $\ell(s)$ is the length of $s$.  Therefore, the denominator in our character formula for the choice $s \Delta^+$ would include the term $\tau_o((-1)^{\ell(s)})$.  However, because our definition of $\epsilon(\tilde\chi, \Delta^+, \tau)$ includes the term $\epsilon(\Delta^+)$, our overall character formula $F(\tilde\chi)$ remains the same regardless of the choice of positive roots.  The same line of reasoning is true for the case of $\delta_{E/F} = 1$ and $PGL(\ell,F)$.

\subsection{The case $\delta_{E/F} = 1$}

Let $E/F$ now be a degree $\ell$ extension such that $\delta_{E/F} = 1$.  This occurs if and only if $\Delta_{\chi} = 1$.  Recall that our proposed character formula reduces in the case $\delta_{E/F} = 1$ to $$F(\tilde\chi)(w) = \epsilon(\tilde\chi, \Delta^+, \tau) \frac{\displaystyle\sum_{s \in W} \chi({}^s w) }{\tau(\Delta^0(w, \Delta^+)) \rho_{\tau}(w)}, \ \ w \in T(F)^{reg}$$

\begin{theorem}\label{charactersmatchingupell2}
$F(\tilde\chi)$ agrees with the character of the supercuspidal representation $\pi_{\chi \Delta_{\chi}}$ on the range $\{ w \in E^* : 0 \leq n(w) \leq r/2 \}$.
\end{theorem}

\proof Since $\tau_o$ is trivial on $N_{EL/L}((EL)^*)$, and since $N_{E/F}(w) = N_{EL/L}(w)$, we have that

\noindent $\tau_o(\Delta^0(w, \Delta^+)) \rho_{\tau_o}(w) = \tau_o((-1)^{\sum_{k=1}^{\frac{\ell-1}{2}} k})$.  Thus, our character formula reduces to  $$F(\tilde\chi)(w) = \epsilon(\tilde\chi, \Delta^+, \tau)' \frac{\displaystyle\sum_{s \in W} \chi({}^s w)}{|\Delta^0(w, \Delta^+) \rho(w)|} = $$ $$\epsilon(\tilde\chi, \Delta^+, \tau)' \frac{\displaystyle\sum_{s \in W} \chi({}^s w)}{|D(w)^{1/2}|},  \ \ w \in E^*$$
\qed

Note again that $F(\tilde\chi)$ is independent of the choice of $\tau$.  That is, all that matters is $\tau_o|_{L^*}$, which we have required from the outset is $\aleph_{EL/L}$.  Also note again that we have chosen $\Delta^+$ to be the standard set of positive roots, which implies that $\epsilon(\Delta^+) = 1$.  Moreover, $F(\tilde\chi)$ is independent of the choice of $\Delta^+$ for the same reasoning as in the previous section.

\subsection{On whether there are two positive depth character formulas coming from the same Cartan}

In the next two sections we show that a positive depth supercuspidal representation of $GL(\ell,F)$ is uniquely determined by the restriction of its distribution character to the $n(w) = 0$ range.  In this section, we show that if the distribution characters of two positive depth supercuspidal representations, both coming from the same Cartan, agree on the $n(w) = 0$ range, then they are isomorphic.

\begin{theorem}\label{samecartansell}
Suppose $(E/F, \chi_1)$ and $(E/F, \chi_2)$ are admissible pairs such that $F(\tilde\chi_1)(w) = F(\tilde\chi_2)(w) \ \forall w \in E^* : n(w) = 0$.  Then, $\chi_1 = \chi_2^{\upsilon}$ for some $\upsilon \in Aut(E/F)$.
\end{theorem}

We will split the proof of this theorem into several cases.  Note that as in Section \ref{othercartan}, it is sufficient in this and the next section to consider admissible pairs rather than regular pairs.

\begin{proposition}\label{ramifiedgalois}
Let $E/F$ be ramified Galois.  Suppose $(E/F, \chi_1), (E/F, \chi_2)$ are admissible pairs such that $F(\tilde\chi_1)(w) = F(\tilde\chi_2)(w) \ \forall w \in E^* : n(w) = 0$.  Then $\chi_1 = \chi_2^{\upsilon}$ for some $\upsilon \in Aut(E/F)$.
\end{proposition}

\proof
The proof from Lemma (\ref{samecartanramifiedgl2}) can be adapted to this setting.  In order to use this proof, we need to prove that there exists a $w' \in E^* \setminus F^*(1 + \mathfrak{p}_E)$ such that $[\chi_1](w') \neq 0$.  We will do this in the next section.
\qed

\begin{proposition}\label{ramifiednongalois}
Let $E/F$ be ramified non-Galois.  Suppose $(E/F, \chi_1), (E/F, \chi_2)$ are admissible pairs such that $F(\tilde\chi_1)(w) = F(\tilde\chi_2)(w) \ \forall w \in E^* : n(w) = 0$.  Then $\chi_1 = \chi_2$.
\end{proposition}

\proof
Assume $F(\tilde\chi_1)(w) = F(\tilde\chi_2)(w)$ on the set $\{w \in E^* : n(w) = 0 \}$.  Then since $Aut(E/F) = 1$, by inspecting the formula for $F(\tilde\chi_1), F(\tilde\chi_2)$, we see that there are constants $c_1, c_2$ such that $c_1 \chi_1(w) = c_2 \chi_2(w) \ \forall w : n(w) = 0$.  By \cite[Lemma 5.1]{spice}, this says that $\chi_1(w) = \chi_2(w) \ \forall w \in F^*(1 + \mathfrak{p}_E)$.  We may now proceed as in the proof of Lemma (\ref{samecartanramifiedgl2}), but adapted to this setting.
\qed

\begin{proposition}\label{unramified}
Let $E/F$ be unramified.  Suppose $(E/F, \chi_1), (E/F, \chi_2)$ are admissible pairs such that $F(\tilde\chi_1)(w) = F(\tilde\chi_2)(w) \ \forall w \in E^* : n(w) = 0$.  Then $\chi_1 = \chi_2^{\upsilon}$ for some $\upsilon \in Aut(E/F)$.
\end{proposition}

\proof
See \cite[page 16]{spice}).
\qed

\subsection{On whether there are two positive depth character formulas coming from different Cartans}\label{othercartanell}

In this section we show that the distribution characters of two positive depth supercuspidal representations, coming from different Cartans, can't agree on the $n(w) = 0$ range.  This, together with the results from the previous section, shows that if $(E/F, \chi)$ is an admissible pair, then there is a unique positive depth supercuspidal representation whose distribution character agrees with $F(\tilde\chi)$ on the range $\{w \in E^* : n(w) = 0 \}$.

\begin{theorem}\label{differentcartansell}
Suppose $(E/F, \chi)$ and $(E_1/F, \chi_1)$ are admissible pairs with $E \ncong E_1$.  Then \\ $\exists w \in E^* : n(w) = 0$ such that $F(\tilde\chi)(w) \neq \theta_{\pi_{\chi_1 \Delta_{\chi_1}}}(w)$.
\end{theorem}

\noindent There are several cases to check, and we split them up in a sequence of propositions.

\begin{proposition}\label{case1ell}
Suppose $(E/F, \chi)$ and $(E_1/F, \chi_1)$ are admissible pairs with $E$ ramified Galois and $E_1$ unramified.  Then $\exists w \in E^* : n(w) = 0$ such that $F(\tilde\chi)(w) \neq \theta_{\pi_{\chi_1 \Delta_{\chi_1}}}(w)$.
\end{proposition}

\proof
The following proof is due to Loren Spice.  It is vastly shorter than our original proof.  It is shown in \cite{takahashi1} and \cite{takahashi2} that $\theta_{\pi_{\chi_1 \Delta_{\chi_1}}}(w) = 0 \ \forall w \in E^* : n(w) = 0$.  Thus, if we can find a single element of $\{ w \in E^* : n(w) = 0 \}$ such that $F(\tilde\chi)(w) \neq 0$, then we'd be done.  Suppose by way of contradiction that $F(\tilde\chi)(w) = 0 \ \forall w \in E^* : n(w) = 0$.  Recall that $\{w \in E^* : n(w) = 0 \} = E^* \setminus F^*(1 + \mathfrak{p}_E)$.  Fix $w_1 \in E^* \setminus F^*(1 + \mathfrak{p}_E)$.  Note that if $\alpha \in F^*(1 + \mathfrak{p}_E)$, then $w_1 \alpha \in E^* \setminus F^*(1 + \mathfrak{p}_E)$.  Since $F(\tilde\chi)(w) = 0 \ \forall w \in E^* : n(w) = 0$, then by considering the numerator of $F(\tilde\chi)$, we get that

\[
\Bigl(\displaystyle\sum_{i = 0}^{\ell - 1} c_i\chi^{\upsilon^i}\Bigr)(\alpha)
= \displaystyle\sum_{i = 0}^{\ell - 1} \chi^{\upsilon^i}(w_1 \alpha)
= 0
\]

\noindent for all $\alpha \in F^\times(1 + P_E)$, where $c_i = \chi^{\upsilon^i}(w_1) \ne 0$ for $i = 0, \ldots, \ell - 1$.  Then linear independence of characters gives that $\chi^{\upsilon^i} = \chi^{\upsilon^j}$ on $F^\times(1 + P_E)$ for all $i$ and $j$ and this contradicts admissibility of the pair $(E/F, \chi)$.
\qed

\begin{proposition}\label{case1'ell}
Suppose that either

(1) $E/F$ is ramified non-Galois and $E_1/F$ is unramified,

(2) $E/F$ is unramified and $E_1/F$ is ramified,

(3) $E/F$ is ramified Galois and $E_1/F$ is ramified Galois such that $E \ncong E_1$, or

(4) $E/F$ is ramified non-Galois and $E_1/F$ is ramified non-Galois such that $E \ncong E_1$.

\noindent Suppose $(E/F, \chi)$ and $(E_1/F, \chi_1)$ are admissible pairs.  Then $\exists w \in E^* : n(w) = 0$ such that $F(\tilde\chi)(w) \neq \theta_{\pi_{\chi_1 \Delta_{\chi_1}}}(w)$.
\end{proposition}

\proof
It is shown in \cite{takahashi1} and \cite{takahashi2} that $\theta_{\pi_{\chi_1 \Delta_{\chi_1}}}(w) = 0 \ \forall w \in E^* : n(w) = 0$.  Thus, if we can find a single element of $\{ w \in E^* : n(w) = 0 \}$ such that $F(\tilde\chi)(w) \neq 0$, then we'd be done.  For case (1) this is clear, because the numerator of $F(\tilde\chi)$ is $\chi(w) \Delta_{\chi}(w)$, which takes values in $\mathbb{C}^*$.  For case (2), see \cite[pages 11,12,16]{spice}.  For case (3), the same proof as in Proposition (\ref{case1ell}) works here.  For case (4), the same argument as in case (1) works here.
\qed

Note that we have now checked all cases, because if $F$ is a local field of characteristic zero, then $F$ can't simultaneously have a degree $\ell$ ramified non-Galois extension and a degree $\ell$ ramified Galois extension.
Therefore, we have finished the proof of Theorem (\ref{differentcartansell}).  Summing up, we have altogether shown that if $(E/F, \chi)$ is a regular pair such that $\chi$ has positive level, then there is a unique positive depth supercuspidal representation, $\pi_{\chi \Delta_{\chi}}$, whose character, on the range $\{z \in T(F)^{reg} : 0 \leq n(z) \leq r/2 \}$, agrees with $F(\tilde\chi)$.  There is one minor point here to resolve.  Is there possibly a depth zero supercuspidal representation whose character, on the range $\{z \in T(F)^{reg} : 0 \leq n(z) \leq r/2 \}$, also equals $F(\tilde\chi)$?  We will prove in the next chapter that if $(E_1/F, \chi_1)$ is a regular pair corresponding to a depth zero supercuspidal representation $\pi$ via Theorem (\ref{bijection2}), then its character formula, on the range $\{z \in T(F)^{reg} : 0 \leq n(z) \leq r/2 \}$, is $$F(\tilde\chi_1)(w) = (-1)^{\ell+1} \frac{deg(\pi)}{ deg(\sigma)} \tau_o((-1)^{\sum_{k=1}^{\frac{\ell-1}{2}} k}) \left(\frac{\displaystyle\sum_{i=0}^{\ell-1} \chi_1(\upsilon^i(w))}{ \tau(\Delta^0(w,\Delta^+)) \rho_{\tau}(w)}\right), \ w \in E_1^* \setminus F^*(1 + \mathfrak{p}_{E_1})$$ where $\upsilon$ is a generator of $Aut(E_1/F)$.  Then, the same arguments as in Theorems (\ref{samecartansell}) and (\ref{differentcartansell}) show that the character of $\pi$ cannot equal $F(\tilde\chi)$, on the range $\{z \in T(F)^{reg} : 0 \leq n(z) \leq r/2 \}$, unless $\pi \cong \pi(\tilde\chi)$.

Therefore, combining Theorems (\ref{samecartansell}), (\ref{differentcartansell}), and (\ref{charactersmatchingupell1}) and (\ref{charactersmatchingupell2}), we obtain the following result.

\begin{theorem}
The assignment

\begin{eqnarray}
\left\{
\begin{array}{rl}
 irreducible \ \varphi : W_F \rightarrow GL(\ell,\mathbb{C})
\end{array} \right\} & \mapsto & \tilde\chi \in \widehat{T(F)}_{\tau \circ \rho} \mapsto \pi(\tilde\chi) \nonumber
\end{eqnarray}

\noindent from Section (\ref{preliminaries}) is the Local Langlands correspondence for positive depth supercuspidal representations of $GL(\ell,F)$, where $\pi(\tilde\chi)$ is the unique supercuspidal representation whose character, on the range $\{ z \in T(F)^{reg} : 0 \leq n(z) \leq r/2 \}$, is $F(\tilde\chi)$.
\end{theorem}

\section{Depth zero supercuspidal character formula for $PGL(\ell,F)$ and $GL(\ell,F))$, $\ell$ an odd prime}

\subsection{On the proof that our conjectural character formulas agree with depth zero supercuspidal characters}

In the following two sections, we prove Theorems (\ref{gl2theorem1}) and (\ref{gl2theorem2}) for the case of depth zero supercuspidal representations of $GL(\ell,F)$, where $\ell$ is an odd prime.  Since regular pairs in the setting of depth zero supercuspidal representations of $GL(\ell,F)$ are of the form $(E/F, \chi)$ where in particular $E/F$ is unramified, we have that $\delta_{E/F} = 1$.  Let us recall from the previous chapter that the proposed character formula simplifies to $$F(\tilde\chi)(w) = \epsilon(\tilde\chi, \Delta^+, \tau) \frac{\displaystyle\sum_{s \in Aut(E/F)} \chi({}^s w)}{\tau(\Delta^0(w, \Delta^+)) \rho_{\tau}(w)}, \ \ w \in T(F)^{reg}$$
For depth zero representations, we define $\epsilon(\tilde\chi, \Delta^+, \tau) := (-1)^{\ell+1} \frac{deg(\pi)}{ deg(\sigma)} \tau_o((-1)^{\sum_{k=1}^{\frac{\ell-1}{2}} k}) \epsilon(\Delta^+)$, where $\epsilon(\Delta^+)$ is as in Section (\ref{epsilonforell}) and $deg(\pi), deg(\sigma)$ are as in Theorem (\ref{depthzerocharacters}).

We will again need Theorem (\ref{depthzerocharacters}).  Let $(E/F, \chi)$ be a regular pair such that $\chi$ has level zero and $E/F$ is degree $\ell$. Thus, $E/F$ is unramified and $\chi|_{U_E}$ gives rise to a character $\theta$ of the multiplicative group of the residue field $\mathbb{F}_{q^{\ell}}$ of $E$.  Note that when $E/F$ is unramified, $(E/F, \chi)$ is regular if and only if $(E/F, \chi)$ is admissible.  Let $\mathbb{G} := GL(\ell, \overline{\mathbb{F}_q})$.  Let $\mathbb{T}$ be the maximal torus of $\mathbb{G}$ defined over $\mathbb{F}_q$ such that $\mathbb{T}^{\Phi} = \mathbb{F}_{q^{\ell}}^*$ is the elliptic torus in $GL(\ell, \mathbb{F}_q)$.  Then, by Deligne-Lusztig theory, the pair $(\mathbb{T}, \theta)$ yields a generalized character $R_{\mathbb{T}, \theta}$ of $\mathbb{G}(\mathbb{F}_q) = GL(\ell, \mathbb{F}_q)$.

\begin{proposition}
$$R_{\mathbb{T}, \theta}(s) = (-1)^{\ell+1} \displaystyle\sum_{i=0}^{\ell-1} \theta(\upsilon^i(s))$$ for all regular semsimple s in $\mathbb{T}^{\Phi}$, where $\upsilon$ is a generator of $Gal(\mathbb{F}_{q^{\ell}}/\mathbb{F}_q)$.
\end{proposition}

\proof
Exactly as in the beginning of the proof of Proposition (\ref{gl2characterformuladepthzerolala}), we get
$$R_{\mathbb{T}, \theta}(s) = \frac{\epsilon_{\mathbb{T}} \epsilon_{C^0(s)}}{ |C^0(s)^{\Phi}|_p}  \displaystyle\sum_{i=0}^{\ell -1} \theta(\upsilon^i(s))$$ since the relative Weyl group is $W(\mathbb{G}(\mathbb{F}_q),\mathbb{T}(\mathbb{F}_q)) = Aut(\mathbb{F}_{q^{\ell}}/\mathbb{F}_q)$.  It remains to calculate the constants in front.

Now, since $s \in \mathbb{T}^{\Phi}$ is regular semisimple, then $|C^0(s)^{\Phi}|_p = 1$.  Moreover, $\epsilon_{\mathbb{T}} = -1$ and $\epsilon_{C^0(s)} = (-1)^{\ell}$.  Therefore, $$R_{\mathbb{T}, \theta}(s) = \epsilon_{\mathbb{T}} \epsilon_{C^0(s)} \displaystyle\sum_{i=0}^{\ell-1} \theta(\upsilon^i(s)) = (-1)^{\ell+1} \displaystyle\sum_{i=0}^{\ell-1} \theta(\upsilon^i(s))$$
\qed

Our character formula is defined on the unramified elliptic torus $E^*$.  We wish to show that our character formula agrees with a depth zero supercuspidal character on the sets where they are both defined i.e.  the set $(F^* K_0 \setminus F^* K_1) \bigcap E^*$.

\begin{lemma}
$(F^* K_0 \setminus F^* K_1) \bigcap E^* = E^* \setminus F^*(1 + \mathfrak{p}_E)$
\end{lemma}

\proof
The proof is not difficult.
\qed

\begin{theorem}\label{depthzerocharactersmatchingupell}
$F(\tilde\chi)$ agrees with the character of the depth zero supercuspidal representation $\pi_{\chi \Delta_{\chi}}$ on the range $E^* \setminus F^* ( 1 + \mathfrak{p}_E) = \{z \in T(F)^{reg} : 0 \leq n(z) \leq r/2 \}$.
\end{theorem}

\proof
Recall that since $E/F$ is unramified, $\Delta_{\chi} \equiv 1$.  Therefore, we need to show that

$$\tau_o((-1)^{\sum_{k=1}^{\frac{\ell-1}{2}} k}) \left(\frac{\displaystyle\sum_{i=0}^{\ell-1} \chi(\upsilon^i(w))}{  \tau(\Delta^0(w,\Delta^+)) \rho_{\tau}(w)}\right) = \displaystyle\sum_{i=0}^{\ell-1} \chi(\upsilon^i(w)) \ \ \ \forall w \in E^* \setminus F^*(1 + \mathfrak{p}_E)$$

Recall that $\tau(\Delta^0(w,\Delta^+)) \rho_{\tau}(w) = \tau_o(\Delta^0(w, \Delta^+)) \rho_{\tau_o}(w) |\Delta^0(w, \Delta^+) \rho(w)|$.  Let $w \in E^*$.
Recall that $$\Delta^0(w, \Delta^+) \rho(w) = (-1)^{\sum_{k=1}^{\frac{\ell-1}{2}} k} \frac{N_{EL/L}((w-\upsilon(w))(w - \upsilon^2(w))(w - \upsilon^3(w))...(w - \upsilon^{\frac{\ell-1}{2}}(w)))}{N_{E/F}(w)^{\frac{\ell-1}{2}}}$$

Therefore, $\tau_o(\Delta^0(w,\Delta^+)) \rho_{\tau_o}(w) = \tau_o((-1)^{\sum_{k=1}^{\frac{\ell-1}{2}} k})$ for any element $w \in E^*$.  Now, let $w \in E^* : n(w) = 0$.  We claim that $|\Delta^0(w, \Delta^+) \rho(w)| = 1$.  For let $w = p^n u$, where $u \in \mathfrak{o}_E^*$ and $n \in \mathbb{Z}$.  Then
$$|\Delta^0(w, \Delta^+) \rho(w)| = $$ $$\left|\frac{N_{EL/L}((p^n u -\upsilon(p^n u))(p^n u - \upsilon^2(p^n u))(p^n u - \upsilon^3(p^n u))...(p^n u - \upsilon^{\frac{\ell-1}{2}}(p^n u)))}{N_{E/F}(p^n u)^{\frac{\ell-1}{2}}} \right| = $$
$$\left|\frac{N_{EL/L}(p^n)^{\frac{\ell-1}{2}} N_{EL/L}(( u -\upsilon( u))( u - \upsilon^2( u))( u - \upsilon^3(u))...( u - \upsilon^{\frac{\ell-1}{2}}( u)))}{N_{EL/L}(p^n)^{\frac{\ell-1}{2}} N_{E/F}( u)^{\frac{\ell-1}{2}}} \right| = $$
$$\left|\frac{ N_{EL/L}(( u -\upsilon( u))( u - \upsilon^2( u))( u - \upsilon^3(u))...( u - \upsilon^{\frac{\ell-1}{2}}( u)))}{N_{E/F}( u)^{\frac{\ell-1}{2}}} \right| $$ Now, since $n(w) = 0$, this means that the leading coefficient of $u$ is in $k_E \setminus k_F$, where $k_E$ is the residue field of $E$ and $k_F$ is the residue field of $F$.  Therefore, $u - \upsilon^i(u) \in \mathfrak{o}_E^* \ \forall i = 1,2,... \frac{\ell-1}{2}$.  Therefore, $$\frac{ N_{EL/L}(( u -\upsilon( u))( u - \upsilon^2( u))( u - \upsilon^3(u))...( u - \upsilon^{\frac{\ell-1}{2}}( u)))}{N_{E/F}( u)^{\frac{\ell-1}{2}}} \ \in \mathfrak{o}_F^*$$ and therefore it's absolute value is 1.
\qed

The same argument as in the case of positive depth supercuspidal representations shows that our overall character formula $F(\tilde\chi)$ remains the same regardless of the choice of positive roots.

\subsection{On whether there are two character formulas coming from the same Cartan}

In this section, we show that if the distribution characters of two depth zero supercuspidal representations, both coming from the unramified Cartan, agree on the $n(w) = 0$ range, then the supercuspidal representations are isomorphic.  Note that $(F^* K_0 \setminus F^* K_1) \bigcap E^* = E^* \setminus F^*(1 + \mathfrak{p}_E) = \{w \in E^* : n(w) = 0 \}$.

\begin{theorem}\label{depthzerosamecartansell}
Suppose $(E/F, \chi_1), (E/F, \chi_2)$ are admissible pairs such that $F(\tilde\chi_1)(w) = F(\tilde\chi_2)(w)$ on the set $E^* \setminus F^*(1 + \mathfrak{p}_E) = \{w \in E^* : n(w) = 0 \}$.  Then $\chi_1 = \chi_2^{\upsilon}$ for some $\upsilon \in Aut(E/F)$.
\end{theorem}

\proof
See \cite[page 16]{spice}.
\qed

Summing up, we have altogether shown that if $(E/F, \chi)$ is a regular pair such that $\chi$ has level zero, then there is a unique depth zero supercuspidal representation, $\pi_{\chi \Delta_{\chi}}$, whose character, on the range $\{z \in T(F)^{reg} : 0 \leq n(z) \leq r/2 \} = \{z \in T(F)^{reg} : n(z) = 0 \}$, agrees with $F(\tilde\chi)$.  As in the argument at the end of Section \ref{othercartanell}, one also sees that there is no positive depth supercuspidal representation whose character, on the range $\{z \in T(F)^{reg} : n(z) = 0 \}$, also equals $F(\tilde\chi)$.
Therefore, combining Theorems (\ref{depthzerosamecartansell}) and (\ref{depthzerocharactersmatchingupell}), we obtain the following result.

\begin{theorem}
The assignment

\begin{eqnarray}
\left\{
\begin{array}{rl}
 irreducible \ \varphi : W_F \rightarrow GL(\ell,\mathbb{C})
\end{array} \right\} & \mapsto & \tilde\chi \in \widehat{T(F)}_{\tau \circ \rho} \mapsto \pi(\tilde\chi) \nonumber
\end{eqnarray}

\noindent from Section (\ref{preliminaries}) is the Local Langlands correspondence for depth zero supercuspidal representations of $GL(\ell,F)$, where $\pi(\tilde\chi)$ is the unique supercuspidal representation whose character, on the range $\{ z \in T(F)^{reg} : 0 \leq n(z) \leq r/2 \} = \{z \in T(F)^{reg} : n(z) = 0 \}$, is $F(\tilde\chi)$.
\end{theorem}

\section{Appendix}

Here we list some various helpful lemmas that we needed throughout the paper.

\begin{lemma}{\label{infinitesquareroot}}
Let F be a local field whose residual characteristic is not 2.  Suppose $\lambda$ is a character of $F^*$ whose order is a power of 2. Then $\lambda|_{U_E^1} \equiv 1$.
\end{lemma}

\proof
Let $w \in U_E^1$.  Then $w$ is a square since the leading term in the power series expansion is $1$.  Moreover, one of the square roots of $w$ is in $U_E^1$.  We may then proceed inductively to conclude that $w$ is a $2^n$-th power for any $n \in \mathbb{N}$.
\qed

\begin{lemma}\label{serresquare}
Suppose the residual characteristic of the local field $F$ is not 2.  Let $x = p^n u$ be an element of $F^*$, with $n \in \mathbb{Z}, u \in U_F$.  For $x$ to be a square, it is necessary and sufficient that $n$ is even and the image $\overline{u}$ of $u$ in $\mathbb{F}_p^* = U_E / U_E^1$ is a square.
\end{lemma}

\proof
See \cite[Section 3.3]{serre}
\qed

We now collect some basic properties of the Weil index, which we need for various computations (see \cite{rao} for more details).  Throughout, $F$ denotes either a local field or a finite field of characteristic $\neq 2$.  Let $\eta$ be a nontrivial additive character of $F$.  For any $a \in F$, we write $a \eta$ for the character $a \eta : x \mapsto \eta(ax)$.

\begin{definition}\label{rao0}
Define $$\gamma_F(\eta) := \mathrm{the Weil \ index \ of \ the \ map} \ x \mapsto \eta(x^2)$$ $$\gamma_F(a, \eta) := \gamma_F(a \eta) / \gamma_F(\eta) \ \ \ a \in F^*$$
\end{definition}

\begin{lemma}\label{rao1}
$(1) \gamma_F(a c^2, \eta) = \gamma_F(a, \eta)$ and $\gamma_F(ab, \eta) \gamma_F(a, \eta)^{-1} \gamma_F(b, \eta)^{-1} = (a,b)_F$.

$(2) \gamma_F(-1, \eta) = \gamma_F(\eta)^{-2}$

$(3) \{\gamma_F(a, \eta) \}^2 = (-1, a)_F = (a,a)_F$
\end{lemma}

Let $Q$ be a nondegenerate quadratic form of degree $n$ over $F$.

\begin{definition}\label{rao2}
The Hasse invariant $h_F(Q)$ is defined as follows:  $$h_F(Q) = \gamma(\eta \circ Q) \{\gamma_F(\eta) \}^{-n} \{\gamma_F(det Q, \eta) \}^{-1}$$ Here $\gamma(\eta \circ Q)$ is the Weil index of $x \mapsto \eta(Q(x))$.
\end{definition}

\begin{lemma}\label{rao3}
$(1)$ If $n = 2$, and $Q = a_1 x_1^2 + a_2 x_2^2, a_1, a_2 \in F^*$, then $h_F(Q) = (a_1, a_2)_F$.
\end{lemma}

\end{document}